\documentclass[preprint,11pt, 3p]{elsarticle}

\journal{Stochastic Processes and their Applications
}

\usepackage[english]{babel}
\AtBeginDocument{
  \addtocontents{toc}{\relax \fontsize {10}{11}\selectfont }}

\usepackage{amsmath,amsfonts,amssymb}
\usepackage{mathrsfs}
\usepackage{amsthm}
\usepackage{enumerate}
\usepackage[utf8]{inputenc}
\usepackage{fontenc}
\usepackage{xfrac}
\usepackage{amsbsy,amscd}
\usepackage{hyperref}
\usepackage{colortbl}
\usepackage{graphics,epsfig}
\usepackage{graphicx}
\usepackage{enumitem}

\newtheorem{theo}{Theorem}[section]
\newtheorem{prop}{Proposition}[subsection]
\newtheorem{lem}[prop]{Lemma}
\newtheorem{fact}[prop]{Fact}
\newtheorem{cor}[prop]{Corollary}
\newtheorem{rem}[prop]{Remark}
\newtheorem{hypot}{Assumption} 

\newcommand{\Ex}[1]  {\exists \, #1 ,\;}
\newcommand{\frl}[1]  {\forall \, #1 ,\;}

\newcommand{\Exq}[1]  {\exists \, #1 ,\quad}
\newcommand{\frlq}[1]  {\forall \, #1 ,\quad}

\newcommand{\idc}[1]  { \mathbf{1}_{\left \lbrace #1 \right \rbrace}  } 

\newcommand{\limInf}[1]  {\underset{#1\rightarrow \infty}{\lim}}

\newcommand{\medcup}  {\mathbin{\scalebox{1.5}{\ensuremath{\cup}}}}%
\newcommand{\medcap}  {\mathbin{\scalebox{1.5}{\ensuremath{\cap}}}}%
\newcommand{\ltm}	{\mathbin{\scalebox{.8}{\ensuremath{\! \times \!}}}	}
\newcommand{\Ndiv}[1]  {\text{\large{$\,^1\!/_#1$} } }

\newcommand{\hcm}[1]  {\hspace{#1 cm}}

\newcommand{\Tsum}[1]  { \ensuremath{ \textstyle{\sum_{\{#1\} }  }\, }     }

\newcommand{\Tsup}[1]  { \ensuremath{ \textstyle{\sup_{\{#1\} }  }\, }     }

\newcounter{eq} 
\newcounter{REFeq} 
\newcommand{\EQn}[2]  {\tag{$\arabic{section}.\arabic{eq}$} \label{#1} \ensuremath{ \stepcounter{eq}\stepcounter{REFeq} } }
\newcommand{\Req}[1]  {\textup{(\ref{#1})}}

\newcounter{enm} 

\newcommand{\E}  {\ensuremath{\mathbb{E}}}
\newcommand{\R}  {\ensuremath{\mathbb{R}}}
\newcommand{\N}  {\ensuremath{\mathbb{N}}}
\newcommand{\PR}  {\ensuremath{\mathbb{P}}}

\newcommand{\M}  {\ensuremath{\mathcal{M}}}
\newcommand{\cE}  {\ensuremath{\mathcal{E}}}
\newcommand{\cD}  {\ensuremath{\mathcal{D}}}
\newcommand{\B}  {\ensuremath{\mathcal{B}}}
\newcommand{\cT}  {\ensuremath{\mathcal{T}}}
\newcommand{\cX}  {\mathcal{X}}

\newcommand{\amu} {\ensuremath{ a } }

\newcommand{\np}  {\ensuremath{ n } }

\newcommand{\tp}  {\ensuremath{ t } }

\newcommand{\tDo}  { \ensuremath{t}  }

\newcommand{\tmo}  { \ensuremath{t}  }

\newcommand{\Tex}   {\ensuremath{ T } }

\newcommand{\spr}  {\ensuremath{ s } }

\newcommand{\uV}  { u }

\newcommand{\cp}  {\ensuremath{ c } }

\newcommand{\cY}  {\ensuremath{c} }

\newcommand{\heig}  {\ensuremath{h}}

\newcommand{\ext}  { \ensuremath{ \tau_\partial   }       }

\newcommand{\iMix}  {\ensuremath{ _{M}} }

\newcommand{\iET}  { \ensuremath{_{E}}  }

\newcommand{\iSv}  { \ensuremath{_{S}}  }

\newcommand{\iPs}  { \ensuremath{_{P}}  }

\newcommand{\iXT}  { \ensuremath{_{\circ}}  }

\newcommand{\iRN}  { \ensuremath{_{R}}  }

\newcommand{\iDB}  { \ensuremath{_{D}}  }

\newcommand{\iSB}  { \ensuremath{_{S}}  }

\newcommand{\iG}  { \ensuremath{_{C}}  }

\begin{document}

\begin{frontmatter}

\title{
Unique quasi-stationary distribution,
\\ with a possibly stabilizing extinction
 }
\author{Aurélien Velleret}

\address{Institut für Mathematik, Goethe Universität, Fachbereich 12, 
	60054  Frankfurt am Main, Germany, \\
 email: velleret@math.uni-frankfurt.de}

\begin{abstract}

We establish
sufficient conditions for exponential convergence to a
unique quasi-sta\-tio\-nary distribution in the total variation norm.
These conditions also ensure 
the existence and exponential ergodicity of the
\mbox{Q-process}, the process 
conditionned upon never being absorbed.
The technique relies
 on a coupling procedure 
that is related to Harris recurrence (for Markov Chains). 
It applies to general 
continuous-time and continuous-space Markov processes.
The main novelty is 
that we modulate each coupling step 
depending both
on a final horizon of time 
(for survival)
and on the initial distribution.
By this way, we could notably include 
in the convergence  
a dependency  on the initial condition. 
As an illustration, 
we consider a continuous-time birth-death process with catastrophes 
and a diffusion process
describing a (localized) population adapting to its environment.

\end{abstract}


\end{frontmatter}


\section{Introduction}
\label{U:sec:intro}
\setcounter{eq}{0}

\subsection{Presentation}

Given a continuous-time and continuous-space
Markov process with an absorbing state, 
we are interested in this work
in the long time behavior of the process 
conditionally on not being absorbed 
(not being "extinct").

More precisely, 
our first concern 
is on the marginal ---at time $t$--- 
conditioned on not being extinct ---also at time $t$--- (the \textbf{MCNE} in short).
We wish to highlight key conditions
on the process 
such that these MCNE
converge as $t\rightarrow \infty$ 
to a unique distribution $\alpha$.  
This limiting distribution is
called the \textbf{quasi-stationary distribution} (the \textbf{QSD}) 
---cf Subsections~\ref{U:sec:abs} and \ref{U:sec:ECV}, 
or chapter 2 in \cite{coll} for more details on this notion.
The techniques we use 
allow us to establish not only 
the existence and uniqueness of the QSD, 
but also the exponential convergence in total variation norm,
cf Theorem~\ref{U:AllPho}. 

In addition, 
we deduce, 
under the same conditions,
the existence of a specific eigenfunction $\heig$ of the infinitesimal generator,
with the same eigenvalue as the QSD. 
As time goes to infinity,
the renormalizing factor at time $t$ 
behaves asymptotically 
as $\heig \, \exp[-\lambda t]$, 
cf \Req{U:EtaT} and Theorem \ref{U:EtaECV}.
This convergence motivates the name \textbf{survival capacity}
that we give to $\heig$ 
(sometimes described as the "reproductive value" in ecological models).
Again, the convergence is exponential,
but not uniform over the state space
in our case.
Moreover, we deduce 
the existence 
of the \textbf{Q-process}. 
Its marginal at time $t$ is given by 
the limit (as $T\rightarrow \infty$) 
of the marginal of the original process at time $t$
conditioned on not being extinct at time $T$,
cf Theorem \ref{U:QECV}.
Thus, it is often described 
as the process conditioned to never be absorbed. 
Finally, we deduce 
for the Q-process
the existence and uniqueness of its stationary distribution $\beta$ 
together with a property related to exponential ergodicity.
\\

To deduce these results, 
our aim is 
to combine a large degree of generality 
with conditions as easy to verify as possible. 
A specificity of our approach 
is that it allows to deduce
a coupling procedure depending on the initial condition
that ensures a contraction in total variation
towards the limiting distribution.
It is only for commodity
that we have restricted the analysis 
to cases where there is a unique QSD.
One can find in \cite{AV_GS} 
an application to group selection models
where our procedure of proof
is exploited 
to deduce the convergence to some QSD
in a specific basin of attraction.
Also, the proof can be adapted quasi verbatim 
to discrete time processes.

We exploit the idea, 
first exploited in \cite{ChQSD},
to rely on a more constructive method 
in the form of a strong regeneration condition, 
analogous to Harris' recurrence
(what we can see maybe a bit more clearly in the present work).
At the foundation of our proof
is clearly the characterization given in \cite{ChQSD}
of the uniform exponential convergence to a unique QSD.
As we can see in the applications we present 
(cf Section~\ref{U:sec:Appli})
lack of reversibility is not at all an issue for our proofs.
The hope with these techniques is also 
to include easily more complexity on the stochastic models, 
(for instance time inhomogeneity) 
while relying on the same method
with uniform in time estimates 
(cf \cite{InhomChp}, \cite{BenGenD}, \cite{DV16}). 
\\

The remainder of Section~\ref{U:sec:intro} is organized as follows.
Subsection~\ref{U:sec:not} describes our general notations~;
Subsection~\ref{U:sec:abs} presents our specific setup
of a Markov process with extinction~;
and Subsection~\ref{U:sec:PA} 
the decomposition of the state space
on which we base our assumptions.
Subsection~\ref{U:sec:as} presents
the main set of conditions which we show to be sufficient 
for the exponential convergence to the QSD.
Subsection~\ref{U:sec:ECV}
states the three main theorems of the present paper,
dealing respectively with the QSD, 
the survival capacity and the Q-process.
The conditions that we present are then certainly numerous~; 
yet we believe that they are 
quite convenient to deal with in practice, 
except maybe for $\textup{$(A3)$}$,
for which we can only give a few hints in the present work
(cf Subsection~\ref{U:sec:AF} and \ref{U:sec:Adapt}).
Other issues on the assumptions are discussed 
in Subsection~\ref{U:sec:asR}.
Subsection~\ref{U:sec:comp}
is devoted to the comparison with the literature.
We turn in Section \ref{U:sec:rAs}
to elementary properties 
that relate our assumptions.
Before we deal with the main proofs of the general theorems
in Section \ref{U:sec:Pf},
we present in Section~\ref{U:sec:Appli} two applications of these.

Theses results seem to be new, 
but concern toy-models. 
We hope that they will help the reader get insight
on our approach.
The application of our theorems 
to more significant biological models
is intended for following work
(cf. already \cite{AV_Ada}, \cite{AV_GS}, \cite{MPV21}).

\subsection{Elementary notations}
\label{U:sec:not}
In the following, 
the notation $k\ge 1$ 
has generally to be understood as
$k\in \N$ 
while $\tp \ge 0$ 
(resp. $c>0$)
should be understood as 
$\tp \in \R_+:= [0, \infty)$ 
(resp. $c\in \R_+^* $ $:= (0, \infty)$). 
In this context 
(with $m\le n$),
we denote classical sets 
of integers by:
$\quad \mathbb{Z}_+:= \left \lbrace 0,1,2...\right \rbrace$,$\; \N:= \left \lbrace 1,2, 3...\right \rbrace$,
$\; [\![m, n ]\!]:= \left \lbrace m,\, m+1, ..., n-1,\, n \right \rbrace$,
\noindent
where the notation $:=$ makes explicit that we define some notation by this equality.
For maxima and minima, we usually denote:
$s \vee t:= \max\{s, t\}$,\,  
$s \wedge t:= \min\{s, t\}.$
Accordingly, for a function $\psi$, $\psi^{\wedge}$ 
(resp. $\psi^{\vee}$) will usually be used
for a lower-bound 
(resp. for an upper-bound) 
of $\psi$.

Let $\big(\Omega; (\mathcal{F}_\tp)_{\tp\ge0}; (X_\tp)_{\tp\ge0}; (P_\tp)_{\tp\ge0}; (\PR_x)_{x\in \cX\cup {\partial}} \big)$ be a time homogeneous strong Markov
process with cadlag paths on some Polish space  \mbox{$\cX\cup \{\partial\}$} [\cite{Rog00}, Definition III.1.1], where $(\cX; \B)$ is a measurable space and \mbox{$\partial \notin \cX$}. 
We also assume that the filtration $(\mathcal{F}_\tp)_{\tp\ge0}$ is right-continuous and complete.
We recall that $\PR_x(X_0 = x) = 1$, $P_\tp$ is the transition function of the process satisfying the usual measurability assumptions and Chapman-Kolmogorov equation.
The first entry time (resp. the first exit time) of $\cD$,
for some domain $\cD \subset \cX$, 
will generally be denoted by $\tau_{\cD}$ (resp. by $T_{\cD}$).
%
While dealing with the Markov property
between different stopping times, 
we wish to clearly indicate 
with our notation 
that we introduce a copy of $X$ 
(ie with the same semigroup $(P_\tp)$)
whose dependency upon $X$ 
is limited to its initial condition. 
This copy 
(and the associated stopping times)
is then denoted with a tilde
($\widetilde{X},\, \widetilde{\ext},\,
\widetilde{T}_\cD$ etc.). 
In the notation 
$\PR_{X_{\tau_{E}}} (t- \tau_{E} < \widetilde{\ext})$ 
for instance, 
$\tau_{E}$ and $X_{\tau_{E}}$ 
refer to the initial process $X$ 
while $\widetilde{\ext}$ refers 
to the \mbox{copy $\widetilde{X}$.}

\subsection{The stochastic process with absorption}
\label{U:sec:abs}

We consider a strong
Markov processes absorbed at $\partial$: the cemetery. 
More precisely, we assume that 
$X_s = \partial$ 
implies $X_\tp = \partial$ for all $\tp \ge s$
and that 
the extinction epoch:
$\ext:= \inf\left \lbrace \tp \ge 0;\; X_\tp = \partial \right \rbrace\quad$
is a stopping time. 
Thus, the family $(P_\tp)_{\tp \ge0}$
defines a non-conservative semigroup of operators on the set \mbox{$\B_+(\cX)$}
(resp.  \mbox{$\B_b(\cX)$})
of positive (resp. bounded)
$(\cX,\B)$-measurable functions. 
For any probability measure $\mu$ on $\cX$, 
that is $\mu \in \M_1\left(\cX\right)$, and \mbox{$f\in \B_+(\cX)$} 
(or \mbox{$f\in \B_b(\cX)$}) 
we use the notations:
\begin{align*}
\PR_\mu (.):= \int_{\cX} \PR_x(.) \; \mu(dx), \quad
\langle \mu \, \big| \, f\rangle:= \int_{\cX} f(x) \; \mu(dx).
\end{align*}
We denote by $\E_x$ (resp. $\E_\mu$) the expectation corresponding to $\PR_x$ (resp. $\PR_\mu$).
\begin{align*}
	&\mu P_\tp(dy):= \PR_\mu(X_\tp \in dy),\qquad
	&\langle \mu P_\tp\, \big| \, f\rangle \,&=\, \langle \mu \, \big| \, P_\tp f\rangle 
	\,= \E_\mu[f(X_\tp)],
	\\
	&\mu A_\tp(dy):= \PR_\mu(X_\tp \in dy \, \big| \, \tp < \ext), \qquad 
	&\langle\mu A_\tp\, \big| \, f\rangle \,&= \E_\mu[f(X_\tp)\, \big| \, \tp< \ext],
\end{align*}
$\mu A_\tp$ is what we called the MCNE 
(at time $\tp$, with initial distribution $\mu$).
In this setting, 
the family $(P_\tp)_{\tp\ge0}$ 
(resp. $(A_\tp)_{\tp\ge0}$)
defines a linear but non-conservative semigroup
(resp. a conservative but non-linear semigroup)
of operators on $\M_1\left(\cX\right)$ 
endowed with the total variation norm:
$
\|\mu \|_{TV}:= \sup\left \lbrace |\mu(A)|;\; A \in \B \right \rbrace 
\text{ for } \mu \in \M(\cX).$
A probability measure $\alpha$ is said to be the \textit{quasi-limiting distribution} of an initial condition $\mu$ if:
\begin{align*}
\frl{B \in \B}\quad
\limInf{\tp} \PR_\mu(X_\tp \in B\, \big| \, \tp< \ext):= \limInf{\tp} \mu A_\tp(B) = \alpha(B).
\end{align*}
It is now classical 
(cf e.g. Proposition~1 in \cite{MV12}) 
that $\alpha$ is then a quasi-stationary distribution or QSD, in the sense that:
$\frl{\tp\ge 0}\quad
\alpha A_\tp(dy)
= \alpha(dy). $

Our first purpose will be to prove that 
the assumptions in Subsection~\ref{U:sec:as}
provide sufficient conditions
for the existence of a unique quasi-limiting distribution
$\alpha$, independent of the initial condition.

\subsection{Specification on the state space}
\label{U:sec:PA}

In the following, we will always assume
the following decomposition of $\cX$:
\setcounter{hypot}{-1} 
	\begin{hypot}: 
		\textbf{"Exhaustion of $\cX$"}
		\label{U:as:PA}
		There exists a sequence $(\cD_\ell)_{\ell\ge 1}$ of closed subsets of $\cX$ such that:
\begin{align*}
			\frl{n\ge 1} \cD_\ell \subset \cD_{\ell+1} ^\circ\quad\text{ and } \quad 
			\medcup_{\ell\ge 1} \cD_\ell= \cX.
			\tag*{\textup{$(A0)$}}
		\end{align*}
	\end{hypot}

This sequence will serve as a reference for the following statements.
For instance, we will have control on the process through the fact that the initial distribution belongs to some set of the form:
\newcounter{mnx}\newcounter{Dp} \newcounter{Dps}
\begin{equation*} 
	\M_{\ell,\, \xi}  := \left \lbrace \mu \in \M_1\left(\cX\right);\; 
	\mu\left( \cD_\ell \right) \ge \xi\right \rbrace,
	\hspace{1cm} \text{ with } \xi \in (0, 1).
	\EQn{U:mnx}{\M_{\ell,\, \xi}  }
 \end{equation*}
Note that for any $\xi>0$: $\M_1\left(\cX\right) = \medcup_{\ell\ge 1} \M_{\ell,\, \xi}  .$ Let also: 
\begin{equation*} 
	\mathbf{D}:= \left \lbrace \cD;\; \cD \text{ is closed and there exists } 
	\ell\ge 1 \text{ such that } \cD \subset \cD_\ell\right \rbrace.
	\EQn{U:Dps}{\mathbf{D}}
\end{equation*}

\section{Exponential convergence to the $QSD$}
\label{U:sec:EC}
\setcounter{eq}{0}

\subsection{Hypotheses}
\label{U:sec:as}

We recall that for any set $\cD$,
we defined the first exit and entry times as:
$$T_{\cD}:= \inf\left \lbrace  \tp \ge 0;\; X_\tp \notin \cD \right \rbrace
,\quad
\tau_\cD:= \inf\left \lbrace \tp \ge 0;\; X_\tp \in \cD \right \rbrace.$$

\begin{hypot}: \textbf{"Mixing property"}
	\label{U:as:Mix}\\
	There exists some probability measure $\zeta \in  \M_1\left(\cX\right)$ such that, 
	for any $\ell \ge 1$, 
	there exists $L \ge \ell$,
	$\cp, \tp>0$ such that:
	\begin{align*}
		\tag*{\textup{$(A1)$}}
		&
		\frl{x \in \cD_{\ell}}
		\hspace{.5cm}
		\PR_x \left[ {X}_{\tp}\in dx;\;
		\tp < \ext \wedge T_{\cD_L} \right] 
		\ge \cp\; \zeta(dx).
	\end{align*}
\end{hypot}

\begin{hypot}: \textbf{"Escape from the Transitory domain"}
	\label{U:as:eT}\\
	For given $\rho >0$ and $E \in \mathbf{D}$:
	\begin{equation}
		e_\cT:=
		\underset{x\in \cX}{\sup} \;\E_{x} \left( 
		\exp\left[\rho\, (\ext\wedge \tau_{E}) \right] \right) < \infty.
		\tag*{\textup{$(A2)$}}
	\end{equation}
\end{hypot}

The order  $\rho$ in the previous exponential moment 
is required to be larger 
than the following \textbf{"survival estimate"} 
that involves the measure $\zeta$ in $\textup{$(A1)$}$: 
\begin{align*}
	\rho\iSv
	:= \sup\big\{\rho \in \R;\;
	\sup_{\ell \ge 1} \inf_{t>0} \;
	e^{\rho t}\,\PR_\zeta(t < \ext\wedge T_{\cD_\ell}) 
	= 0	\big\}.
	\EQn{U:rsv}{}
\end{align*}

	\begin{hypot}: \textbf{"Asymptotic comparison of survival"}
		\label{U:as:cd}
		
		For a given $E \in \mathbf{D}$ and $\zeta \in \M_1\left(\cX\right)$: 
\begin{align*}
			\underset{t\rightarrow \infty}{\limsup} \;
			\underset{x\in E}{\sup} \;
			\dfrac{\PR_{x} (\tp< \ext)}
			{\PR_{\zeta} (\tp < \ext)} 
			< \infty.
			\tag*{\textup{$(A3)$}}
		\end{align*} 
	\end{hypot}

We say that Assumption $\mathbf{(A)}$ holds, whenever:\\
"\textup{$(A1)$}\, holds 
for some  \mbox{$\zeta \in \M_1\left(\cX\right)$}
and a sequence $(\cD_\ell)$ that satisfies $\textup{$(A0)$}$.
Moreover, there exist 
$E \in \mathbf{D}$
such that \textup{$(A2)$}, holds
with some $\rho > \rho\iSv$
as well as \textup{$(A3)$}." 
\\

As we shall see in Subsection \ref{U:sec:rAs}, 
\textup{$(A1)$}\, implies that $\rho\iSv < \infty$.
In order to ensure Assumption $\mathbf{(A)}$,
we may not need to estimate precisely $\rho\iSv$:
it is possible (depending on the process)
that \textup{$(A2)$}\, is satisfied
for any potential value for $\rho>0$
(where $E$ is likely 
to depend on $\rho$).
Moreover, 
$\rho\iSv$ as well as  assumption \textup{$(A3)$}\,
actually do not depend of the  choice of $\zeta$ satisfying \textup{$(A1)$}.
\\

\subsection{Main Theorems: the simplest set of assumptions}
\label{U:sec:ECV}

\begin{theo}
	\label{U:AllPho}

	Assume that Assumption $\mathbf{(A)}$ 
	holds.
	Then, there exists a unique QSD $\alpha$. 
	Moreover, we have exponential convergence to $\alpha$ of the MCNE's
	at a given rate  $\gamma >0$. 
	More precisely,	
	for any pair $\ell \ge 1$ and $\xi \in (0, 1)$, 
	there exists $C = C(\ell,\, \xi)>0$ such that:
	\begin{align*}
		\frl{\tp > 0}
		\frlq{\mu \in \M_{\ell,\, \xi}  }
		\|\, \PR_\mu \left[\, X_{\tp} \in dx \; 
		| \; \tp < \ext \right]  - \alpha(dx) \, \|_{TV}
		\le C \; e^{-\gamma \; \tp}.
		\EQn{U:ECvAl}{\textup{ECv:}\,\alpha}
	\end{align*}	
\end{theo}
\noindent
It is classical (cf e.g. 
Theorem~2.2 in \cite{coll})
that, as a QSD, $\alpha$ is associated to some extinction rate $\lambda$:  
\begin{align*}
	\frlq{\tp\ge 0} 
	&\PR_\alpha (\tp< \ext) = e^{-\lambda\, \tp},
	\text{ so that } 
	\alpha P_\tp = e^{-\lambda\, \tp}\, \alpha.
	\EQn{U:LBz}{\lambda}
	\\\text{Let: }
	\quad
	&\heig_\tp (x):= e^{\lambda\, \tp}  \PR_x (\tp< \ext).
	\EQn{U:EtaT}{\heig_\tp}
\end{align*}

\begin{theo}
	\label{U:EtaECV}
	Again under Assumption $\mathbf{(A)}$, 
	we have exponential convergence 
	in the supremum norm
	of $(\heig_\tp)_{\tp\ge 0}$ to 
	a limit $\heig$,
	with the rate $\gamma$ 
	deduced from \Req{U:ECvAl}.
	The function $\heig$, which describes the "survival capacity" of the initial condition $\mu$, 
	has a  positive lower-bound on any $\cD_\ell$, 
	an upper-bound on $\cX$
	and vanishes on $\partial$. 
	It also belongs to the domain 
	of the infinitesimal generator $\mathcal{L}$, 
	associated with the semi-group 
	$(P_\tp)_{\tp \ge0}$ 
	on $(B(\cX\cup\{\partial\} ); \| . \|_\infty)$, 
	and:
	\begin{align*}
		\mathcal{L}\, \heig = -\lambda\, \heig, \qquad 
		so\quad \frl{\tp \ge 0}
		P_\tp \, \heig = e^{-\lambda\, \tp }  \heig.
		\EQn{U:EFeta}{}
	\end{align*}
\end{theo}

\noindent \textsl{Remark:} 
Like in \cite{ChQSD}, 
it is also not difficult to show that there is no eigenvalue of $\mathcal{L}$ 
between $0$ and $-\lambda$, 
and that $\heig$ is the unique eigenvector
associated to $-\lambda$.

\begin{theo}
	\label{U:QECV}
	Under again Assumption $\mathbf{(A)}$,
	we have:

\noindent	\textup{(i) \textbf{Existence of the $Q$-process:}} \\
	There exists a family $(\mathbb{Q}_x)_{x \in \cX}$ of probability
	measures on $\Omega$ defined by:
	\begin{align*}
		\limInf{\tp }
		\PR_x(\Lambda_\spr \, \big| \, \tp  < \ext) = \mathbb{Q}_x(\Lambda_\spr),
		\EQn{Qdef}{\mathbb{Q}_x}
	\end{align*}
	for all $\mathcal{F}_\spr$-measurable set $\Lambda_\spr$. 
	The process $(\Omega;(\mathcal{F}_\tp )_{\tp \ge 0};
	(X_\tp )_{\tp \ge 0};(\mathbb{Q}_x)_{x \in \cX})$
	is an $\cX$-valued homogeneous strong Markov process.
\end{theo}

\noindent\textbf{(ii) Transition kernel:}\\
\textit{The transition kernel of the Markov process $X$ under
	\mbox{$(Q_x)_{x \in \cX}$} is given by:\newcounter{qt}
	\begin{align*}
		&q(x; \tp ; dy) = e^{\lambda\, \tp } \, \dfrac{\heig(y)}
		{\heig(x)}\; p(x; \tp ; dy),
		\EQn{U:qt}{q(x; \tp ; dy)}
	\end{align*}
	where $p(x; \tp ; dy)$ is the transition kernel of the Markov process $X$ under
	\mbox{$(P_x)_{x \in \cX}$}.}
\textit{In other words, for all $\psi \in \B_b(\cX)$ and $\tp  \ge 0$,
	$\quad \langle\delta_x\, Q_\tp  \, \big| \, \psi \rangle 
	= e^{\lambda\, \tp }\,
	\langle\delta_x \,P_\tp  \, \big| \, \heig\times\psi\rangle\, / \,\heig(x), \quad$
	where $(Q_\tp )_{\tp \ge 0}$ is the semi-group of $X$ under $\mathbb{Q}$.}
\\

\noindent \textit{\textbf{(iii) Exponential ergodicity:} \\
	There is a unique invariant distribution 
	of $X$ under $\mathbb{Q}$, 
	defined by:
	$$\beta(dx):= \heig(x)\, \alpha(dx).$$
	Moreover, 
	there exists $\gamma >0$ 
	and $C = C(\ell, \xi)$
	such that:}
%
\begin{align*}
	&\quad\frl{\tp  > 0}
	\frlq{\mu \in \M_{\ell,\, \xi}  }
	\| \mathbb{Q}_{\mu B[\heig]} ( X_{\tp } \in dx)
	- \beta(dx) \|_{\frac{1}{\heig}}
	\le C\; e^{-\gamma \; \tp }.
	\EQn{U:ECvBeta}{\textup{ECv:}\,\beta}
	\\&\quad \text{ where } 
	\|\mu\|_{\frac{1}{\heig}} 
	:= \left \|\dfrac{\mu(dx)}{\heig(x)}\right \|_{TV}
	\ge \dfrac{\|\mu(dx)\|_{TV}}{\|  \heig \|_\infty},
	\quad
	\mu  B[\heig](dx) 
	:= \heig(x)\, \mu(dx)\, 
	/ \langle\mu\, \big| \, \heig\rangle,
	\\&\hcm{1.5}
	\mathbb{Q}_\mu (dw):= \textstyle{\int_\cX }\mu(dx) \, \mathbb{Q}_x (dw).
\end{align*}

\subsection{How to verify \textup{$(A3)$}?}
\label{U:sec:AF}

For discrete space, it is quite natural to deduce \textup{$(A3)$}\, 
from the fact that there exists $\tp$ such that:
$\inf_{x\in E} \PR_{\zeta}(X_\tp = x) > 0.$
We can thus couple
some trajectories 
starting from $\zeta$ and
passing in $x$ at time $t$
to the set of all trajectories starting from $x$.
From this 
we can infer a lower-bound 
of the asymptotic survival ability 
of the former (starting from $\zeta$)
compared to the latter 
(starting from $x$).
For an illustration,
this coupling
is exploited in the birth and death process 
in Subsection \ref{U:sec:bd}.  

For continuous space however, 
the process starting from $\zeta$ will never hit precisely $x$.
We need  to wait a bit for the process starting from $x$ to diffuse before 
the association we expect can be ensured. 
Although it appears quite more complicated, 
our argument is very similar. 
In cases where the Harnack inequality holds 
(notably pure diffusive processes,
cf Subsection \ref{U:sec:Harn}),
one is usually able to prove:\\
$$
\frl{x \in E}
\hspace{0.5cm}
\PR_{x} \big(X_\tp \in dx;\; \tp < \ext \big) 
\le \cp \,\PR_{\zeta} \big(X_{\tp_{\alpha}} \in dx;\; \tp_{\alpha} < \ext\big),\,$$
where $\tp, \tp_{\alpha}, c>0$ are independent of $x$.
Like for the discrete space case,
we then deduce \textup{$(A3)$}\, from the Markov property
and an additional control of the survival on finite time-interval.
Note that when the Harnack inequality  holds, 
it is natural to exploit it already in the proof of $(A1)$,
cf Subsection \ref{U:sec:Harn}.

In a much general setting,
and especially when jumps 
are involved in the process,
the situations might get much tricky 
when one wishes to look for a similar coupling of trajectories.
The issue is notably 
on exceptional behavior 
along which we have poor controls 
(no jump for a long time,
too many jumps, too large etc.).
In \cite{AV_disc}, 
we provide 
a very efficient and more easily verified condition 
which ensures  \textup{$(A3)$}, given the other assumptions.
Since this condition is technical and may appear too abstract 
without the illustration of various examples,
we encourage any interested reader 
to look at \cite{AV_Ada} and \cite{AV_M}, 
besides the simple illustrations 
given in \cite{AV_disc}.

\subsection{Remarks on the Assumptions and the results}
\label{U:sec:asR}

\begin{rem}
	Since $X$ is right-continuous 
	and the filtration is both right-continuous 
	and complete,
	the first entry time of any Borel set 
	is a stopping time,
	cf. Theorem~52 in \textup{\cite{Mey}},
	or more recently Theorem~2.4 in \textup{\cite{B10}}.
	It means in particular that
	the first exit time $T_{\cD_\ell}$ 
	and the first entry time $\tau_{\cD_\ell}$ 
	are stopping times
	(for any $\ell \ge 1$ 
	and any initial condition).
	The result extends in fact 
	to any iterated combination 
	of the kind 
	"next entry time of $\cD_\ell$
	after the first exit time of $\cD_L$
	following the first entry time of $\cD_\ell$".
	For this, 
	we shall use that 
	there is a positive gap 
	between each of the three random times 
	(say $\tau_0< T_1<\tau_1$) involved,
	and that for any $t$, 
	$:(s, \omega)\mapsto \idc{\tau_0(\omega) < s\le t}$
	has left continuous paths
	(and similarly with $T_1$ instead of $\tau_0$
	and possibly so on by induction).
\end{rem}

\begin{rem}
	This property on first entry times
	is the main reason for us to assume $\cX$ Polish.
	The space topology is not much exploited.
	It means notably that one could treat càdlàg processes $X$
	that are known to satisfy the strong Markov property
	only w.r.t. the basic filtration 
	(that is a priori not right-continuous),
	as long as this strong Markov property can be ensured 
	for sequences of entry and exit times of sets 
	from the family $(\cD_\ell)_{\ell\ge 1}$.
	It is for simplicity that 
	we assume 
	that the strong Markov property 
	is fulfilled for a right-continuous filtration,
	noting that it holds for our examples.
	At least, as stated in Theorem 7.7 of \cite{SW13},
	a Feller semigroup on a locally compact and separable space
	generates a process that is strong Markov for the augmented filtration that we consider.
	Yet, the required condition on the family  $(\cD_\ell)_{\ell\ge 1}$
	covers a broader range of processes.
	
	For an example of càdlàg process 
	for which one may easily find suitable sets $(\cD_\ell)_{\ell\ge 1}$
	although the strong Markov property 
	is not satisfied for the augmented filtration,
	one can adapt
	the counter-example provided in page 90 of \cite{SW13}:
	such counter-examples
	are easily produced
	by specifying a partially absorbing set
	whose exit time is non-predictible.
	We can think for instance that the process stays 
	where it hits this sets 
	(in \cite{SW13}, the set is $\{0\}$ for a process on $\R$)
	for an exponential time, before leaving.
\end{rem}	
\begin{rem}	
	As we can see in the illustrations
	of \cite{AV_disc}
	(exploiting this result), 
	it is not required for the process to be strong Feller:
	for jump processes, there may exist 
	bounded measurable function $f$ such that
	$P_t f$ is discontinuous.
	$\heig$ itself might not be continuous, 
	notably for discontinuous jump rate.
\end{rem}


\begin{rem}
	\textup{$(A1)$}\, imposes
	a weak form of irreducibility condition, 
	with this reference measure $\zeta$,
	and a coherence in time to prevent periodicity.	
	
	It may happen 
	that there exists absorbing domains $\cD_A$
	(whose escape can only happen at $\ext$).
	Any MCNE with initial condition $x \in \cD_A$ 
	is necessarily supported in $\cD_A$.
	Any $\zeta$ that satisfies $\textup{$(A1)$}$
	is thus also supported in $\cD_A$.
	Moreover, 
	if these MCNE converge 
	to a unique QSD as in our result, 
	this QSD is necessarily supported 
	on $\cD_A$ as well.
\end{rem}

\begin{rem}
	Assumption \textup{$(A1)$}\, is a stronger version of Doeblin's condition that appears for the convergence of Markov Chains without extinction. 
	It also implies that any border of extinction
	shall be \mbox{approached} by the sequence $\cD_\ell$ while $\ell\rightarrow \infty$,
	but never from inside any $\cD_\ell$, 
	since by Lemma \ref{U:lem:mixbar}: 
	\begin{align*}
		\frl{\ell \ge 1}
		\frl{\tp>0}
		\quad
		\inf\{ \PR_x(\tp<\ext);\;
		x\in \cD_\ell\} \;> 0.
		\end{align*}
\end{rem}

%
\begin{rem}
	When it concerns pure jump processes, 
	one can generally choose $\cD_{L}:= \cD_{\ell}$. 
	For other processes, 
	one often needs "a bit of space" 
	between $\cD_{\ell}$ and $\cD_{L} ^c$
	to obtain a lower bound uniform in $x\in \cD_{\ell}$ 
	over trajectories from $x$ to $\zeta$ staying inside $\cD_{L}$
	(as in our second application with a diffusion).
\end{rem}

\begin{rem}
	To understand \Req{U:ECvBeta}, 
	it is worth noticing that,
	considering some general initial condition 
	in the left-hand side of \Req{Qdef},
	we obtain for the Q-process a biased initial condition:
	\begin{align*}
		&\hspace{1cm}
		\frl{\mu \in \M_1\left(\cX\right)}\quad
		\textstyle{\lim_{\tp \rightarrow \infty}}
		\PR_\mu(\Lambda_\spr \, \big| \, \tp  < \ext) = \mathbb{Q}_{\mu  B[\heig]}(\Lambda_\spr).
		\EQn{U:etaEt}{\mathbb{Q}_{\mu B[\heig]}}
		\end{align*}
	To deduce \Req{U:ECvBeta} from \Req{U:ECvAl}, 
	we reformulate \Req{U:qt}
	in terms of $B[\heig]$, $P_\tp $, $A_\tp $ and $Q_\tp$:
	\begin{align*}
		\frl{\tp \ge 0} 
		\frl{\mu \in \M_1\left(\cX\right)}\quad 
		(\mu B[\heig]) Q_\tp  
		= \left( \mu  P_\tp\right)  B[\heig] 
		= \left( \mu   A_\tp\right) B[\heig].
		\EQn{U:etaQ}{}
		\end{align*}

	Originally, 
	we intended to adapt the proof of Theorem \ref{U:AllPho}
	on the marginal at time $t$ 
	conditioned on survival at time $t+T$ 
	to deduce a control uniform in $T$. 
	This approach is effective but leads 
	to a weaker result where $\|.\|_{\frac{1}{\heig}}$
	is replaced by $\left\|  . \right\|_{TV}$.
	The convergence of the MCNE is here more informative,
	because $\heig$ is bounded.
\end{rem}
\begin{rem}[On the indices]
	Throughout the proof, 
	the constants and sets that we consider 
	will be indexed by a capital letter referring 
	to the property they are involved in.
	The indexes $S$, $E$, $M$, $C$, $A$, $R$, $D$, $P$ and $L$ 
	stand respectively for 
	"Survival" (notably in \Req{U:rsv}, \Req{U:Rg} and \Req{U:eSb}),
	"first Entry" (in \textup{$(A2)$}, where it can also refer to "Escape", and in Lemma \ref{U:T1in}),
	"Mixing" (in \textup{$(A1)$}),
	"Containment" (throughout Section \ref{U:sec:Qt}),
	"Absorption" (in \Req{U:AS}),
	"Renewal" (in \Req{U:Mrn}),
	"Doeblin" (also in \Req{U:Mrn}),
	"Persistence" (in \Req{U:eqPS})
	and finally "Last exit" (in Lemma \ref{U:TteOut}).
	$\circ$ as an index plays a similar role 
	to indicate parameters referring 
	to the core of convergence as specified 
	in Theorem \ref{U:EtRet}.	
\end{rem}

\subsection{Comparison with the literature}
\label{U:sec:comp}
\paragraph{General perspective}
Although there is already 
a vast literature on QSDs
(see notably the impressive bibliography collected by Pollett \cite{QSDbibli}), 
the approach we follow 
seems to have been explored 
only in the very recent years.
For a review on the results
that were previously obtained, 
we refer notably to
general surveys as in \cite{coll}, 
in \cite{DP13} or more specifically 
for population dynamics in 
\cite{MV12}. 
We see already in these surveys 
how essential is 
the role played by the spectral theory. 
The spectral theory is very effective 
both to relate the QSD and the survival capacity 
to the first eigenvector of a diagonalizable operator 
and to identify the convergence rate 
as the gap between the first and the second eigenvalues 
(cf e.g. \cite{CCM16}). 
The principal drawback of the spectral theory is that 
it usually relies on reversibility.
Certainly, for 1 dimensional processes, 
this condition of time-symmetry is
quite easily satisfied ; 
while, more generally,
it can be deduced from conditions easy to verify 
(detailed balance notably).
This may explain why reversibility is 
so extensively studied.
Yet, it is a very restrictive condition for higher dimensions,
as it is well explained in the appendix A of \cite{CCM17}.

Alternative methods are usually much less effective.
In \cite{cmms}, the authors prove 
the existence of the QSD 
via a Tychonov fixed point theorem. 
Another proof for the existence of the QSD is presented 
in \cite{FKMP95} for Markov Chains on $\mathbb{Z}_+$, 
based on compactness arguments and renewal techniques.
In \cite{BP10}, the authors prove,
under quite stringent conditions,
the existence and uniqueness of the QSD
and propose estimations of this QSD 
up to some computable time,
again with renewal arguments. 
The authors of \cite{DM15} 
relate the speed of convergence to QSD 
to the one of a related Doob's transform towards its stationary distribution.
Yet the conditions of the last two papers seem to apply 
essentially to discrete-space processes, 
or at least when the extinction is in some sense uniformly bounded.
The existence of the QSD and the survival capacity 
has also been related, at least for discrete time and discrete space,
to the notion of R-positivity (cf e.g. \cite{T74}, \cite{T74b} or  \cite{SV96}). 
This is especially useful when the process is easily described 
by generating functions
(in particular for Galton-Watson processes)
but seems quite an abstract criterion otherwise.
Still, 
it provides the main principle of 
focusing on the exponential rate of extinction, 
which is at the core of our study.
Our proof can reasonably be judged
as an extension of the one presented in \cite{FKM96}
with a focus on general practical assumptions 
for $R$-positivity,
noting that their analysis is restricted 
to discrete-time and discrete-space Markov processes.

\paragraph{The dependency on the initial condition}
Upper-bounds of the form:
\newcounter{Cmu}
\begin{align*}
	\|\, \PR_\mu \left[\, X_{\tp} \in dx \; 
	| \; \tp < \ext \right]  - \alpha(dx) \, \|_{TV}
	\le C(\mu) \; e^{-\gamma \; \tp}.
	\EQn{U:Cmu}{C(\mu)}
\end{align*}
assume generally 
$C(\mu) = \langle \mu\, \big| \, W \rangle$
in the case of $\alpha$ being a stationary distribution 
(i.e. for processes without extinction).
The use of such a reference function $W$
has been thoroughly studied in \cite{MTbook} 
in the case of Markov Chains, 
or in e.g. \cite{DMT95}, \cite{CGZ}, \cite{CG16} 
for continuous time processes. 
The condition on $W$ 
is what characterizes it as a Lyapunov function 
and relate a priori 
to a control of the first entry time $\tau_{E}$.
Different probabilistic bounds have generally been proposed,
although, including extinction, 
exponential moments appear compulsory 
(all the more since $\lambda$, the limiting rate of extinction, 
is not precisely known).
In a loose version,
and still for $P_t$ as a conservative semi-group,
such exponential control may take the form:
\begin{align*}
P_{t_0} W \le e^{-\rho_W t_0} W + C_W,
\quad \text{ where } \rho_W, C_W >0. 
\end{align*}
$E$ is then generally chosen as the set $\{W\le d_W\}$
for $d_W$ sufficiently large for 
\begin{align*}
P_{t_0} W \le e^{-\rho'_W t_0} W + C'_W \mathbf{1}_{E}
\end{align*}
to hold (for a smaller value $\rho'_W < \rho_W$
and some $C_W' >0$).
For $E$ to be convenient with respect to the other criteria
(mixing or comparison of survival)
and especially when $\rho_S$ is not estimated precisely,
it is usually assumed that $W$ is proper
(i.e. $W(x)$ converges to infinity as $\|x\|$ tends to infinity).

An extension of this assumption 
for the non-conservative case
has been recently proposed in \cite{BCGM19}
through their Assumption $\textbf{A}$.
Another usual version involves the infinitesimal generator $\mathcal{L}$ and assumes the following:
\begin{align*}
\mathcal{L} W \le -\rho_W  W + C_W,
	\quad \text{ where } \rho_W, C_W >0. 
\end{align*}
The inequality is stronger than the previous one for proper $W$.
A related estimate for non-conservative cases
is also proposed in Proposition 2.10 of \cite{BCGM19}.
When considering extinction, we lose also 
the property of linearity over the initial condition. 
This explains 
why upper-bounds like $\langle \mu\, \big| \, W \rangle$ are not so general
and why we focus on general initial distributions 
and not only Dirac Masses.
\\

The conclusion that we present is quite natural 
for the models we have in mind, 
where extinction plays a stabilizing role, 
preventing transient dynamics. 
In the perspective of natural selection, 
we expect to observe 
the prevalence of trajectories 
leading to and gravitating 
around some basin of attraction, 
notably compared to those dragged away 
in deadlier regions.
Although the burden of mal-adaptation  
may seem light in the short run, 
if it is too hard for the process 
to escape from less adapted areas, 
one can presume 
that the process cannot have been there for long.
In particular, 
the trajectories 
starting from favorable initial conditions 
may outcompete
what remains of the distribution, 
so that it becomes the leading part 
in the convergence to the QSD.

Other expressions of $C(\mu)$
have been presented in \cite{CVly2}, \cite{FRS19}
and \cite{BCGM19},
with similar interpre\-ta\-tion.
Note that the proofs in \cite{CVly2}
and \cite{BCGM19}
concern the convergence towards a unique Yaglom distribution,
which may not be unique as a QSD, 
for any initial conditions 
with a light enough tail.  
In \cite{CVly2},
\Req{U:Cmu}
is obtained  
with a non-linear dependency of the form 
$C(\mu) 
	= C\, \langle\mu \, \big| \, \psi_1\rangle		\,/\, 
	\langle\mu \, \big| \, \psi_2\rangle$.
As can be seen 
in our following paper\cite{AV_disc},
the dependency we introduce 
implies \Req{U:Cmu} with 
$C(\mu) = C / \langle\mu \, \big| \, \heig\rangle$,
for some $C>0$.
So the former extends our result
by including 
the more classical dependency 
through a Lyapunov function.
In  \cite{FRS19} and \cite{BCGM19}, 
the convergence is stated 
in a weighted norm
involving a weight function $W$ (resp.$V$) 
related to the previous $\psi_1$.
The dependency $C(\mu)$ 
stated in their analog of \Req{U:Cmu} 
is implicitly related to both 
$ \langle\mu \, \big| \, \heig\rangle$ 
and $\langle\mu \, \big| \, W\rangle$ 
(their function $h$ plays the same role as our).
A dependency on $\psi_1$ (or on $W$)
is neglected in our article:
\textup{$(A2)$}\, ensures in a way 
that we can find some upper-bounded  $\psi_1$ 
(we refer e.g. to Lemma 3.6 in \cite{CVly2}).   

A hint to connect the current techniques to their setting
is to adjust the probabilities of transition of the Markov process
according to such weight function $W$.
The set of positive measures $\mu$
such that $\int W(x) \mu(dx)\le 1$ 
for instance takes the place of the set of probability measures
(for more details, we refer to Subsection 2.3.2 of \cite{AV_disc}).
In general practice, 
it does not seem so clear to us
how to find such Lyapunov functions
especially when one wishes to combine simple bounds
on different parts of the space.
So  we believe that our assumption $(A2)$
 is more natural to verify in many examples
(cf e.g. our second application),
while easier to interpret.

\paragraph{The assumptions}
If one can relate our set of assumptions 
to the ones proposed in \cite{FKM96},
the similarity is clearly greater 
with \cite{CVly2}, \cite{FRS19} and \cite{BCGM19},
because of the introduction of a continuous-space setting.
It is proved in \cite{BCGM19}
that their conditions are not only sufficient, but also necessary,
and similarly in Theorem 2.3 of \cite{CV20}
with regards to \cite{CVly2}.
Similar reciprocal results are obtained in \cite{AV_disc}
for the current setting.

Due to our trajectorial approach, 
we require additional confinement properties
(with restrictions of the probabilities
upon the events $\{t<T_{\cD_L}\}$).
As explained above, Lyapunov functions 
are also not directly exploited in our approach.
The survival estimates presented in \cite{CVly2}, \cite{FRS19}  and \cite{BCGM19}
appears less natural to interpret 
than the condition on survival that we propose 
($\rho> \rho\iSv$):
 they require the introduction of a function, 
say $\psi_2$,
for which the time behavior of $\psi_2(X_t)$ 
must satisfy a certain minimization property.
Given Lemma 3.2 in \cite{CVly2},
their condition appears more general
and can certainly be convenient for specific models.
Nonetheless Theorem 2.3 in \cite{AV_disc}
proves that $\rho_S$ takes the optimal value $\lambda$
provided that the convergence results \Req{U:ECvAl}
 and of the survival capacity  hold.

\paragraph{The approach}
The techniques exploited in \cite{CVly2},
\cite{FRS19} and \cite{BCGM19}
are quite different from ours. 
In the steps of the R-theory,
the study of the $h$-transformed process
is at the core of \cite{FRS19},
with a weighted norm.
Contraction estimates under similar weighted norms
are exploited in \cite{CVly2} and \cite{BCGM19}.
Our proofs are much more constructive 
and rely on a control 
on entry times of core sets
thanks to the competition between different behaviors.
It extends to models 
where the uniqueness of the QSD 
does not hold due to transitivity conditions,
as one can observe in the applications
we have in \cite{AV_GS}.
In particular, 
our work offers a new constructive perspective 
even for the results in \cite{ChQSD} 
(cf Subsection~\ref{U:sec:coupling})
since the coupling steps 
which we introduce 
apply directly to the MCNE
(and not to their linearized versions).

\section{Several implications of \textup{$(A1)$}}
\label{U:sec:rAs}
\setcounter{eq}{0}

\begin{lem}
	\label{U:alcInd}
	Assume that \textup{$(A1)$}\, holds for two probability measures
	$\zeta^1$ and $\zeta^2$.
	Then, the associated values for $\rho\iSv$ coincide.
	Moreover, the sets $E$ for which assumption \textup{$(A3)$}\,
	holds are the same for both measures.
\end{lem}

\noindent
\textbf{"From mixing to regeneration, then survival":} 
\textup{$(A1)$}\, trivially implies,
for any $\ell$ such that $\zeta(\cD_\ell) > 0$
the following \textbf{regeneration estimate}: 

There exists $t, c>0$, $L\ge\ell$ such that, with $\cD\iSv:= \cD_\ell \subset \cD_L$:
\begin{align*}
	&\frl{x\in \cD\iSv}
	\hspace{0.5cm}
	\PR_{x} ( X_{\tp} \in \cD\iSv;\;
	\tp < \ext \wedge T_{\cD_{L}})
	\ge \cp.
	\EQn{U:Rg}{R_g}
\end{align*}

\begin{lem} 
	\label{U:rg:sv}
	Assume that  \Req{U:Rg} holds
	(for $\tp,\,\cp > 0$, $\cD\iSv \subset \cD_L$
	and $\zeta(\cD\iSv) > 0$).
	Then, $\rho\iSv \le -\frac{1}{t} \ln(c)$.
\end{lem}

In particular, 
we deduce that \textup{$(A1)$}\, implies $\rho\iSv < \infty$.

\begin{lem}
	\label{U:lem:mixbar}
	\textup{$(A1)$}\, 
	is equivalent to the apparently stronger version
	(with the same $\zeta$):\\
	For any $\ell \ge 1$,
	and $t_\veebar> 0$, 
	there exists $L\ge \ell$,
	$t\ge t_\vee$ and  $\cp>0$ such that:
	\begin{align*}
		\frl{x \in \cD_{\ell}}
		\hspace{.5cm}
		\PR_x \left[ {X}_{\tp}\in dx;\;
		\tp < \ext  \wedge T_{\cD_{L}} \right] 
		\ge \cp\, \zeta(dx).
		\tag{$\overline{A1}$} 
	\end{align*}
\end{lem}

\subsection{Proof of Lemma \ref{U:alcInd}}
Assume that $\rho^1\iSv$ is associated to a first choice 
of $\zeta^1$ satisfying $\textup{$(A1)$}$,
and consider another choice $\zeta^2$.
By $\textup{$(A0)$}$, 
there exists $\ell\ge 1$ such that $\zeta^2 (\cD_\ell) \ge 1/2$. 
By $\textup{$(A1)$}$ applied to $\zeta^1$,
for some $c_J, t_J>0$ and $L> \ell$:
\begin{equation}
	\PR_{\zeta^2}(X_{t_J} \in dx;\; t_J< \ext\wedge T_{\cD_{L}})
	\ge c_J \zeta^1(dx).
	\EQn{U:cJ}{}
\end{equation}
By definition of $\rho^1\iSv$, 
for any $\rho > \rho^1\iSv$,
there exists $c_S, t_S>0$ and $L' \ge L$ such that:
\begin{equation}
	\frl{t\ge t_S}
	\PR_{\zeta^1}(t< \ext\wedge T_{\cD_{L'}})
	\ge c_S \exp[-\rho t].
	\EQn{U:cS}{}
\end{equation}
By combining \Req{U:cJ}, \Req{U:cS} and the Markov property,
we deduce:
\begin{equation*} 
	\limsup_{t>0} \dfrac{\exp[-\rho (t+t_J)]}
	{\PR_{\zeta^2}(t+t_J < \ext\wedge T_{\cD_L'})}
	\le \big( c_J\, c_S\, \exp[\rho t_J]\big)^{-1} < \infty.
\end{equation*}
By optimizing in $\rho$,
we deduce $\rho^2\iSv\le \rho^1\iSv$ 
and the equality by symmetry.
\\

Concerning assumption \textup{$(A3)$},
\Req{U:cJ} and the Markov property 
imply that for any $t\ge 0$ and $x\in \cX$:
\begin{align*}
	&
	\PR_{\zeta^2}(t+t_J< \ext)
	\ge c_J \PR_{\zeta^1}(t< \ext)
	\\\text{Thus} \qquad
	&
	\dfrac{\PR_x(t+t_J< \ext)}
	{\PR_{\zeta^2}(t+t_J< \ext)}
	\le (c_J)^{-1}\, \dfrac{\PR_x(t< \ext)}
	{\PR_{\zeta^1}(t< \ext)}.
\end{align*}
If assumption \textup{$(A3)$}\, holds for $\zeta^1$ and $E$, 
it thus holds also for $\zeta^2$ and the same $E$.
\hfill $\square$

%

\subsection{Proof of Lemma \ref{U:rg:sv}}
\label{U:sec:rgsv}
Assume \Req{U:Rg}. Let $x \in \cD\iSv$, 
$\rho:= -\frac{1}{t_{RG}} \ln(c_{RG}),$ 
$T_{L}:=\inf \left \lbrace \tp \ge 0,\, X_\tp \notin \cD_{L} \right \rbrace.$
By induction over $k \in \N$ and the Markov property:
\begin{align*}
\frl{k\ge 1}\qquad
\inf_{x\in \cD\iSv} 
\PR_{x} ( k \, t_{RG} < T_{L} ) 
\ge \exp(-  \rho \,k \, t_{RG}).
\end{align*}

Thus, for a general value of $t>0$:
\begin{align*}
	&\inf_{x\in \cD\iSv} \PR_{x} ( \tp < T_{L} ) 
	\ge \inf_{x\in \cD\iSv}
	\PR_{x} \left( \left \lceil \tfrac{\tp}{t_{RG}} \right \rceil t_{RG} < T_{L} \right)
	\ge \exp\left(-  \rho \,\left \lceil \frac{\tp}{t_{RG}} \right \rceil \, t_{RG}\right)\\
	& \hspace{2cm}
	\ge \exp(-  \rho \, (\tp + t_{RG})) 
	= c\iSv\; e^{-\rho\, \tp}
	\quad \text{ with }  c\iSv:= \exp(-  \rho \, t_{RG}) = c_{RG}.
	\tag*{$\square$}
\end{align*}

\subsection{Proof of Lemma \ref{U:lem:mixbar}: }
\label{U:mix2}

Let $\ell \ge L\iSv$
for which we apply $\textup{$(A1)$}$. 
By induction with the Markov property, 
it is quite straightforward to extend the property 
$\textup{$(A1)$}$ on $\cD_{\ell}$ with the same $L\iMix$,
$\tp^{(k)}:= k\times \tp$, 
$\cp^{(k)}:= \cp \times (\cp\,  \zeta(\cD\iSv)) ^{k-1}$.
Then, for any $t_\veebar>0$, 
we only need to apply this extension for some $k\ge1$ such that $\tp^{(k)} \ge t_\veebar$.
On the other hand, $(\overline{A1})$ clearly implies \textup{$(A1)$}\, (take $t_\veebar =0$), so that we have indeed proved
$\textup{$(A1)$}\Leftrightarrow (\overline{A1})$.
\hfill $\square$

\section{Two models to which our results apply}
\label{U:sec:Appli}
\setcounter{eq}{0}

\subsection{Birth-and-death process with catastrophes}
\label{U:sec:bd}

We choose to illustrate 
our result with this example 
for its clear simplicity.
%
In this birth and death process, 
the population can get extinct punctually at any time
during what we call a catastrophe.
These events happen 
at a rate depending on the current number 
of alive individuals.
Otherwise, the process gets extinct when there is only a unique individual 
that ends up dying.
To ensure uniqueness of a QSD, 
we will impose that 
the catastrophe rate is large enough 
when the population size is large. 
Biologically, we could imagine that 
the population is under the threat 
of some voracious predators, 
but can stay hidden 
as long as the population size 
is not too large.

In fact, one has now 
quite a complete description of quasi-stationarity 
for birth-and-death processes. 
It is proved in \cite{MSV14} that 
there exists a unique QSD 
for one dimensional birth and death processes 
if and only if \Req{U:Cmu}\, holds with a uniform constant $C(\mu) = C >0$. 
This equivalence is probably due to the fact that in these models, 
extinction can only occur once the process is inside some given compact set 
(i.e. once it has descended from infinity), 
as suggested in Theorem~19 in \cite{DP13}. 
Like in \cite{ChQSD} and as we will do, 
the authors of \cite{DP13} 
include direct extinction 
from any state of the birth-and-death process 
(what is called a "catastrophe").
Theorem~19 in \cite{DP13} 
states that the behavior of the process is 
the same if catastrophe 
only happens in a compact set.
In Theorem~4.1 of \cite{ChQSD}, 
the authors prove that, 
for a bounded catastrophe rate, 
there is descent from infinity 
(see notably \cite{BCM17}) 
iff \Req{U:Cmu}\, holds with a uniform constant $C(\mu) = C$.
This does not exclude however that 
\Req{U:Cmu}\, could hold 
without descent from infinity, 
which we prove with our technique.

\subsubsection{Description of the process}
\label{U:descBD}
$X$, the population size, 
is a time-homogeneous Markov Chain on $\mathbb{Z}_+$
where $\partial = 0$ is the absorbing state and $\cX = \N$. 
Given $X_0 = n \ge 1$, there is 
a death with rate $d_n>0$ (leading to $X = n-1$),
a birth with rate $b_n>0$ (leading to $X = n+1$)
and a catastrophe with rate $c_n \ge 0$ (leading to $X = 0$).
Since $c_1$ and $d_1$ play the same role 
(the transition is from $X=1$ to $X= \partial$),
we assume w.l.o.g. $c_1 = 0$. 
Actually, 
$d_1>0$ is not required in the following statements.

\begin{theo}
	\label{U:th:BD}
	\newcounter{Cn}
	
	Assume that: $\quad $
	for some $n \ge1$  (thus for all $n$)
	$\PR_n (\ext < \infty) = 1$
	\begin{align*} 
		\text{and} \qquad
\underset{n\rightarrow \infty}{\liminf} \,\, c_n > \inf_{k\ge 1} (b_k + d_k + c_k).
		\EQn{U:Cn}{c_n}
	\end{align*} 
	Then, the conclusions of Theorem~\ref{U:AllPho}, \ref{U:EtaECV} and \ref{U:QECV} hold.
\end{theo}

At least for some of the models, the speed of convergence towards the QSD cannot be uniformly bounded over all initial conditions, since:
\begin{prop}
	\label{U:BD.NU}
	We can define some positive values 
	for $(b_n, d_n, c_n)_{n \ge 1}$
	such that 
	\Req{U:Cn} holds
	and for which,
	whatever large the time $t>0$, 
	and whatever small the similarity threshold
	$\epsilon\in (0,1)$, 
	we can still find some initial condition $x\in \cX$ such that:
	\begin{align*}
		&\left\| \PR_x\left( X_\tp \in dy\, \big| \, \tp < \ext \right) 
			- \PR_{1}\left( X_\tp \in dy\, \big| \, \tp < \ext \right) \right\|_{TV}
		\ge  1- \epsilon.
	\end{align*}
\end{prop}

\noindent
The proof of Theorem~\ref{U:th:BD} and Proposition \ref{U:BD.NU} are achieved resp.
in Subsection~\ref{U:sec:thBD}-3.

\begin{rem}
	\label{U:SvEst}
Explicit values for $\rho_S$ can hardly be obtained
except for very specific models. 
Yet, it might be of interest to find,
depending on the specific model under consideration,
more precise estimates as our value 
$\inf_k (b_k + d_k + c_k)$.
This upper-bound 
comes from a survival estimate
of the simplest form:
the process reaches 
some position $k$
(as optimal as we need)
on which to stay up to (large) time $t$.
We are a priori very far 
from a necessary and sufficient condition:
it seems hardly possible to infer generically
the level of catastrophe rate
that affects the process 
as it evolves at large values.
\end{rem}

\begin{rem}
	\label{U:nonExpl}
	By the condition $\PR_n (\ext < \infty) = 1$, we mean that the process is non-explosive. With:
	$\quad 
	T_\infty:= \textstyle{\lim_{n \rightarrow \infty} } \inf\left \lbrace \tp\ge 0;\; X_\tp \ge n \right \rbrace,\quad$
	our condition means that "for some $n \ge1$ (thus for all),   
	$\PR_n(T_\infty = \infty) = 1$".
	Clearly, 
	this property holds
	provided it holds
	in the associated birth-and-death process without catastrophe
	(i.e. imposing $c_n\equiv 0$).
	The simplest case is then 
	when the sequence $b_n /n$ is upper-bounded.
	%
	We refer to Theorem~5.5.2 in \cite{M16} for a more general condition 
	(still deduced from the case without catastrophe).
	
	We do not exclude that it could hold more generally.
	Yet, for catastrophe to play a role in this condition
	would require the family of catastrophe rate $(c_n)$ to quickly reach very large values
	as $n$ goes to infinity. 
	It does not seem likely in practical applications.
\end{rem}

\begin{rem}
	\label{U:expl}
	Considering $\bar{\ext}:=  \inf\left \lbrace \tp \ge 0;\; X_\tp = 0\right \rbrace \wedge T_\infty$
	as the extinction epoch,
	our theorem extends to the case $T_\infty < \infty$.
	It also extends to models  
	where catastrophes do not entirely exterminate the population.
	Assume for instance that after a catastrophe, 
	from a population of size larger than some $K\ge 1$, 
	only $K$ individuals are to survive.
	We can keep the extinction for population of size initially lower than $K$, 
	but it's not very significant here.
	Then \textup{$(A2)$}\, can easily be adapted with $K \in E = [\![1, \ell\iET ]\!]$. 
	The proof of the other assumptions remains the same.
\end{rem}

\begin{rem}
	The alternative conditions given in \cite{CVly2} 
	seem also very efficient to obtain 
	Theorem \ref{U:th:BD}. 
	Since $\rho\iSv$ is finite,
	this will certainly not be the case 
	for the ones in \cite{FRS19}. 
\end{rem}

\subsubsection
[Proof of the exponential convergence]
{Proof of Theorem~\ref{U:th:BD}}
\label{U:sec:thBD}
By \Req{U:Cn}, let $k, \ell\iET\ge 1$ and $\rho\iET>0$ be such that:
\newcounter{RhoBD}
\begin{align*}
	&0 < \widetilde{\rho}\iSv:= b_k + d_k + c_k
	\,<\, \rho\iET
	\,<\, \textstyle{\inf_{\{n \ge \ell\iET +1\}}} c_n  
	:= \widetilde{\rho}\iET.
	\EQn{U:RhoBD}{BD}
	\\
	&\hspace{1cm}
	\zeta:= \delta_k,\quad
	\cD\iSv:= \left \lbrace k \right \rbrace, \quad
	E:=  [\![1, \ell\iET ]\!].
\end{align*}
Let $\cD_\ell = [\![1, \ell\vee k]\!]$, for $\ell\ge 1$. 
In the following, we ensure $\mathbf{(A)}$ (where $\zeta(\cD\iSv) >0$ is obvious).
First, \textup{$(A0)$}\, is obvious. 

\paragraph{Proof of  \textup{$(A1)$}\, and \textup{$(A3)$}}	
\textcolor{white}{:}

Let $n \ge k$. Consider 
\begin{align*}
	\partial^n:= \left \lbrace 0 \right \rbrace \cup [\![n+1, \infty [\![, \qquad
	\ext^n:= \inf\left \lbrace \tp\ge 0;\; X_\tp \in \partial^n \right \rbrace.
\end{align*}
Then the process $Y_\tp:= X_\tp\, \mathbf{1}_{\{ \tp<\ext^n \}}$ 
is a Markov Chain on the finite space $[\![0, n ]\!]$, absorbed at $\partial = 0$.
Since $\frl{\ell\ge 1} d_\ell >0, b_\ell >0$, this Chain $(Y_\tp)$ is irreducible and it is elementary to prove that: \newcounter{MixBD}
\begin{align*}
	\frl{t_Y>0} \Ex{c_Y >0}\quad 
	\frl{i, j \in [\![1, n ]\!]} \quad 
	\PR_i(Y_{t_Y} = j) \ge c_Y.
	\EQn{U:MixBD}{Y_\tp}
\end{align*}

With $j:= k$ and $n:= \ell\iMix \vee k$, 
\Req{U:MixBD}\, clearly implies \textup{$(A1)$}\,
(with parameters $\ell = \ell\iMix, L= L\iMix, c = c\iMix, t = t\iMix$).
We can indeed choose  $\zeta:= \delta_k$, $L\iMix:= n$,
$t\iMix = 1$ (arbitrary),
and $c\iMix$ the value of $c_Y$ 
associated to the choice of $t_Y = t\iMix$.

With $i:= k$ and $\ell = \ell\iET$, 
\Req{U:MixBD}\, 
and the Markov property imply $\textup{$(A3)$}$
for any $E$,
because:
\begin{align*}
	\frl{t>0}
	\frl{j \in [\![1, \ell\iET ]\!]} \quad 
	\PR_j(t< \ext)
	& \le (1/c_Y)\ltm \PR_k(X_{t_Y} = j;\; t+t_Y< \ext)
	\\&
	\le (1/c_Y)\ltm\PR_k(t< \ext).
\end{align*}

\paragraph{Proof of  \textup{$(A2)$}}
\label{U:BDReT}
\textcolor{white}{:}

By \Req{U:RhoBD},
the catastrophe rate is larger than $\widetilde{\rho}\iET$
as long as the process remains outside $E$.
It implies that we can upper-bound $\ext \wedge \tau_{E}$
by an exponential variable 
with rate $\widetilde{\rho}\iET$. 
Thus: 
\begin{align*}
	\frl{\tp>0} \frl{n \ge \ell\iET +1}
	\quad
	\PR_n(\tp < \tau_{E} \wedge \ext) 
	\le \exp(-\widetilde{\rho}\iET \, \tp).
\EQn{A2_1}{}
\end{align*} 
It is classical -by Fubini Theorem, and the integral expression of the exponen\-tial- to relate the exponential moment with the repartition function by:
\begin{align*}
	&
	\E_{n} \left( 
	\exp\left[\rho\iET\, (\tau_{E}\wedge \ext) \right] \right) 
	= 1 + \rho\iET \, \int_{0}^{\infty} \exp[\rho\iET\, \tp] \, \PR_n(\tp < \tau_{E}\wedge \ext) \, dt.
	\EQn{A2_2}{}
\end{align*}
By \Req{A2_1}\,and \Req{A2_2}, we conclude:
\begin{align*}
	&\frl{n \ge \ell\iET +1} \quad
	\E_n \left(\exp\left[\rho\iET\, (\tau_{E}\wedge \ext) \right] \right)  
	\le 1 + \rho\iET \, \int_{0}^{\infty} \exp[-(\widetilde{\rho}\iET -\rho\iET)\, \tp]\, dt
	\\ & \hspace{7cm}
	= 1 +\ \{ \rho\iET \ /\ (\widetilde{\rho}\iET -\rho\iET)\} < \infty.
\end{align*}
\paragraph{Proof that for any $k\ge 1$:
	$b_k + d_k + c_k\ge \rho\iSv$} 
\textcolor{white}{:}

Immediately,  by \Req{U:RhoBD}: 
\begin{align*}
	&\frl{\tp \ge 0}\hspace{.5cm}
	\PR_{k} (X_\tp = k;\; \tp < \ext
	\wedge T_{\cD_k})
	\ge \PR_{k} (\frl{\spr\le \tp} X_\spr = k) 
	= \exp(-  \widetilde{\rho}\iSv \, \tp).
\end{align*}

\subsubsection
[Proof of the non-uniformity]
{Proof of Proposition \ref{U:BD.NU}}
\label{U:sec:BD.NU}

We consider one of the simplest choice, 
which is to take $b_n$, $d_n$ linear in $n$ (the classical Malthus' growth model,
without competition)
and $c_n$ constant for $n \ge 2$. 
We can then choose arbitrarily:
\begin{align*}
	&b_1, d_1, \bar{b}, \bar{d} \in (0, \infty)^5, \quad 
	c_2 > (b_1 + d_1),
	\EQn{Bdc}{b,d,c}\\
	\quad \text{ with } &c_1 = 0,\qquad
	\frl{n \ge 2} \quad b_n:= \bar{b}\, n,\, 
	d_n:= \bar{d}\, n,\, 
	c_n:= c_2.
\end{align*}

\Req{U:Cn} 
is clearly satisfied. 
There is no explosion for this model, 
so that extinction happens a.s. (note Remark \ref{U:nonExpl} of Subsection \ref{U:descBD} on this aspect).

We shall only need to consider transitions between values of the form $2^n, n \ge 2$.
Let:
\begin{align*}
	T_n &:= \inf \left \lbrace \tp\ge 0;\; 
	X_\tp \le 2^{n-1} \, or\, X_\tp \ge 2^{n+1} \right \rbrace,
	\EQn{Tn_1}{}	
	\\
	\tau_n &:= \inf \left \lbrace \tp\ge 0;\; 
	X_\tp \le 2^n \right \rbrace.
	\EQn{U:TTn}{}
\end{align*}

\noindent
We use the following lemma, 
whose proof is deferred 
after the one of Proposition \ref{U:BD.NU}:
\begin{lem} 
	\label{U:TnTv}
	For some $\uV > 0$, it holds: 
	$\quad\limInf{n} \PR_{2^n} (T_n \le \uV) = 0.$
\end{lem}

For given $\tp, \epsilon >0$, let 
$\quad K:= \left \lfloor \tp/\uV \right \rfloor + 1 \quad$
%
and $N \ge 1$ (by Lemma \ref{U:TnTv}) such that: 
\begin{align*}
	&\hcm{2}
	\PR_1 (X_\tp\le 2^N	\, \big| \,	 \tp< \ext) \ge 1 - \epsilon / 2, 
	\EQn{U:X2te}{}
	\\&
	\frl{n \ge N} \quad \PR_{2^n} (T_n \le \uV) 
	\le \epsilon \times e^{-(c_2-d_1)\tp} / (4\,K).
	\EQn{TnTv_1}{}
\end{align*}
With initial condition $\quad  x:= 2^{N+K+1},\quad $ 
in order that $X$ reaches $2^N$ before time $t \le K\, \uV$, 
it must at least once have got from $2^{N+k}$ to $2^{N+k-1}$ during a time-interval less than $\uV$, 
for some $1\le k \le K+1$. 
With the Markov property, this implies, with \Req{Tn_1}, \Req{U:TTn}, \Req{TnTv_1}\,
and the fact that the extinction rate is always lower-bounded by $d_1$:
\begin{align*}
	&\PR_x ( \tau_N \le \tp ) 
	\le \textstyle{\sum_{\{k \le K+1\}}} \PR_{2^{N+k}} (T_{N+k} \le \uV )
	\le \epsilon \; e^{-(c_2-d_1)\tp} / 2.
	\\&
	\PR_x ( \tau_N \le \tp;\; \tp < \ext) \le e^{-d_1\, \tp}\, \PR_x ( \tau_N \le \tp )
	\le \epsilon \; e^{-c_2\,\tp} / 2.
\end{align*}
Since the extinction rate is upper-bounded by $c_2$:
$\PR_x ( \tp< \ext ) 
\ge  e^{-c_2\, \tp}$.
\\
This implies 
$\PR_x ( \tau_N \le \tp \, \big| \, \tp < \ext) \le \epsilon/2.$
Therefore, with also $\Req{U:X2te}$: 
\begin{align*}
	\| \delta_1 A_\tp  - \delta_x A_\tp \|_{TV}
	&\ge \PR_1 (X_\tp\le 2^N	\, \big| \,	 \tp< \ext) - \PR_x (X_\tp\le 2^N	\, \big| \,	 \tp< \ext) &&
	\\&
	\ge 1- \epsilon /2  - \epsilon/2 \ge 1- \epsilon.
	&\square&
\end{align*}

\paragraph{Proof of Lemma \ref{U:TnTv}}:\\
With initial condition $2^n$, we can decompose $X$ as a semi-martingale, 
up to time $\tp\wedge T_n$:
\begin{align*}
	\frl{\tp>0}\quad
	X_{\tp\wedge T_n}:= 2^n + \textstyle{\int_{0}^{\tp\wedge T_n}} (\bar{b} - \bar{d}) \, X_s\, ds
	+ M_{\tp\wedge T_n}, 
	\EQn{TnTv_2}{}
\end{align*}
where $(M_{\tp\wedge T_n})_{\tp}$ is a martingale 
with bounded quadratic variation, with $\Req{Tn_1}$:
\begin{align*}
	<M> _{\tp\wedge T_n} = \textstyle{\int_{0}^{\tp\wedge T_n}}  (\bar{b} + \bar{d}) \, X_s\, ds
	\le (\bar{b} + \bar{d})\, 2^{n+1}\, \tp.
	\EQn{U:MQV}{}
\end{align*}
Let $\quad \uV:= (8\, |\bar{b} - \bar{d}|\vee 1)^{-1}\quad$ so that, 
by $\Req{Tn_1}$, a.s.:
\begin{align*}
	\frl{\tp\le \uV}\quad
	\, \Big| \, \textstyle{\int_{0}^{\tp\wedge T_n}}  (\bar{b} - \bar{d}) \, X_s\, ds \, \Big| \,
	\le |\bar{b} - \bar{d}|\, 2^{n+1}\, \uV \le 2^{n-2}.
	\EQn{U:bdX}{}
\end{align*}
\begin{align*}
	\hcm{.3}
\PR_{2^n} (T_n \le \uV)  
	&	\le \PR_{2^n} \left( \textstyle{\sup_{\{\tp\le \uV\}}} M_{\tp\wedge T_n} \ge 2^{n-2} \right)
	\hcm{1} \text{ by } \Req{TnTv_2} \text{ and } \Req{U:bdX}
	\\&
	\le 2^{-(2\, n - 4)}\; \E_{2^n} ( <M> _{\uV\wedge T_n} )
	\text{ by Doob's inequality }
	\\&
	\le 2^{-(2\, n - 4)}\; (\bar{b} + \bar{d})\, 2^{n+1}\, \uV
	\hcm{1} \text{ by } \Req{U:MQV}
	\\&\hcm{0.5}
	= \dfrac{ 4\, (\bar{b} + \bar{d}) }{|\bar{b} - \bar{d}|\vee 1} 
	\, 2^{-n} 
	\underset{n\rightarrow \infty} {\longrightarrow} 0
	\quad \text{with the definition of } \uV.
	\tag*{$\square$}
\end{align*}

\subsection{Adaptation of a population to its environment: application to a diffusion process}
\label{U:sec:Adapt}

In  this illustration, 
the notion of being in a mal-adapted region 
is quite intuitive
and the criteria for the exponential convergence 
to a unique QSD rather natural.
Again, the general proof
for this illustrative example
is unclear without our techniques,
except maybe with those of \cite{CVly2}.
Yet, in this case, 
it is presumably quite technical 
to find a proper Lyapunov function
(although our argument proves in fact that they exist).
In fact,
our control is deduced
from local bounds
ensuring both a rapid escape from several specific local domains
together with sufficiently low transition rates
between these domains.

\subsubsection{Presentation  of the model}
\label{U:sec:presDiff}

We consider a simple coupled process describing 
the eco-evolutive dynamics of a population.
We model the population size by a logistic Feller diffusion $(N_t)_{t\ge0}$
where the growth rate $(r(X_t))_{t\ge0}$ 
is changing randomly.
Namely, 
the adaptation of the population 
and the change of the environment
are assumed to act on a hidden process $(X_t)$ in $\R^d$,
from which the growth rate is deduced.
For simplicity, 
we will assume that $X_t$ evolves as a continuous Markov process
driven by some Brownian Motion and a drift 
(possibly depending on $N$ and $X$).
For very low values of $r(X_t)$,
it is expected that the population shall vanish very quickly.
It would thus not change much of the result 
to introduce an absorbing boundary 
at some threshold of mal-adaptation.
Yet, we want our result to be independent 
of any such truncation of the trait space
and say that this large extinction is sufficient in itself 
to bound the mal-adaptation,
while highlighting that the initial condition indeed matters here.

In a general setting, 
the process can be described as:
\begin{align*}
	(S) \left\{
	\begin{aligned}
		&dN_t = (r(X_t) - c\ N_t)\ N_t\ dt + \sigma \ \sqrt{N_t}\ dB^N_t
		\\&
		dX_t =  b(X_t, N_t)\ dt + \theta(X_t, N_t) \ dB^X_t
	\end{aligned}
	\right.
\end{align*}
with initial conditions $(n, x)$, $B^N$ and $B^X$ two independent Brownian Motions,
$c, \sigma >0$, 
and $r, b, \theta$ being locally H\"older continuous functions.
We also require that $\theta$ is locally elliptic,
in the following sense: 
for any compact set $K$ of $\R^{d+1}$,
there exists $\bar{\theta}>0$ 
such that for any $(n, x)\in K$
and $\xi\in \R^d$:
$\sum_{i, j} \theta_{i, j}(n,x) \xi_i \xi_j 
\ge \bar{\theta}\, |\xi|^2$.

\begin{theo}
	\label{U:th:ECVX}
Consider the process $(X, N)$ with the notations specified above
	 and the assumption that
	$\limsup_{\|x\| \rightarrow \infty} r(x)= -\infty$.
	Then, 
 all the results of Subsection \ref{U:sec:ECV} hold. 
 In particular, 
 there is exponential convergence 
 in total variation 
 of the MCNEs to the unique QSD.
\end{theo}

It is also not much more costly to introduce catastrophes,
arising at rate $\rho_c(x, n)$,
leading to the complete extinction of the population.
Partial extinction of the population (with jumps on the population size),
are however quite more technical to deal with
(because the Harnack inequality is not as obvious).
In \cite{AV_disc},
where the focus is on processes with jumps,
we shall present techniques that makes it much more manageable.

The main issue for this model is to specify the conditions 
for \textup{$(A2)$}\, to hold.
We provide in Subsection \ref{U:sec:etD}
a way to prove it 
in a strong case where it holds for any $\rho$,
i.e. for the following Theorem \ref{U:th:ECVXN}.
For diffusions like this, 
\textup{$(A1)$}\, and \textup{$(A3)$}\, 
may be deduced quite roughly 
thanks to the Harnack inequalities, 
as presented in the next subsection.
In order to satisfy assumption \textup{$(A3)$},
there is no additional restriction 
on the set $E$, 
so that the only requirement on $E$
is for \textup{$(A2)$}.

\subsubsection{Harnack inequalities for \textup{$(A1)$}\, and \textup{$(A3)$}}
\label{U:sec:Harn}

In the following, we say that a process $(Y_t)$ on $\mathcal{Y}\subset \R^d$
with generator $\mathcal{L}$
(including possibly an extinction rate $\rho_c$)
satisfies Assumption $(H)$ if the following property holds:

\textsl{
	Consider any path-connected open relatively compact sets $\cD, \cD'\subset \mathcal{Y}$,
	such that $\overline{\cD}\subset \cD'$,
	with $C^\infty$ boundaries,
	and such that for any point $x\in \partial \cD'$, there exists a closed ball $C \in \R^d$ 
	(of non-empty interior)
	such that $C \cap \overline{\cD'} =\{y\}$.
	For any  $0<t_1<t_2$ and non-negative $C^2$ constraints:
	$u_{\partial \cD'}: (\{0\} \times \cD') \cup ([0, t_2] \times \partial \cD')
	\rightarrow [0, \infty)$,
	there exists a unique solution $u\in C^{1, 2}((0, t_2)\times \cD')\cap C^{0}([0, t_2]\times \overline{\cD'})$ to the problem:}
\begin{align*}
	&
	\partial_t u (t, y) = \mathcal{L} u (t, y) 
	&&\text{ on } [0, t_2] \times \cD';\;
	\\&
	u(t, y) = u_{\partial \cD'}(y) 
	&&\text{ on } (\{0\}\times \cD') \cup ([0, t_2] \times \partial \cD').
\end{align*}
\textsl{It is non-negative on $int(\cD')$
	and it satisfies,
	for some $C = C(\mathcal{L}, t_1, t_2, \cD, \cD')>0$ independent of $u_{\partial \cD'}$:}
\begin{align*}
	{\textstyle \inf_{y\in \cD} u(t_2, y) 
		\ge C\, \sup_{y\in \cD} u(t_1, y).}
\end{align*}

\begin{prop}
	\label{U:th:ECVXN}
	Assume that Assumption $(H)$ holds,  and
	$\limsup_{\|x\| \rightarrow \infty} r(x) = -\infty$.
	Then, all the results of Subsection \ref{U:sec:ECV} hold, and we have in particular exponential convergence in total variation of the MCNE to the unique QSD.
\end{prop}

Assumption $(H)$ is crucial 
for the proofs of both $(A1)$ and $(A3)$,
yet not at all for the one of $(A2)$.
As stated in the next proposition, 
the form of the equation for $N$
is the main ingredient.

\begin{prop}
	\label{A2XN}
Assume that $(X, N)$ is a càdlàg process
on $\R^d\times \R_+$ such that $N$ is solution to:
$$dN_t = (r(X_t) - c\ N_t)\ N_t\ dt + \sigma \ \sqrt{N_t}\ dB^N_t,$$
where $B^N$ is a Brownian motion.
Assume that $\ext$ is upper-bounded by $\inf\{t\ge 0; N_t =0\}$.
Provided that $\limsup_{\|x\| \rightarrow \infty} r(x) = -\infty$,
it holds that for any $\rho>0$, there exist $n>0$ such that:
$$\underset{x\in \cX}{\sup} \;\E_{x} \left( 
\exp\left[\rho\, (\ext\wedge \tau_{E}) \right] \right) < \infty,$$
where $E:= \bar{B}(0, n)\times [1/n, n]$, 
with  $\bar{B}(0, n)$ the closed ball of $R^d$ centered in $0$
and of radius $n$.
\end{prop}
The main elements of the proof are given in Subsection  \ref{U:sec:etD},
with the most elementary arguments deferred to the Appendix.

\paragraph{From Proposition \ref{U:th:ECVXN} to Theorem \ref{U:th:ECVX}}
Looking at the system $(S)$ of equations,
 Assumption $(H)$ 
holds for the generator: 
\begin{align*}
	\mathcal{L} f(x, n) 
	&:= [r(x) - c\, n]\, \partial_n f(x, n) + b(x, n)\, \partial_x f(x, n)
	\\&\hcm{2}
	+ n \ \sigma^2 / 2 \times \Delta_n f(x, n)
	+ \theta^2(x,n) / 2
	\times \Delta_x f(x, n).
\end{align*}

The proof for existence and uniqueness 
of the solution $u$
for such second-order partial operator
with H\"older coefficients 
and elliptic diffusion coefficient
can be found for instance 
in Corollary 2, Section 4, Chapter 3
of \cite{FR}. 
It also ensures that the solution 
has two continuous $x$-derivatives 
and one continuous $t$-derivative.
The fact it is non-negative 
is then a consequence of the Maximum principle
(cf. e.g. Theorem 1, Chapter 2 in \cite{FR}).
Finally, the comparison comes 
from the parabolic Harnack inequality,
exploiting the regularity of $u$ 
and the fact it is non-negative.
For its proof,
we refer to Theorem 1.1 in \cite{KS80}.
The Harnack inequality 
on the open and path-connected set $\cD$
is not too difficult to deduce
from the local Harnack inequalities these authors provide.
One essentially needs to cover $\bar{\cD}$
by a finite number of balls 
included in $\cD'$
on which the local inequalities can be applied.
For any points $x, y \in\bar{\cD}$,
we can then construct a path between them
for which the number of visited such balls
is uniformly bounded.
The interval $[t_1, t_2]$
shall then be split 
into as many time-intervals
and the local Harnack inequalities 
applied recursively to conclude 
that Assumption $(H)$ holds.
\\

Assumption $(H)$ shall hold more generally, 
notably under the Hörmander condition 
instead of the condition of ellipticity
(cf e.g. \cite{Horm}).
Lots of articles are dedicated 
to prove such estimates under various conditions.

\subsubsection{Proof of Theorem \ref{U:th:ECVXN}}
In order to conveniently exploit Assumption $(H)$,
we choose to consider $(\cD_\ell)$  as some sequence 
(can be anyone) 
 of strictly increasing compact and path-connected sets 
with $C^\infty$ boundaries 
whose union is $\mathcal{Y}:= \R^d \times \R_+^*$.

Such a sequence clearly satisfies $(A0)$
and any set of the form $\bar{B}(0, n)\times [1/n, n]$
is a subset of $\cD_\ell$ for $\ell$ sufficiently large.
$(A2)$ thus also hold for this process.

\paragraph{Assumption $(H)$
	with $Y_t = (X_t, N_t)$
	implies \textup{$(A1)$}}
\textcolor{white}{:}

For some non-negative $C^\infty$ function $f$ with support in $\cD_\ell$, 
we apply Assumption $(H)$ with 
$\cD:= \cD_\ell $,
$\cD':= \cD_{\ell+1}$
with $u_{\partial \cD'}(t, y):= f(x)$ on $\{0\}\times \cD_\ell$
and $u_{\partial \cD'}(t, y):= 0$ on $\R_+^* \times \partial\cD_{\ell+1}$.
The solution $u$ we obtain is identified thanks to Itô formula as:
$u(t,y):= \E_y\left( f(Y_t);\; t< \ext^{\ell+1}\right)$, 
with an additional extinction when the process exits $\cD_{\ell+1}$.
Applying Harnack inequalities implies thus that for any $y\in \cD_\ell$,
and some reference $y_0\in \cD_1$:
\begin{align*}
	\E_y\left( f(Y_{t_2});\; t_2< \ext^{\ell+1}\right)
	\ge c\; \E_{y_0} \left( f(Y_{t_1});\; t_1< \ext^{1}\right).
\end{align*}
Since it is classical that $\PR_{y_0} \left( Y_t \in \cD_1;\; t< \ext^{1}\right) >0$, 
we can obtain 
a probability measure $\zeta$, independent of $\ell$, such that
(since $c$ does not depend on $f$):
\begin{align*}
	\frlq{y\in \cD_\ell} 
	\E_y\left( Y_t \in dy;\; t< \ext^{\ell+1}\right)
	\ge c\; \zeta(dy).
\end{align*}

\paragraph{Assumption $(H)$
	with $Y_t = (X_t, N_t)$
	implies \textup{$(A3)$}}
\textcolor{white}{:}

The proof of \textup{$(A3)$}\, is a bit similar but
much more technical because the reference measure 
is now in the upper-bound, 
so that we can no longer neglect trajectories 
exiting $\cD'$. 
W.l.o.g., 
we consider $E$ 
to be of the form $\cD_\ell$ 
for $\ell$ sufficiently large. 
Since the support of $\zeta$ is included in $\cD_1$, 
we wish to prove 
that there exists $\cp >0$
such that for any 
$y_1 \in \cD_1$ and $y_E\in E$:
\begin{align*}
	&
	\PR_{y_E}\left( Y_{\tp}\in dy;\; \tp < \ext \right)
	\le \cp \, \PR_{y_1}\left( Y_{\tp_{\alpha}}\in dy;\; 
	\tp_{\alpha} < \ext \right),
	\EQn{Reach}{}
\end{align*}
where we can choose here $0< \tp_{\alpha}< \tp$ arbitrary 
($\cp$ depending on this choice). 
$\Req{Reach}$ directly implies \textup{$(A3)$}\,
with the functions
$f_s(y) = \PR_y(s - t_\alpha <\ext)$, 
and the Markov property.

In the step 4 of the proof given in Section~4 of \cite{ChpLyap}, \mbox{N. Champagnat} and D. Villemonais used a trick to obtain results such as \Req{Reach}. 
Their idea is to apply the parabolic Harnack inequality on some regular and compact domain  $\mathcal{R}$ such that
$E \subset \mathcal{R} \subset  \cX$ and $d(E, \partial \mathcal{R}) >0$ while
approximating the function:
\begin{align*}
		u(\tp, y):= \E_{y}\left(  f(Y_{\tp});\; \tp < \ext \right),
\quad \text{ with } \; \tp  \ge t_\mathcal{R},\; y \in \mathcal{R},
\end{align*}
defined for some non-negative $f \in C^\infty(\cX)$ 
and any choice of $0<t_\mathcal{R}< \tp_{\alpha}$.
Although we can prove (as they do) that $u$ is continuous, 
it is a priori not regular enough to apply Harnack inequality directly.
Thus, we approximate it
on the parabolic boundary \mbox{$[t_\mathcal{R},\, \infty) \times \partial\mathcal{R}\,$}
\mbox{$\bigcup\, \{t_\mathcal{R}\} \times \mathcal{R}$}
by the family $(U_k)_{k\ge 1}$ of smooth non-negative functions ($\mathcal{C}_+^\infty$ w.l.o.g.).
We then deduce approximations of $u$ in $[t_\mathcal{R},\, \infty) \times \mathcal{R}$ 
by (smooth) solutions of:
\begin{align*}
	&
	\partial_\tp u_k (\tp, y) - \mathcal{L} u_k (\tp, y)= 0,\qquad
	\tp \ge t_\mathcal{R},\; y \in \mathcal{R}^\circ
	\\& u_k(\tp, y) = U_k(\tp, y),
	\hcm{2} \tp \ge t_\mathcal{R},\; y \in \partial\mathcal{R}, \quad 
	\text{ or } \tp = t_\mathcal{R},\, y \in \mathcal{R}.
\end{align*}
By Assumption $(H)$, 
the constant involved in the Harnack inequality 
does not depend on the values on the boundary.
Thus, it applies with the same constant 
for the whole family of approximations $u_k$.
We refer to the proof in \cite{ChpLyap} to state that
the Harnack inequality then extends to the approximated function $u$,
where the regularity of $u \in C^{1, 2}$ 
is required to apply the Itô formula on the process $u(t-s, Y_s)$.
Thus, \Req{Reach}\, indeed holds (where we could have chosen any $y\in E$).

Now we have concluded that Assumption $\mathbf{(A)}$ 
holds,
so that the conclusions of Theorems \ref{U:AllPho}
to \ref{U:QECV} hold.
\hfill $\square$

\subsubsection{Proof of Proposition \ref{A2XN}: escape from the transitory domain}
\label{U:sec:etD}

The purpose of this section is 
to demonstrate how to prove \textup{$(A2)$}\,
when, depending on the position 
in the transitory domain,
there are
various reasons for a quick escape.
To combine several local estimates,
dealing with suprema 
in the initial condition 
of exponential moments 
appear much more convenient
than Lyapunov estimates,
see Appendix A.
Moreover, 
these exponential moments 
can be naturally deduced
from probabilities
of retention (or transfer) 
in the transitory domain 
for a finite given time,
see Appendix B, C or D.
\\

\noindent
\textbf{Decomposition of the transitory domain} 

\begin{figure}[h]
	\begin{center}
		\includegraphics[width = 16cm, height = 5cm]{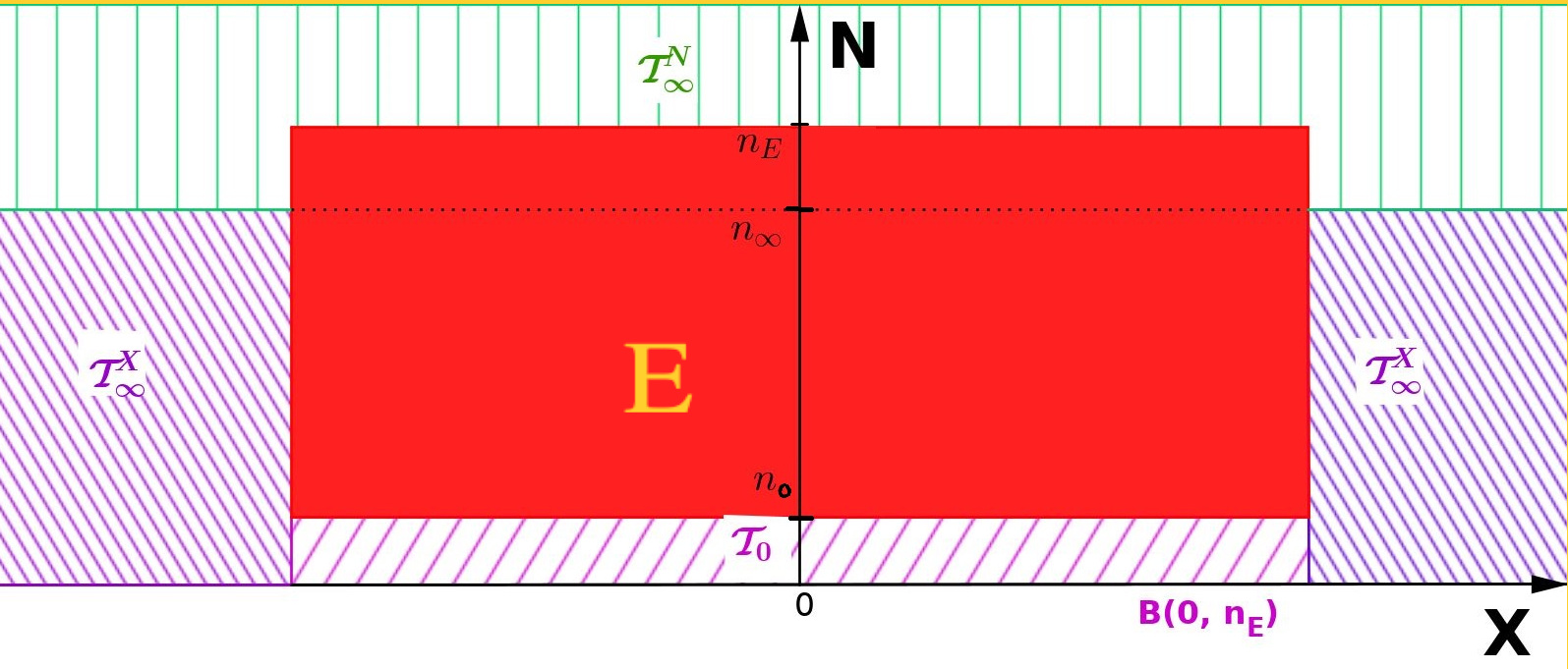}
		\caption{subdomains for \textup{$(A2)$}}
		\label{U:eTD}
	\end{center}
\end{figure}

In our choice for $E$,
with three parameters $n_0<n_\infty<n_E$ to be fixed,
its complementary $\cT = \cX \setminus E$ 
is made up of 3 subdomains:
"$y=\infty$", "$y=0$", and "$\|x\| = \infty"$, according to figure \ref{U:eTD}.
They are formally defined as follow:
\begin{itemize}
	\item $\cT^N_\infty:= 
	\left \lbrace \R^d \setminus B(0,  n\iET) \right \rbrace
	\times (n_{\infty}, \infty)$
	$\bigcup \; B(0,  n\iET)\times (n\iET, \infty)\hfill $("$y= \infty$"),
	\item $\cT_0:= B(0,  n\iET) \times [0,n_0]\hfill$ ("$y = 0$"),
	\item $\cT^X_\infty:= \left \lbrace \R^d \setminus B(0,  n\iET) \right \rbrace \times (n_0, n\iET]
	\hfill$ 	 	 ("$\|x\| =\infty$").
\end{itemize}

Essentially, we will need to choose 
$n_{\infty}$ sufficiently large 
to have the property of descent from infinity for $\cT^N_\infty$;\;
$n\iET > n_{\infty}$ sufficiently large to have a growth rate so low that the population cannot maintain itself in $\cT^X_\infty$;\;
$n_0$ sufficiently small to prove that
the population can hardly survive after entering $\cT_0$.
Thus, the process will escape each region with an exponential moment.
Yet, we also need to prove that the process will not circulate
between the different transitory areas. 
Therefore, 
we will set these areas such that 
at least some of the transitions 
(those associated with an increase of the  population size) 
happens with so low probability that Theorem \ref{U:eTD} holds true.\\ 

For each of these domains,
we define the following exponential moments
that we shall relate by specific inequalities.
Let $t_{h}>0$ 
(a threshold needed to ensure the boundedness)
and 
$\widehat{\tau}_{E}:= \tau_{E} \wedge \ext\wedge t_{h}$
(remember that $\tau_{E}$ 
is the hitting time of $E$):
\begin{itemize}
	%
	\item $\cE^N_\infty 
	:= 	\sup \{
	\E_{(x, n)}[\exp(\rho \; \widehat{\tau}_{E})];\;
	(x, n) \in \cT^N_\infty\}$,  	
	%
	\item $\cE^X_\infty 
	:= \sup \{\E_{(x, n)}[\exp(\rho \; \widehat{\tau}_{E})];\;
	(x, n) \in \cT^X_\infty\}$,
	\item  $\cE_0 
	:= \sup \{\E_{(x, n)}[\exp(\rho \; \widehat{\tau}_{E})];\;
	(x, n) \in \cT_0\}$.
\end{itemize}
Implicitly, we assume $\rho$ to be given. 
Then, $\cE^N_\infty,\, \cE^X_\infty$ and $\cE_0$ can be seen as functions of 
$n_0$, $n_{\infty}$ and  $n\iET$.
These values are to be specified depending on $\rho$
for the proof of $(A2)$ to hold.
The dependency in $t_{h}$
shall be negligible as $t_{h} \rightarrow \infty$. 
\\

\noindent
\textbf{A set of inequalities associating the local bounds}

The local exponential moments that we introduce are related 
thanks to the three following propositions,
obtained from local bounds mentioned in the following three lemmas. 
We refer to Appendices A, B, C and D
to see first how to deduce \textup{$(A2)$}\, 
from the three propositions that follow,
and then respectively 
for the (technical) proofs of the propositions
(including the lemmas):

\begin{prop}
	\label{U:cEinf}
	Given any $\rho>0$,
	we can define $n_{\infty} > 0$ and $C\ge 1$ such that,
	whatever $n\iET> n_{\infty}$ and $t_{h}>0$:
	$\quad \cE^N_\infty \le C \, \left( 1+ \cE^X_\infty \right).$
\end{prop}

\begin{prop}
	\label{U:cEXinf}
	Given any $\rho,\, \epsilon,\, n_{\infty}>0$,
	we have, for some $C\ge 1$ (in fact independent of any parameters),
	and any $n\iET$ sufficiently large and $t_{h}>0$:\\
	$\quad \cE^X_\infty \le C \, \left( 1 + \cE_0 \right)
	+ \epsilon \,\cE^N_\infty.$
\end{prop}

\begin{prop}
	\label{U:cE0}
	Given any $\rho,\, \epsilon,\, n_{\infty}>0$,
	we have, for some $C \ge 1$,
	any $n_0$ sufficiently small,
	 any $n\iET\ge n_{\infty}$
	 and any $t_{h}>0$:
	$\quad \cE_0 \le C  + \epsilon \, \left( \cE^N_\infty +\cE^X_\infty \right).$
\end{prop}

\noindent
\textbf{The associated elementary bounds on finite time}

The main ingredient for these propositions
are simple comparison properties 
that are specific 
to each of the transitory domain.
By focusing on each of the domains separately
(with the transitions between them),
we can highly simplify our control 
on the  dependency of the processes. 
Specific autonomous one-dimensional processes
indeed act as upper-bounds for each of the domains.
The values of $(X_t)$ 
do not affect these auxilliary processes, 
but only the regions on which they act as upper-bounds. \\

Propositions \ref{U:cEinf} and \ref{U:cEXinf}
are deduced from the estimates 
given in Lemmas \ref{U:ydP} and \ref{U:N^DP}
on autonomous processes of the form:
\begin{align*}  
	& N^D_t
	:= n +  \int_{0}^{t} 
	(r - c\ N^D_s)\ N^D_s\ ds 
	+\int_{0}^{t}  \sigma \ \sqrt{N^D_s}\ dB_s.
	\EQn{U:YP2}{}
\end{align*}

\noindent
Propositions \ref{U:cEinf} relies on the property of descent from infinity valid for any value of $r$:
\begin{lem}
	\label{U:ydP}
	Let $N^D$ be the solution of \Req{U:YP2}, 
	for some $r \in \R$ and $\cY>0$, 
	with $n$ the initial condition.
	Then, for any $\tDo, \epsilon > 0$
	there exists $n_{\infty} >0$
	such that:
	\begin{align*}
		\textstyle
		\sup_{n>0}
		\PR_{n} (\tDo < \tau^D_{\downarrow} ) 
		\le \epsilon \quad
		\quad \text{ with } \tau^D_{\downarrow} := 
		\inf \left \lbrace s\ge 0, \; N^D_s \le n_{\infty} \right \rbrace.
	\end{align*}
\end{lem}

Proposition \ref{U:cEXinf} relies on the strong negativity on the drift term:
\begin{lem}
	\label{U:N^DP}
	Considering any $\cY, \tDo >0$, with 
	$\ext^D:= \inf \left \lbrace t\ge 0,\, N^D_t = 0 \right \rbrace$:
	\begin{align*}
		&\textstyle 
		\sup_{ n>0} \,\PR_{n} \left(  \tDo < \ext^D \right)
	\underset{r\rightarrow -\infty}{\longrightarrow}	0.
	\end{align*}
	Moreover, for any $\np, \epsilon>0$,
	there exists $n_c$ such that, 
	for any $r$ sufficiently low,
	with $T^D_{\infty}:= \inf \left \lbrace t\ge 0,\, N^D_t \ge n_c \right \rbrace$: 
	$\quad
	\PR_{\np} \left( T^D_{\infty} \le \tDo \right) 
	+ \PR_{\np} \left( N^D_{\tDo} \ge \np \right) 
	\le \epsilon.$
\end{lem}

Finally, Proposition \ref{U:cE0}\  
relies on an upper-bound given 
as a Continuous State Branching Process, 
for which the extinction rate 
is much more explicit.
It is clearly as strong as needed 
for sufficiently small initial condition.

We recall that the complete proofs
(of Proposition \ref{A2XN} from the propositions 
and of the propositions themselves)
are deferred to respectively
Appendices A, B, C and D.
With this,
we have concluded
the proof of Proposition \ref{A2XN}.


\section{Proof of Theorems \ref{U:AllPho}-3}
\label{U:sec:Pf}
\setcounter{eq}{0}

In Subsection~\ref{U:sec:coupling}, we present 
the general principles of our coupling 
that concludes the proof of Theorem~\ref{U:AllPho}.
These principles would alone end the proof
in the context of the Assumption $(A)$ 
in \cite{ChQSD}.
Yet, with our more general assumptions, 
these principles require 
the results of the two previous subsections.
First, we prove in Subsection~\ref{U:sec:Qt} that 
the MCNE will keep in the long run some mass on some specific set $\cD\iXT$ 
(which is weaker but related in some sense to the tension of the laws); 
then we prove in Subsection~\ref{U:sec:pers} that 
\textup{$(A3)$}\, holds in fact for $\cX$ instead of just $E$.
At the end of Subsection~\ref{U:sec:coupling},
the proof of Theorem~\ref{U:AllPho} is then complete. 
The following Subsection \ref{U:sec:eta} and \ref{U:sec:QECv}
then prove respectively Theorem \ref{U:EtaECV} and \ref{U:QECV}.

\subsection{Stabilization of the process in the long run}
\label{U:sec:Qt}

The main purpose of this section is to prove:

\begin{theo}
	\label{U:EtRet}
	Assume that $\mathbf{(A)}$ holds. 
	Then, there exists $\M\iXT = \M_{\ell\iXT,\, \xi\iXT}$ 
	(with $\ell\iXT \ge 1$, $\xi\iXT > 0$) such that for any $\ell \ge 1$ and $\xi \in (0, 1)$, 
	there exists 
	$t\iXT = t\iXT(\ell, \xi)>0$ such that:
	\newcounter{Mxt}
	\begin{align*}
		&
		\frl{\mu \in \M_{\ell,\, \xi}  }
		\frl{\tp \ge t\iXT}\qquad 
		\mu A_\tp \in \M\iXT.
		\EQn{U:Mxt}{\M\iXT}
	\end{align*}
\end{theo}

\begin{rem}
	Assumption \textup{$(A3)$}\, is not involved in the proof
	of Theorem \ref{U:EtRet}.
	This will be important 
	in \cite{AV_disc}
	since we will exploit 
	Theorem \ref{U:EtRet}
	to provide an alternative criterion to $\textup{$(A3)$}$.
\end{rem}

\noindent$  $
\textbf{Proof of Theorem \ref{U:EtRet}}:

According to \textup{$(A2)$}, let 
$\ell\iET \ge 1$ and $\rho\iET > \rho\iSv$ be
such that:
\newcounter{eqDbk}
\begin{align*}
	&\text{with }E:= \cD_{\ell_{\iET}}, \; 
	\tau\iET^1:= \inf\left \lbrace \tp \ge 0;\; X_\tp \in E \right \rbrace, \\
	&\hspace{.5cm}
	\textstyle
	e_\cT:= \sup_{x} \E_{x} \left( \exp\left[ \rho\iET\, (\tau\iET^1\wedge \ext) \right] \right)< \infty. 
	\EQn{U:eqDbk}{E}
\end{align*}
From \Req{U:rsv}, i.e. the definition of $\rho\iSv$,
there exists $\widetilde{\rho}\iSv \in (\rho\iSv, \rho\iET)$,
$c\iSv>0$ and $\ell_S \ge 1$,
such that:
\begin{equation}
	\frl{\tp \ge 0}
	\PR_{\zeta} (\tp < \ext\wedge T_{\cD_{\ell_S}})
	\ge \cp\iSv \; \exp(-  \widetilde{\rho}\iSv \, \tp).
	\EQn{U:RsvI}{}
\end{equation}
We then apply \textup{$(A1)$}\, with  $\ell = \ell_{\iET}$ to state that there exists
$L\iG \ge \ell_S\vee  \ell\iET$,
$t\iG, c\iG > 0$ such that,
with $\cD\iG:= \cD_{L\iG}$:
\begin{align*}
	&
	\frl{x \in E} \quad
	\PR_{x} \Big[ {X}_{t\iG}\in dy,\;
	t\iG < \ext\wedge T_{\cD\iG}\Big]
	\ge c\iG\, \zeta(dy).
	\EQn{U:cDiG}{\cD\iG}
\end{align*}
W.l.o.g., we are allowed to replace,
in the following usage of \Req{U:RsvI},
$\cD_{\ell_S}$ by $\cD\iG$.

In order to conclude the proof of Theorem~\ref{U:EtRet}, we need the following three Lemmas,
for which we define 
by induction over $i\in \N$:
\begin{align*}
	&T^i\iG:= \inf\left \lbrace \tp \ge \tau\iET^i;\;
	X_\tp \notin \cD\iG \right \rbrace,\; 
	T^0\iG:= 0,\\
	&\tau\iET^{i+1}:= \inf\left \lbrace \tp \ge T^i\iG;\; X_\tp \in E \right \rbrace.
	\end{align*}

\begin{lem} \textbf{First entry in $E$: }
	\label{U:T1in}
	\noindent
	Assume that \Req{U:eqDbk}, \textup{$(A1)$}\, and \Req{U:RsvI} hold.
	Then, for any $\ell, \xi$,
	there exists $C_{E} = C_{E}(\ell,\, \xi) > 0$ such that:
	\begin{equation*} 
		\frl{t_{h} > 0} 
		\frl{\mu\in  \M_{\ell,\, \xi}  }\quad
		\PR_\mu ( t_{h} \le \tau\iET^1 \, \big| \,   t_{h} < \ext )
		\le  C_{E}\, e^{-(\rho\iET - \widetilde{\rho}\iSv)\, t_{h}}.
	\end{equation*}
\end{lem}

\begin{lem} \textbf{Containment of the process after $T^i\iG$: }
	\label{U:TjOut}
	
	\noindent
	Suppose that \Req{U:RsvI} and \textup{$(A1)$}\, hold. Then,
	there exists $\ell_\circ \ge L\iG$, 
	and $c_\circ> 0$ such that:
	\begin{align*}
		&\frl{x\in \cD\iG}
		\frl{t > 0}
		\hspace{0.5cm}
		\PR_{x} \left( 
		\tp < T\iG^1\wedge T_{\cD_{\ell\iXT}}\wedge \ext \right) 
		\ge c\iXT\, \exp[- \widetilde{\rho}\iSv\, t].
	\end{align*}
\end{lem}

\begin{lem} \textbf{Last exit from $\cD\iG$:}\\
	\label{U:TteOut}
	Suppose that \textup{$(A0)$}, \textup{$(A1)$}, \Req{U:eqDbk}, \Req{U:RsvI} hold
	(with $E \subset \cD\iG$)
	and $\rho\iET > \widetilde{\rho}\iSv$. 
	Then, there exists $C_{L} > 0$, 
	such that for any $\mu\in  \M_1\left(\cX\right)$ with $t_{h} > \tmo>0$:
	\begin{align*}
		&\hspace{.5cm}
		\PR_\mu \left( T\iG^{I(t_{h})} \le t_{h} - \tmo,\;
		t_{h} \le \tau\iET^{I(t_{h})+1},\; 
		\tau\iET^1 < t_{h}  \, \Big| \,   t_{h} < \ext \right) 
		\le  C_{L} \, e^{-(\rho\iET - \widetilde{\rho}\iSv)\, \tmo},\\
		&\hspace{2.5cm}\text{ with }
		I(t_{h}):= \max\left \lbrace i\ge 0;\; T^i\iG \le t_{h}\right \rbrace 
		(< \infty \text{ a.s.}).
	\end{align*}
\end{lem}

The proofs of these Lemmas are deferred,
in the order of occurrence,
after the proof that they imply Theorem~\ref{U:EtRet}.

%
\subsubsection
[Proof  of \Req{U:Mxt}]
{Proof that Lemmas \ref{U:T1in}-3
	imply Theorem~\ref{U:EtRet}} 

With Lemma \ref{U:T1in} and \ref{U:TteOut}, we obtain an upper-bound 
(with high probability) on how much time the process may have spent outside $\cD\iG$. 
Thus, we can associate most of trajectories ending outside $\cD\iG$
to others ending inside $\cD\iG$. From this association, we deduce a lower-bound on the probability to see the process in $\cD\iG$.\\

Let us first define $\cD\iXT$ according to Lemma \ref{U:TjOut}. In the following, we will define:
\\
$\M\iXT:= \left \lbrace \mu\in \M_1\left(\cX\right);\; \mu(\cD\iXT) \ge \xi\iXT\right \rbrace 
\quad$ for a well-chosen $\xi\iXT.$
Thanks to Lemma \ref{U:TteOut}, we choose some $\tmo> 0$ sufficiently large 
to ensure: $ \frl{t_{h} > \tmo}
\frl{\mu \in \M_1\left(\cX\right)}$
\begin{align*}
	& \PR_\mu \left( T\iG^{I(t_{h})} \le t_{h} - \tmo,\;
	t_{h} \le \tau\iET^{I(t_{h})+1},\; 
	\tau\iET^1 < t_{h}  \, \Big| \,   t_{h} < \ext \right) 
	\le \Ndiv{4}.
	\EQn{LeD}{Le}
\end{align*}
Let $\ell\ge 1,\, \xi \in (0, 1)$.
Thanks to Lemmas \ref{U:T1in}, we know that for some 
\mbox{$t\iXT \ge \tmo > 0$:}
\begin{align*}
	\frl{t_{h} \ge t\iXT} 
	\frl{\mu\in  \M_{\ell,\, \xi}  }\quad
	\PR_\mu \left( t_{h} \le \tau\iET^1 \, \Big| \,   t_{h} < \ext \right) 
	\le  \Ndiv{4}.
	\EQn{Fa_1}{Fa}
\end{align*}
Let $\mu \in \M_{\ell,\, \xi}  $. Let us first assume that:
\begin{align*}
	\PR_\mu \left(  \tau\iET^{I(t_{h})+1}\le t_{h}  
	\, \Big| \,   t_{h} < \ext \right) 
	\ge \Ndiv{4}. 
\EQn{TEtH}{}
\end{align*}
By definition of $I(t_{h})$, on the event $\left \lbrace \tau\iET^{I(t_{h})+1}\le t_{h} \right \rbrace \medcap
\left \lbrace t_{h} < \ext  \right \rbrace$, we know that the process stays in 
$\cD\iG$ in the time-interval $[\tau\iET^{I(t_{h})+1}, t_{h}]$. In particular:
\begin{align*}
	& \mu A_{t_{h}} (\cD\iXT)
	\ge \mu A_{t_{h}} (\cD\iG)
	\ge \PR_\mu \left(  \tau\iET^{I(t_{h})+1}\le t_{h}  
	\, \Big| \,   t_{h} < \ext \right) 
	\ge \Ndiv{4}.
	\EQn{U:maT}{}
\end{align*}
where we recall that $\ell\iXT \ge L\iG$ by Lemma \ref{U:TjOut}.

Now that this case has been easily treated, we consider the complementary:
\begin{align*}
\PR_\mu \left(  \tau\iET^{I(t_{h})+1}\le t_{h}  
\, \Big| \,   t_{h} < \ext \right) 
< \Ndiv{4}.
\end{align*}
Thus, by \Req{LeD}\, and \Req{Fa_1}: 
$\PR_\mu \left( t_{h} - \tmo < T\iG^{I(t_{h})},\;
\tau\iET^1 \le t_{h}  \, \Big| \,   t_{h} < \ext \right) 
\ge \Ndiv{4}.$
By defining the stopping time: 
$\tau\iG:= \inf\{s\ge t_{h} - t;\; X_s \in \cD\iG\},$ we deduce:
\begin{align*}
	&\PR_\mu \left( t_{h} - \tmo \le \tau\iG < t_{h}  \, \Big| \,   t_{h} < \ext \right) 
	\ge \Ndiv{4}.
	\EQn{U:condC}{}
\end{align*}
By the Markov property, then Lemma \ref{U:TjOut}:
\begin{align*}
	&\PR_\mu \left( X_{t_{h}} \in \cD\iXT,\;
	t_{h} - \tmo \le \tau\iG < t_{h} ,\;   t_{h} < \ext \right) 
	\\&\hcm{0.5}
	\ge \E_\mu \left[
	\PR_{X_{\tau\iG}} \left(
	\widetilde{X}_{t_{h}- \tau\iG} \in \cD\iXT,\; 
	t_{h}- \tau\iG < \widetilde{\ext}\right);\;
	t_{h} - \tmo \le \tau\iG < t_{h} \wedge \ext \right]
	\\&\hcm{0.5}
	\ge c\iXT \exp[-\widetilde{\rho}\iSv \tmo]\,
	\PR_\mu \left[
	t_{h} - \tmo \le \tau\iG < t_{h} \wedge \ext \right]
	\\&\hcm{0.5}
	\ge c\iXT \exp[-\widetilde{\rho}\iSv \tmo]\,
	\PR_\mu \left[
	t_{h} - \tmo \le \tau\iG < t_{h} \, \big| \, t_{h} < \ext \right]
	\times \PR_\mu \left[t_{h} < \ext \right].
\end{align*}
So \Req{U:condC} indeed implies
$\mu A_{t_{h}} \left( \cD\iXT\right) \ge \xi\iXT$
with $\xi\iXT:= c\iXT\, e^{-\widetilde{\rho}\iSv\, \tmo} /4$.
With $\M\iXT:= \left \lbrace \mu\in \M_1\left(\cX\right);\; \mu(\cD\iXT) \ge \xi\iXT\right \rbrace$ 
($\xi\iXT$ given by the previous formula 
does not depend on $\ell$, $\xi$ or $\mu$),
we indeed prove \Req{U:Mxt} 
for the case where \Req{TEtH} does hold false.
Recall that the first case where \Req{TEtH} does hold true
 is directly concluded in  \Req{U:maT}.
\hfill $\square$

\subsubsection
[First arrival in $E$]
{Proof of Lemma \ref{U:T1in}} 
%
By \Req{U:eqDbk} and the Markov inequality:
\begin{align*}
	\frl{\mu}\frl{t_{h}>0} \quad
	\PR_\mu \left( t_{h} \le \tau\iET^1\wedge \ext \right) \le e_\cT \; e^{-\rho\iET \, t_{h}}. \quad 
	\EQn{tauE1}{\tau\iET^1}
\end{align*}
Let $\ell\ge 1$, $\xi \in (0, 1)$.
We apply \textup{$(A1)$}, \Req{U:RsvI} and the Markov property to obtain that there exists $c>0$
(that we can decompose as $c=\xi\, c_S\, c_M(\ell)$ 
for some constant $c_M(\ell)>0$ only depending on $\ell$)
 such that:
$\frl{\mu \in \M_{\ell,\, \xi}  }\frl{t_{h}>0}$
\begin{align*}
	&\PR_\mu \left( t_{h} < \ext \right) 
	\ge c\, e^{-\widetilde{\rho}\iSv\, t_{h}}.
	\EQn{U:extR}{}
\end{align*}
Thus, by \Req{tauE1}, \Req{U:extR}, with: 
$ C_{E}:= 
e_\cT / c > 0,$
\begin{align*}
	\frl{\mu \in \M_{\ell,\, \xi}  }\frlq{t_{h}>0}
	\PR_\mu \left( t_{h} \le \tau\iET^1 \, \Big| \,   t_{h} < \ext \right) 
	\le  C_{E} \, \exp[-(\rho\iET - \widetilde{\rho}\iSv)\, t_{h}].
	\tag*{$\square$}
\end{align*}


\subsubsection
[Containment of the process after $T^i\iG$]
{Proof of Lemma \ref{U:TjOut}} 


Thanks to \textup{$(A1)$},
applied with $\ell = L\iG$, 
there exists some $\cD\iXT$, $t\iET, c\iET > 0$ such that:
\begin{align*}
	&\hspace{.5 cm}\frlq{x \in \cD\iG }
	\PR_{x} \left( \tau\iET^1 \le t\iET\wedge T_{\cD\iXT}\wedge \ext \right) \ge c\iET.
	\EQn{ciE}{c\iET}
\end{align*}
Recalling \Req{U:cDiG}, 
we deduce that conditionally on $\mathcal{F}_{\tau\iET^{1}}$
and on the event 
$\left \lbrace \tau\iET^1 \le t\iET\wedge T_{\cD\iXT} \right \rbrace$:
\begin{align*}
	\PR_{X_{\tau\iET^{1}} } \left( 
	\widetilde{X}_{t\iG}\in dx,\;
	t\iG < \widetilde{T\iG^1} \wedge \widetilde{T}_{\cD\iXT} 
	\wedge \widetilde{\ext} \right) 
	\ge c\iG\, \zeta(dx).  
	\EQn{ciG}{c\iG}
\end{align*}
By combining \Req{ciE}, \Req{ciG}, 
\Req{U:RsvI} and the Markov property, we deduce that:
\begin{align*}
	& \PR_{x} \left( 
	t_{h} < T\iG^1\wedge T_{\cD\iXT}\wedge \ext\right) 
	\ge c\iXT\, \exp[- \widetilde{\rho}\iSv\, t_{h}],\\
	&\hspace{1.5cm}
	\text{ with } c\iXT:= c\iET\, c\iG\, c\iSv\, 
	\exp[ \widetilde{\rho}\iSv\, (t\iET + t\iG) ]>0. \tag*{$\square$}
\end{align*}

\subsubsection
[Last exit from $\cD\iG$]
{Proof of Lemma \ref{U:TteOut}} 
\label{U:sec:Tteout}

The idea 
is to use that
it is very unlikely
for the process to still be alive 
after experiencing an excursion outside $E$ 
for a long time (and still be there). 
Indeed, 
compared to trajectories 
that stay inside $E$
(in particular those reaching quickly $\cD\iSv$,
for which \Req{U:Rg} holds, and not leaving $\cD_L$)
the probabilities 
of the associated events
vanish with a larger rate: 
$\rho\iET > \widetilde{\rho}\iSv$.
It would have been convenient if we could
 initiate the comparison
just before $T^{I(t_{h})}\iG$,
where the process 
exits $E$ for the last time  before $t_{h}$.
Yet, it is not a stopping time,
so that the Markov property is not directly applicable
and the proof gets more technical.

Let us first prove that $I(t_{h}) < \infty$. 
Since $X$ has càdlàg paths, we would have on the event $\left \lbrace I(t_{h}) = \infty \right \rbrace$:
$\sup_{j} T\iG^i = \sup_{j} \tau\iET^i = T < t_{h}$
with $X_{T-} \in E 
\cap \overline{\cX \setminus \cD\iG}.$
Yet, by \textup{$(A0)$}, this set is empty, so that a.s. $I(t_{h}) < \infty$.
Then, exploiting a discretization of time 
in time-intervals of length $t_L$ to be fixed later:
\begin{align*}
	&P:=
	\PR_\mu \left( T\iG^{I(t_{h})} \le t_{h} - \tmo,\;
	t_{h} \le \tau\iET^{I(t_{h})+1}\wedge \ext,\; 
	\tau\iET^1 < t_{h}  \right) 
	\\&\hspace{1.5cm} 
	= \Tsum{i\ge 1} 
	\PR_\mu \left( T\iG^{i} \le t_{h} - \tmo,\;
	t_{h} \le \tau\iET^{i+1}\wedge \ext \right) 
	\\&\hspace{1.5cm} 
	\le \Tsum{i\ge 1} \Tsum{k\ge 0} 
	\idc{k\, t_L \le t_{h} - \tmo }
	\PR_\mu \left( 
	T\iG^{i} \in (k\, t_L, (k+1)t_L],\;
	t_{h} \le \tau\iET^{i+1}\wedge \ext\right) 
	\\&\hspace{1.5cm} 
	= \Tsum{i\ge 1} \Tsum{k\ge 0} 
	\idc{k\, t_L \le t_{h} - \tmo }
	\E_\mu \big[
	\PR_{X_{(k+1)t_L} } \left( 
	t_{h} - (k+1)t_L  
	\le \widetilde{\tau}\iET^{1}\wedge \widetilde{\ext} \right)
	;\; \\&\hspace{5.5cm} 
	T\iG^i \in (k\, t_L, (k+1)t_L],\; 
	(k+1)t_L \le \tau\iET^{i+1}\wedge \ext \big],
\end{align*}
where we used the Markov property.
Exploiting \Req{U:eqDbk}:
\begin{align*}
	P 
	&\le e_T \Tsum{k\ge 0} 
	\idc{k\, t_L \le t_{h} - \tmo }
	\exp[-\rho\iET\, (t_{h} -(k+1) t_L)]\;
	\\&\hcm{2}
	\times\Tsum{i\ge 1} \PR_\mu \left[ 
	T\iG^i \in (k\, t_L, (k+1)t_L],\; 
	(k+1)t_L \le \tau\iET^{i+1}\wedge \ext \right].
	\EQn{U:decTC}{}
\end{align*}

The trick is to observe that, 
by definitions of $\tau\iET^{i} < T\iG^i$,
one shall have $X_s\in \cD\iG$ 
for any $s\in [\tau\iET^{i}, T\iG^i)$, 
in particular on some vicinity 
to the left of $T\iG^i$.
Defining for $k\ge 0$:
\begin{equation*} 
	\tau\iG^k  
	:= \inf\left \lbrace s\ge k t_L: X_s \in \cD\iG\right \rbrace,
\end{equation*}
we see that the events 
$\left \lbrace T\iG^i \in (k\, t_L, (k+1)t_L]\right \rbrace
\cap\left \lbrace (k+1)t_L \le \tau\iET^{i+1}\wedge \ext\right \rbrace$ are disjoint (for $k$ fixed)
and included in the event 
$ \left \lbrace \tau\iG^k < (k+1)t_L \wedge \ext\right \rbrace$.
On the other hand, 
exploiting the Markov property 
together with Lemma \ref{U:TjOut}:
\begin{equation*} 
	\PR_\mu \left[ 
	t_{h} < \ext\right]
	\ge c\iXT \, \exp[-\widetilde{\rho}\iSv (t_{h} - k\, t_L)]
	\, \PR_\mu \left[ 
	\tau\iG^k < (k+1)t_L \wedge \ext \right].
\end{equation*}
Coming back to \Req{U:decTC},
we deduce:
\begin{align*}
	P
	&\le \dfrac{e_T\, e^{\rho\iET\,t_L}}{c\iXT}
	\PR_\mu \left[ 
	t_{h} < \ext\right]
	\times \Tsum{k\ge 0} 
	\idc{k\, t_L \le t_{h} - \tmo }
	\exp[-(\rho\iET-\widetilde{\rho}\iSv)
	\times (t_{h} -k t_L)].
\end{align*}
The sum over $k$ is upper-bounded by:
\begin{align*}
	\exp[-(\rho\iET-\widetilde{\rho}\iSv) t]\times 
	\Tsum{\ell \ge 0} 
	\exp[-\ell (\rho\iET-\widetilde{\rho}\iSv) t_L]
	\le \dfrac{e^{-(\rho\iET-\widetilde{\rho}\iSv) t}}{
		1- e^{-(\rho\iET-\widetilde{\rho}\iSv) t_L}}.
\end{align*}
This concludes the proof of the Lemma, 
with: 
$C_L:= \dfrac{e_T\,e^{\rho\iET\,t_L} } {c\iXT (1- e^{-(\rho\iET-\widetilde{\rho}\iSv) t_L})}.$
The choice of $t_L$ is free, so that we can fix it to optimize this constant.
\hfill $\square$

\subsection{Persistence}
\label{U:sec:pers}

\subsubsection{Theorem~\ref{U:PersW}}

For the proof of the following Theorem~\ref{U:PersW},
we need the following Corollary of Theorem~\ref{U:EtRet}:
\begin{cor} \textbf{"Stability"}:\\
	\label{U:Sb}
	Under Assumption $\mathbf{(A)}$, 
	there exists $t\iSB, c\iSB' >0$ 
	and $ \widetilde{\rho}\iSv \in (\rho\iSv, \rho\iET)$
	such that: 
	\newcounter{Sb}
	\begin{align*}
		&\frl{ u\ge 0}
		\frl {\tp \ge u + t\iSB}
		\hspace{0.5cm}
		\PR_{\zeta}(\tp - u < \ext) 
		\le c\iSB'\; e^{\widetilde{\rho}\iSv u}\;  \PR_{\zeta}\left(  \tp < \ext \right).
		\EQn{U:eSb}{Sb}
	\end{align*}
\end{cor}

\begin{theo}
	\label{U:PersW}
	Assume that there exists $\rho\iET > \widetilde{\rho\iSv}$, $\cD\iSv \subset \cX$, $E \subset \cX$ and $\zeta \in \M_1\left(\cX\right)$ such that  
	\textup{$(A3)$}, \textup{$(A2)$}, \Req{U:RsvI}  and 
	\Req{U:eSb}\, hold.
	Then, there exists $t\iPs, c\iPs>0$ such that:
	\begin{align*}
		&\frl{x\in \cX}
		\frl{\tp \ge t\iPs}
		\quad
		\PR_{x} \left( \tp < \ext\right)  
		\le c\iPs\; \PR_{\zeta} \left( \tp < \ext\right).
		\EQn{U:eqPS}{}
	\end{align*}
\end{theo}

\subsubsection
[Proof of \Req{U:eSb}:]
{Proof of Corollary \ref{U:Sb}:}

By \Req{U:Mxt} and \textup{$(A1)$}, 
there exists $c>0$ such that 
for any $v$ sufficiently large:
$\zeta A_v \ge c\, \zeta$. 
with the Markov property, it implies 
for any $u\ge 0$:
\begin{equation*} 
	\PR_{\zeta}(v+u< \ext)
	\ge c\, \PR_{\zeta}(v< \ext)\, \PR_{\zeta}(u< \ext).
\end{equation*}
Exploiting \Req{U:RsvI} with $t=u$,
we deduce Corollary \ref{U:Sb}
with $v = t-u$ and $c\iSB' = (c\,c\iSv)^{-1}$.

\hfill $\square$

\subsubsection{Proof of Theorem~\ref{U:PersW}}
\label{U:persw}
From \textup{$(A3)$}, there exists $t_{A}, c_{A}>0$ such that:
\begin{equation}
	\frl{ \tp \ge t_{A}}
	\frl{ x\in E} \quad  
	\PR_{x} (\tp< \ext) 
	\le c_{A} \;\PR_{\zeta} (\tp < \ext).
	\EQn{U:AS}{}
\end{equation}

This proof is very close to the one in \cite{ChpLyap} 
(p13:"Step 2: Proof of \textup{$(A1)$}"), 
except that, in \Req{U:eSb},
$t-u$ shall be larger than some value,
and similarly for $t$ in \Req{U:AS}.
To compare the notations, our $e_\cT$, $c_{A}$ and $c\iSB$ refer resp.
to their $M$, $C_m$ and $4/ c_1\, D_m\, D_{n_1}$.
Thus, we won't detail it much and refer to \cite{ChpLyap}.

Let $\zeta \in \M_1(\cX)$, 
$\tp \ge t\iPs:= t\iSB \vee t_{A}$ and $x\in \cX$.
\begin{align*}
	&\PR_x(\tp < \ext) 
	\le c_{A} \, \E_x\left[ \PR_{\zeta} (t- \tau_{E} < \widetilde{\ext});\;
	\tau_{E} < (t-t\iPs)\wedge \ext \right]
	+ \PR_x (t-t\iPs \le \tau_{E}\wedge \ext) 
	\EQn{A3_1}{\tau_{E}< t- t\iPs}
\end{align*}
thanks to property \textup{$(A3)$}, since $t-\tau_{E}\ge t\iPs \ge t_{A}$ on 
$\left \lbrace  \tau_{E} < (t-t\iPs)\wedge \ext\right \rbrace$. 
By \textup{$(A2)$}\, (with the Markov inequality) and Corollary \ref{U:Sb}, 
with $u = t_{A}$ for the first term of \Req{A3_1}\, and 
$u = \tp - t\iSB$ for the second:$\quad\forall \tp \ge 0,\, \forall x\in \cX,$
\begin{align*}
	\PR_x(\tp < \ext) 
	\le \left( c_{A} + e^{\widetilde{\rho}\iSv (t\iPs- t\iSB)}\,
	/ \, \PR_{\zeta}\left(  t\iSB < \ext \right) \right) \times 
	c'\iSB \,e_\cT\;
	\, \PR_{\zeta}(\tp < \ext)
	\tag*{$\square$}
\end{align*}

\subsection{Coupling procedure: proof of Theorem \ref{U:AllPho}}
\label{U:sec:coupling}

\subsubsection{Definition of the uncoupled part} 
\label{U:sec:amu}
\newcounter{Mrn}

With a given set of parameters $t\iDB,\, c\iDB,\, t\iPs,\, c\iPs>0$ 
(cf following subsection) 
we define for $t_h > t\iPs$:
\begin{equation*} 
	J(t_h)
	:= \left \lfloor (t_h - t\iPs)/t\iDB \right \rfloor.
	\EQn{U:Jtobs}{J(t_h)}
\end{equation*}
For $\tp \ge 0$, $\mu \in \M_1\left(\cX\right)$,
$t_h > t\iPs$, and $k\in \N$, 
let:\newcounter{AMU}\newcounter{weta}
\begin{align*}
	&\amu(k, t) = a_{\mu}^{t_h}(k, t)
	:= \idc{k\le J(t_h),\, k\, t\iDB \le \tp}
	\times \text{\large{$\,^{c\iDB}\!/_{c\iPs}$} }\, ( 1-\text{\large{$\,^{c\iDB}\!/_{c\iPs}$} } )^{k-1}
	\\&\hspace{4cm}
	\times \dfrac{\PR_\mu (t_h < \ext)}{\PR_\mu (\tp < \ext)}
	\times \dfrac{\PR_{\zeta} (\tp -k\, t\iDB< \ext)}
	{\PR_{\zeta} (t_h -k\, t\iDB< \ext)}.
	\EQn{U:AMU}{\amu(k, t)}
\end{align*}
\begin{rem}
	As we can see in the proof of Fact \ref{U:ser}, 
	$\amu(k, t)$ corresponds to 
	the mass associated with the $k$-th step of coupling, 
	considered at time $t$ 
	with the constraint that it must represent 
	a fixed proportion of $\mu A_{t_h}$ (at time $t_h$).
	We refer to Figure 2
	for a presentation of the coupling architecture.
\end{rem} 

Let $\quad r_j:= 1 -  \Tsum{k \le j} \amu(k, j\, t\iDB).\quad$
Under the condition $r_j >0,$
that we will prove to be true by induction over 
\mbox{$j\le J(t_h)$}, we define:
\begin{equation*} 
	\nu_j(dx)
	:= (^1/r_j) \ltm 
	\left[ \mu A_{j\, t\iDB}(dx) - 
	\Tsum{k \le j} 
	\amu(k, j\, t\iDB)\, \zeta A_{(j-k)\, t\iDB}(dx) \right],
	\EQn{U:nuJ}{\nu_j}
\end{equation*}
\noindent
with the convention $\nu_0:= \mu$.
Remark that this definition ensures $\nu_j(\cX) = 1.$
\begin{rem}
	$\nu_j$ shall correspond to the marginal of the process conditioned 
	of not being already coupled at time $j\, t\iDB$. 
	We normalize what remains of $\mu A_{j\, t\iDB}$ 
	when we subtract the contribution of each coupling step 
	(only those up to the $j$-th will contribute to the sum). 
	The main difficulty will be to prove that, 
	under suitable conditions, 
	$\nu_j$ is indeed a positive measure, 
	thus a probability measure.
	In Figure 2,
	the associated coupling procedure 
	is presented for the more general case 
	where we compare two initial conditions in some $\M_{\ell,\, \xi}  $ 
	rather than already in $\M\iRN$.
\end{rem}

\begin{rem}
The procedure extends to cases 
where the QSD is not unique
provided that
$\mu^{(1)}$ and $\mu^{(2)}$ 
are in the same basin of attraction,
as we can see in \cite{AV_GS}.
This procedure already deals 
with an inhomogeneity in time 
due to the conditioning for survival at time $t_h$,
so that adaptations
to time-inhomogeneous processes
are likely to be easy.
\end{rem}

\begin{figure}[!h]
	\rule{\textwidth}{0.4pt}
	\begin{center}
		\includegraphics[width = 14cm, height = 5.5cm]{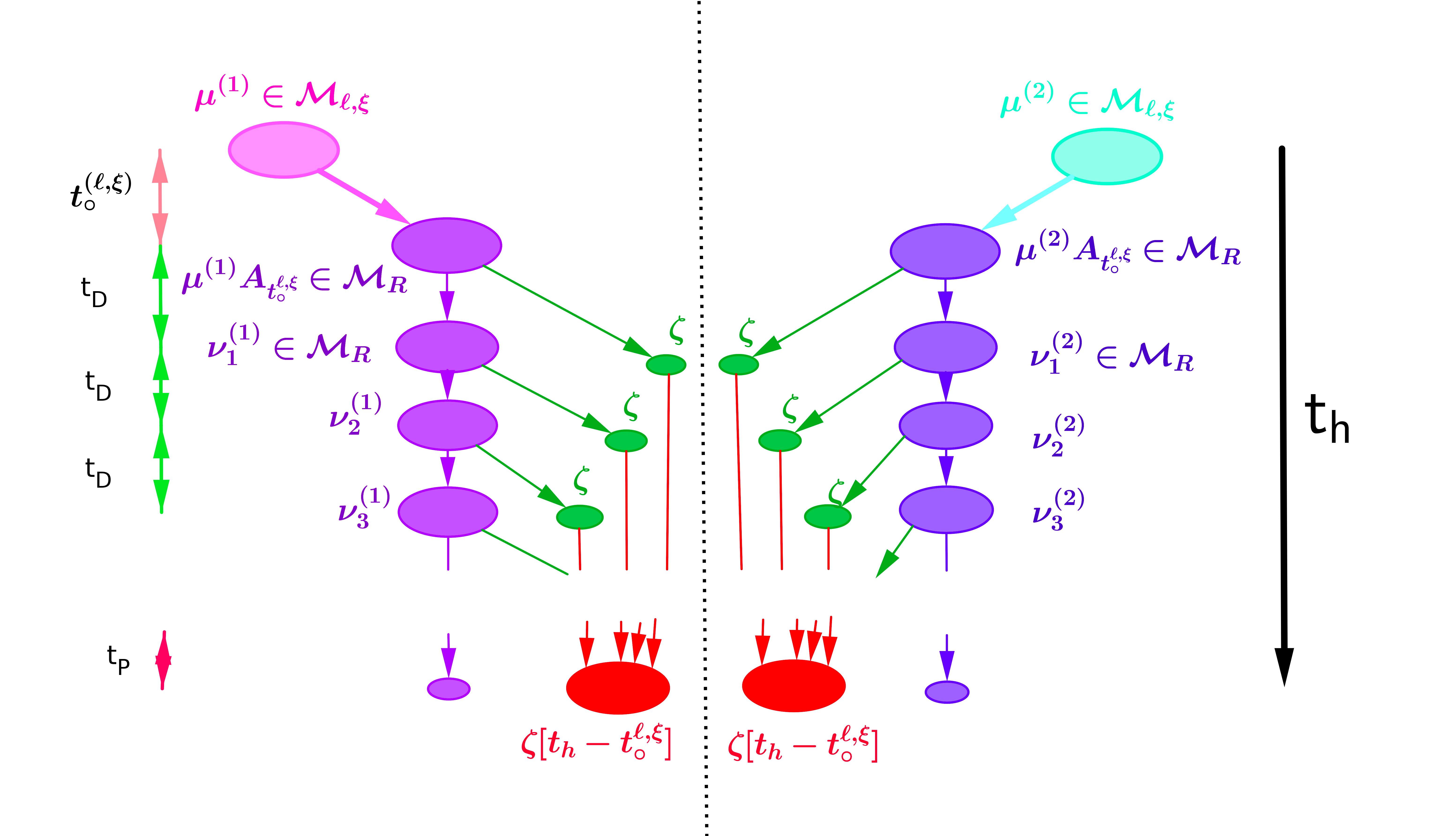}
		\label{U:cpl}
		\caption{Illustration of the coupling procedure}
	\end{center}
	The figure illustrates the coupling procedure 
	on two initial conditions $\mu^{(1)}$ and $\mu^{(2)}$.
	We can observe by symmetry 
	how the MCNE are progressively decomposed
	with time descending along the vertical axis.
	By construction, 
	the middle red parts (at time $t_h$) 
	are common for both initial conditions
	(both its distribution $\zeta[t_h - t\iXT^{\ell, \xi}]$
	and the amount of mass).	
	\rule{\textwidth}{0.4pt}
\end{figure}

\noindent
\textbf{Explanation of the procedure as presented in Figure \ref{U:cpl}: }
Since both initial conditions 
belong to the same $\M_{\ell,\, \xi}  $, 
the time $t\iXT = t\iXT^{\ell, \xi} $ 
needed to reach
$\M\iRN$ can be chosen uniformly.
Then, 
after every time-interval of size $t\iDB$ 
and as long as time
$t_h - t\iPs$ is not reached,
we shall exploit 
property \Req{U:Mrn}.
We split the "remaining MCNE"
(the $\nu_j$ after $j$ splitting steps)
in order to extract a component 
whose contribution to the MCNE at time $t_h$
is explicit. 
This contribution 
(in the expression 
of $\zeta[t_h - t\iXT^{\ell, \xi}]$)
is proportional 
to $\zeta A_{t_h - t}$ 
for a splitting at time $t$.
For this contribution at time $t_h$
to be fixed, 
note that the contribution to the MCNE
at time $t$
has to depend
both on the remaining time $t_h - t$
and the specific value of $\nu_j$.

\subsubsection{Definition of the constants involved}

For clarity, we denote by $t_h$ (for horizon of time) the time $t$ that appears in Theorem~\ref{U:AllPho}. During this coupling procedure, it will stay fixed, and won't appear in the other sections. 
The constants $c\iPs, t\iPs>0$ come from Theorem \ref{U:PersW},
while $c\iDB, t\iDB>0$ come from this corollary of Theorem \ref{U:EtRet}:
\begin{prop} \textbf{"Coupling and Renewal"}\\
	\label{U:c1}
	Suppose that \textup{$(A1)$}\,holds
	and \Req{U:Mxt}\, also for some 
	$\M\iXT:= \M_{\ell\iXT,\, \xi\iXT}$.
	Then, with $\ell\iRN:= \ell\iXT$,\,
	\mbox{$\xi\iRN:= \xi\iXT / 2$}, $\quad \M\iRN:= \M_{\ell\iRN,\, \xi\iRN}$,
	$ \quad \cD\iRN:= \cD_{\ell\iRN}$,
	there exits $c\iDB \in (0,1)$ and \mbox{$t\iDB \ge t\iXT(\ell\iRN,\, \xi\iRN)$} such that: 
	\begin{equation*} 
		\frl{\mu \in \M\iRN }\hspace{.5cm}
		\mu A_{t\iDB}(dx)\ge c\iDB\; \zeta(dx)
		\text{ and } \dfrac{\mu A_{t\iDB}(\cD\iRN) - c\iDB }{1-c\iDB} \ge \xi\iRN.
		\EQn{U:Mrn}{\M\iRN}
	\end{equation*}
\end{prop}

\noindent
\begin{rem}
	The subscript $D$ refers to "Doeblin's" condition, 
	since we will likewise iteratively couple a proportion at most $c\iDB$
	of the distribution.
	The properties \Req{U:Mrn} and \Req{U:eqPS}
	make us able to prove the induction:
	$\nu_j \in \M\iRN \Rightarrow \nu_{j+1} \in \M\iRN$.
\end{rem}

\paragraph
[\Req{U:Mxt}\, + \textup{$(A1)$}\, imply (\ref{U:Mrn})]
{Proof of Proposition \ref{U:c1}}
\label{U:sec:mrn}
\textcolor{white}{:}

We apply \Req{U:Mxt}\, with $\ell = \ell\iRN$ 
and $\xi = \xi\iRN$.
Thus, with $t\iRN:= t\iXT(\ell\iRN,\, \xi\iRN)$:
\begin{align*}
	&\hspace{.5cm}
	\frl{\mu \in \M\iRN}
	\frl{\tp \ge t\iRN}\qquad 
	\mu A_\tp \in \M\iXT,
	\; \text{i.e. } \mu A_\tp (\cD\iRN) \ge 2\, \xi\iRN
	\EQn{xiRN}{\xi\iRN}
\end{align*}
We can then define 
$c\iMix\in (0,\, 1)$, $t\iMix \ge t\iRN$
thanks to \textup{$(A1)$}, 
cf Subsection~\ref{U:mix2},
such that: 
\begin{align*}
	&\frl{x\in \cD\iRN}
	\PR_x \left[ X_{t\iMix} \in dx;\; t\iMix < \ext \right] 
	\ge  c\iMix\; \zeta(dx).
	\EQn{ciRN}{Mix:\cD\iRN}
\end{align*}
We can then choose $c\iDB:= c\iMix\; \xi\iRN\in (0,1),\;
t\iDB:= t\iMix\ge t\iRN$
(for the statement of the proposition), 
and observe that:
\begin{align*}
	\frl{\mu \in \M\iRN }
	\mu A_{t\iMix}(dx) 
&	\ge \mu(\cD\iRN)\,  c\iMix \, \zeta(dx) 
	\hcm{1} \text{ by } \Req{ciRN}\,\\
	\hspace{3.6cm}
&	\ge  \xi\iRN\, c\iMix \, \zeta(dx)
	= c\iDB\; \zeta(dx)
	\qquad \text{ because } \mu \in \M\iRN.\\
	\dfrac{\mu A_{t\iDB}(\cD\iRN) - c\iDB }{1-c\iDB}
	&\ge  \dfrac{2\xi\iRN - c\iDB }{1-c\iDB}
	= \dfrac{2-c\iMix}{1-\xi\iRN \,c\iMix} \xi\iRN \ge \xi\iRN,
\end{align*}
by  \Req{xiRN}, where we used $1\ge (1-\xi\iRN)\, c\iMix$  (of course $c\iMix\in (0,1)$ and $\xi\iRN >0$).
\hfill $\square$
\\

\noindent
\textbf{Remark:}\textsl{
	Our choice for $\xi\iRN$ and $c_D$ is done for simplicity
	and can certainly be improved regarding the convergence rate $\gamma$.
	What we require is rather: 
	$\xi\iRN \le \frac{c\iMix}{c_D}\wedge \frac{\xi\iXT - c_D}{1-c_D}$.
}

\subsubsection{Lower-bound on the marginals}
\label{U:sec:lb}
At time $t_h$, 
for any initial condition $\mu \in \M\iRN$,
the MCNE shall be lower-bounded by:
\begin{align*}
	\zeta[t_h](dx) 
	&
	:= \Tsum{j \le J(t_h)} (\text{\large{$\,^{c\iDB}\!/_{c\iPs}$} })\ltm ( 1-\text{\large{$\,^{c\iDB}\!/_{c\iPs}$} } )^{j-1}
	\, \zeta A_{t_h -j\, t\iDB}(dx) \ge 0.
	\EQn{U:AlcO}{\zeta[t_h] }
\end{align*}
\noindent
\textbf{Remark:}\textsl{
The definition of $(\zeta[t])_{\tp \ge 0}$ implicitly depends on $c\iDB$, $t\iDB$, $c\iPs$ and $t\iPs$, but not on $\mu$, $\ell$ or $\xi$.}
\\

The proof of Theorem \ref{U:AllPho} will be completed thanks 
to Theorems \ref{U:EtRet}, \ref{U:PersW} 
and the following proposition:
\begin{prop}
	\label{U:LowerBound}
	Suppose \Req{U:Mxt}, \textup{$(A1)$}\, and \Req{U:eqPS} hold, with $c\iDB$ and $t\iDB$ chosen according to Proposition \ref{U:c1}, $c\iPs$, $t\iPs$ according to \Req{U:eqPS}. 
	Then, to any pair $\ell \ge 0$ and $\xi \in (0, 1)$, 
	we can associate a time
	\mbox{$t\iXT = t\iXT(n, \xi) >0$} such that:
	\begin{align*}
		\frl{\mu\in  \M_{\ell,\, \xi}  }
		\frl{t_{h, 2}\ge t_{h, 1} \ge t\iXT}
		\hspace{.5cm}
		\mu A_{t_{h, 2}} \ge \zeta[t_{h, 1} - t\iXT].
		\EQn{U:th12}{}
	\end{align*}
\end{prop}

But as  a first step to conclude the proof of Theorem~\ref{U:EtRet}, we prove the following lemma:
\begin{lem}
	\label{U:dec}
	Assume that for some $j\le J(t_h) -1$,
	$r_j > 0$, $\nu_j \in \M_1\left(\cX\right)$ and \Req{U:PSnuj} holds.
	Then:
	\begin{align*} 
		r_{j+1} > 0 \quad \text{ and }
		\Exq{0< c_{j} \le c\iDB}
		\nu_{j+1}(dx) 
		=  \left( \nu_j A_{t\iDB}(dx) - c_{j}\ \zeta(dx) \right) \, /\, (1- c_{j}).
		\end{align*}
\end{lem}
Then, in order to achieve the induction 
"$\nu_j \in \M\iRN$ implies $\nu_{j+1} \in \M\iRN$" 
we ensure iteratively:
\begin{align*}
	\PR_{\nu_j} \left( t_h - j\, t\iDB < \ext \right) 
	&\le c\iPs \,\PR_{\nu_j} \left( t\iDB < \ext \right) \, 
	\PR_{\zeta} \left( t_h - [j+1]\, t\iDB < \ext \right),
	\EQn{U:PSnuj}{}
	\\
	\text{and } 
	\nu_j A_{t\iDB} 
	&\ge c\iDB\, \zeta,
	\EQn{U:DBnuj}{}
	\end{align*}
that come respectively from \Req{U:eqPS} and  \Req{U:Mrn}.
The proof of Proposition \ref{U:LowerBound}
is achieved in this second step, 
while the third one
concludes the proof of Theorem \ref{U:AllPho}.

\paragraph{Step 1: proof of Lemma \ref{U:dec}}
\textcolor{white}{:}

First of all, 
we need to relate 
$1 - \sum_{k \ge 1} \amu(k, j\, t\iDB)$ 
to the repartition of mass at time $t_h$, 
which is done in the proof 
of the following lemma, 
whose proof 
(similar to the next paragraph, yet much simpler)
is reported in Appendix E:
\begin{fact}
	\label{U:ser}
	Assume that for some $j\le J(t_h) -1$: 
	$\quad r_j > 0, \ \nu_j \in \M_1\left(\cX\right)$.
	\begin{align*}
		& \text{Then } 
		r_j
		= \left[ 1- \dfrac{c\iDB}{c\iPs} \right]^{j} \times
		\dfrac{\PR_\mu (t_h < \ext)}{\PR_\mu (j\, t\iDB< \ext)}
		\times \dfrac{1}  {\PR_{\nu_j} (t_h -j\, t\iDB< \ext)}.
		\end{align*}
\end{fact}

By the definition of $\nu_j$, cf \Req{U:nuJ}:
$\mu A_{j\, t\iDB} =
\sum_{k = 1}^j \amu(k, j\, t\iDB)\, \zeta A_{(j-k)\, t\iDB}
+r_j \, \nu_j$.

\begin{align*}
	&\mu A_{[j+1]\, t\iDB} =  \dfrac{\PR_\mu (j\, t\iDB< \ext)}{\PR_\mu ([j+1]\, t\iDB< \ext)}
	\; \mu A_{j\, t\iDB} \cdot P_{t\iDB}\\
	&
	= \Tsum{k \le j} \amu(k, j\, t\iDB)\times 
	\dfrac{\PR_\mu (j\, t\iDB< \ext)
		\times \PR_{\zeta} ([j+1 -k]\, t\iDB< \ext)}
	{\PR_\mu ([j+1]\, t\iDB< \ext)
		\times \PR_{\zeta} ( [j -k]\, t\iDB< \ext)} 
	\zeta A_{[j+1 -k]\, t\iDB}\\
	&\hspace{1.5cm}
	+ \ell_j \;
	\nu_j A_{t\iDB},
	\EQn{mulJ}{}\\
	&\hspace{1cm}
	\text{where } \ell_j:= r_j
	\times \dfrac{\PR_\mu (j\, t\iDB< \ext)}{\PR_\mu( [j+1]\, t\iDB< \ext)} 
	\times \PR_{\nu_j} (t\iDB< \ext).
	\EQn{U:ljD}{\ell_j}
	\end{align*}
By \Req{U:AMU}, 
i.e. the definition of 
$\amu(k, j\, t\iDB)$:
\begin{align*}
	&\amu(k, j\, t\iDB)\times 
	\dfrac{\PR_\mu (j\, t\iDB< \ext)}{\PR_\mu ([j+1]\, t\iDB< \ext)} 
	\times \dfrac{\PR_{\zeta} ([j+1 -k]\, t\iDB< \ext)}
	{\PR_{\zeta} ( [j -k]\, t\iDB< \ext)}
	= \amu(k, [j+1]\, t\iDB).
	\\&\text{Thus }
	\hspace{.5cm} 
	1 = \Tsum{k \le j}  \; \amu(k, [j+1]\, t\iDB) + \ell_j
	\quad\text{i.e.}\quad
	r_{j+1} 
	= \ell_j - \amu(j+1, [j+1]\, t\iDB),
	\EQn{rjDef}{r_{j+1}}
	\end{align*}
by evaluating \Req{mulJ}\, on $\cX$  and the definition of $r_{j+1}$,
cf.\Req{U:nuJ}.

By \Req{U:ljD}, \Req{U:AMU}\, 
and by Fact \ref{U:ser}:
\begin{align*}
	c_{j}:=& \amu(j+1, [j+1]\, t\iDB) / \ell_j
	\EQn{U:cjD}{c_{j}}
	\\&
	= \left[ 1- \dfrac{c\iDB}{c\iPs} \right]^{-j} 
	\times \dfrac{\PR_\mu (j\, t\iDB< \ext)}{\PR_\mu (t_h < \ext)}
	\times \dfrac{\PR_{\nu_j} (t_h - j\, t\iDB< \ext)}{\PR_{\nu_j} (t\iDB< \ext)}  
	\times \dfrac{\PR_\mu ([j+1]\, t\iDB< \ext)}{\PR_\mu (j\, t\iDB< \ext)} \\
	&\hspace{.4cm} 
	\times \dfrac{c\iDB}{c\iPs}\, \left( 1-\dfrac{c\iDB}{c\iPs} \right)^{j} \times 
	\dfrac{\PR_\mu (t_h < \ext)}{\PR_\mu ([j+1]\, t\iDB< \ext)}
	\times \dfrac{1}
	{\PR_{\zeta} (t_h -[j+1]\, t\iDB< \ext)}\\
	&\hspace{0.5cm}  
	= \dfrac{c\iDB}{c\iPs}\,  \dfrac{\PR_{\nu_j} (t_h -j\, t\iDB< \ext)}
	{\PR_{\nu_j} (t\iDB< \ext)\times\PR_{\zeta} (t_h -[j+1]\, t\iDB< \ext)}.
	\end{align*}
Thanks to $\Req{U:PSnuj}$:
$\quad \bullet\quad 0< c_{j} \le  c\iDB. $

Since $c\iDB < 1$, using \Req{rjDef}\,and \Req{U:cjD}: 
$\quad
\bullet\quad  r_{j+1} 
= \ell_j\, \left( 1- c_{j}\right) > 0.$

Finally, by \Req{U:nuJ}, i.e. the definition of $\nu_{j+1}$, \Req{U:cjD}\, and \Req{mulJ}:
\begin{align*}
	&\nu_{j+1} = (1/r_{j+1})\times \left[ \mu A_{[j+1]\, t\iDB} 
	- \Tsum{k \le j+1} \amu(k, (j+1)\, t\iDB)\, \zeta A_{(j+1-k)\, t\iDB} \right]
	\\&\hspace{2.4cm} 
	=\left( \nu_j A_{t\iDB} - c_{j}\; \zeta\right)
	\times  \ell_j \, / r_{j+1}
	\\&\hspace{1cm}
	\bullet\quad \nu_{j+1}
	=  \left( \nu_j A_{t\iDB} - c_{j}\; \zeta\right) / (1-c_{j}).
	\tag*{$\square$}
	\end{align*}

\paragraph
{Step 2: proof of Proposition \ref{U:LowerBound} 
	with Lemma \ref{U:dec}}
\textcolor{white}{:}

We first define $\M\iRN$
thanks to Proposition~\ref{U:c1} together with \textup{$(A1)$}\, and \Req{U:Mxt}\, such that \Req{U:Mrn} holds.

\subparagraph{Step 2.1:}\!under the assumption that $\mu \in \M\iRN$.
Then, by induction over $j\le J(t_h)$, we state 
$(I_j)$:~$r_{j} > 0\text{ and }\nu_j \in \M\iRN.\ $
We initialize at $j =0$, with $r_0 = 1$ and $\nu_0:= \mu \in \M\iRN$ by hypothesis.

Assume $(I_j)$ for some $j\le J(t_h) -1$. 
Then, by $(I_j)$ and \Req{U:Mrn}, $\Req{U:DBnuj}$ holds.
$j \le J(t_h) - 1$ means notably $t_h - [j+1]\, t\iDB \ge t\iPs,$ thus thanks to \Req{U:eqPS}:
\begin{align*}
		& \PR_{\nu_j} \left( t_h - j\, t\iDB < \ext \right)
		= \E_{\nu_j} \left[  \PR_{X_{t\iDB}} \left( t_h - [j+1]\, t\iDB < \ext \right);\; t\iDB < \ext \right]\\
		&\hspace{1cm}
		\le c\iPs\;\PR_{\zeta} \left( t_h - [j+1]\, t\iDB < \ext \right)
		\times \PR_{\nu_j} \left[ t\iDB < \ext \right].
\end{align*}
Thanks to Lemma \ref{U:dec} together with \Req{U:DBnuj}: $\nu_j \ge 0$ 
thus $\nu_j \in \M_1\left(\cX\right)$. 
Moreover, for any  measurable set $\cD$:
\begin{align*}
		\hspace{2cm}
		\nu_{j+1}(\cD) 
		&\ge  \left( \nu_j A_{t\iDB}(\cD) - c_{j} \right) \, / (1- c_{j})
		= 1 - (1 - \nu_j A_{t\iDB}(\cD))\, /\, (1- c_{j})
		\\
		\hspace{3.7cm} 
		&\ge  \left( \nu_j A_{t\iDB}(\cD) - c\iDB \right)\, / \, (1- c\iDB),
\end{align*}
since $\zeta(\cD)\vee \nu_j A_{t\iDB}(\cD) \le 1$,\; 
$c_{j} \le c\iDB$ and 
$:c \rightarrow 1 - (1 - \nu_j A_{t\iDB}(\cD))\, /\, (1- c)$  is decreasing.
In particular, $\nu_{j+1}\in \M_1\left(\cX\right)$ holds true 
and thanks again to  \Req{U:Mrn}
we prove finally:
\begin{align*}
\nu_{j+1}(\cD\iRN) \ge \frac{1}{1-c\iDB} \, \left( \nu_j A_{t\iDB} (\cD\iRN) - c\iDB\right)
\ge \xi\iRN.
\end{align*}
Therefore, $(I_{j+1})$ holds.

By induction, we get $(I_{J(t_h)})$ thus $r_{J(t_h)} >0$ and  $\nu_{J(t_h)} \in \M\iRN \subset \M_1\left(\cX\right)$. 
By \Req{U:nuJ}, i.e. the definition of $\nu_{J(t_h)}$, and since $(A_\tp)$ is a semigroup:
\begin{align*}
		&\mu A_{t_h } 
		= \dfrac{\PR_\mu (J(t_h)\, t\iDB < \ext)}{\PR_\mu (t_h < \ext)} \;
		\mu A_{ J(t_h)\, t\iDB} P_{t_h -  J(t_h)\, t\iDB}
		\\&\hspace{1cm}
		\ge \dfrac{\PR_\mu (J(t_h)\, t\iDB < \ext)}{\PR_\mu (t_h < \ext)} \;
		\left[ \Tsum{k \le J(t_h)} \amu(k, J(t_h)\, t\iDB)\, 
		\zeta A_{[J(t_h) - k]\, t\iDB} \;
		P_{t_h -  J(t_h)\, t\iDB} \right]
		\\&\hspace{1cm}
		\ge  \sum_{k \le J(t_h)} \dfrac{\PR_\mu (J(t_h)\, t\iDB < \ext)}{\PR_\mu (t_h < \ext)} \;
		\amu(k, J(t_h)\, t\iDB)\, 
		\dfrac{\PR_{\zeta} (t_h -k\, t\iDB< \ext)}
		{\PR_{\zeta} ([J(t_h) -k]\, t\iDB< \ext)}
		\zeta A_{t_h - k\, t\iDB}.
\end{align*}
Finally, thanks to \Req{U:AMU}\, and \Req{U:AlcO}, 
we conclude:
$\qquad \mu A_{t_h} \ge \zeta[t_h].$

\subparagraph{Step 2.2:} $\mu \in \M_1\left(\cX\right)$.
For general initial conditions, recall that in Proposition \ref{U:c1}, we constructed $\M\iRN$ such that \mbox{$\M\iXT \subset \M\iRN$}. 
Thus, \Req{U:Mxt} holds with $\M\iRN$ instead of $\M\iXT$.
Since $t_{h, 2} \ge t_{h, 1}$, we obtain
$\, \mu A_{t\iXT + t_{h, 2}-t_{h, 1}} \in \M\iRN.$
As $\mu A_{t_{h, 2}} = \mu A_{t\iXT + t_{h, 2}-t_{h, 1}} A_{t_{h, 1} -t\iXT}$, we 
finally deduce from \Req{mulJ}: 
$$
\mu A_{t_{h, 2}} \ge \zeta[t_{h, 1} - t\iXT].$$
\hfill $\square$

\paragraph{Step 3: conclusion of the proof of Theorem \ref{U:AllPho}} 
\label{U:sec:prTh}
\textcolor{white}{:}

Thanks to Proposition \ref{U:LowerBound}:\newcounter{ZETA}
\begin{align*}
	&
	\frl{\ell \in \N}
	\frl{\xi \in (0, 1)}
	\frl{(\mu_1,\, \mu_2) \in (\M_{\ell,\, \xi}  )^2}
	\frl{t_{h, 2}\ge t_{h, 1} \ge  t\iXT} \\
	&\hspace{1cm}
	\mu_1 A_{t_{h, 1}} 
	\ge \zeta[t_{h, 1} - t\iXT] \quad
	\text{ and }\quad
	\mu_2 A_{t_{h, 2}}
	\ge \zeta[t_{h, 1} - t\iXT]\\
	& \text{ thus } \|\mu_2 A_{t_{h, 2}} - \mu_1 A_{t_{h, 1}}\|_{TV} 
	\le \|\mu_2 A_{t_{h, 2}} - \zeta[t_{h, 1} - t\iXT]\|_{TV}  + 
	\|\mu_1 A_{t_{h, 1}} - \zeta[t_{h, 1} - t\iXT]\|_{TV} \\
	&\hspace{4.7cm}
	\le 2 \times [ 1- \zeta[t_{h, 1} - t\iXT](\cX)].
	\EQn{muA21}{}
	\end{align*}
\begin{align*}
	&1- \zeta[t_h](\cX) 
	= 1- \Tsum{j \le J(t_h)} (\text{\large{$\,^{c\iDB}\!/_{c\iPs}$} })\, \left[ 1-\text{\large{$\,^{c\iDB}\!/_{c\iPs}$} } \right]^{j-1}
	\hcm{1} \text{ by } \Req{U:AlcO},
	\\&\hspace{2cm}
	=  \left[ 1-\text{\large{$\,^{c\iDB}\!/_{c\iPs}$} } \right]^{J(t_h)}
	\le   \exp[- \gamma \, (t_h - t\iPs - t\iDB)]
	\hcm{1} \text{ by } \text{\Req{U:Jtobs},}
	\\&\hspace{3cm}
	\EQn{U:ZETA}{\gamma}
	\quad \text{ with }\; \gamma:= -\frac{1}{t\iDB} \, \ln \left[ 1-\dfrac{c\iDB}{c\iPs} \right]
	\end{align*}
Finally, with \Req{muA21}, \Req{U:ZETA}\, and 
$\quad C = C(\ell,\, \xi):= 2\,  \exp[ \gamma \, (t\iPs+t\iDB + t\iXT(\ell, \xi))]$:
\newcounter{CNX}
\begin{align*}
	&\|\mu_2 A_{t_{h, 2}} - \mu_1 A_{t_{h, 1}}\|_{TV} \le 
	C\; e^{-\gamma\, t_{h, 1}}.
	\EQn{U:CNX}{C(\ell,\, \xi)}
	\end{align*}
This states that for any $\mu \in \M_{\ell,\, \xi}  $, 
$(\mu A_{t_h})_{\left \lbrace t_h\ge 0\right \rbrace}$ 
is a Cauchy-sequence 
for the total variation distance.
Thus, 
it converges for this distance to some distribution $\alpha^{\ell, \xi}$. 
Since for any $\ell \le \ell'$ and $\xi \ge \xi'> 0$,
it is clear by definition 
that $\M_{\ell,\, \xi}   \subset  \M_{\ell',\, \xi'},$
we deduce $\alpha^{\ell, \xi} 
= \alpha^{\ell \vee \ell',\, \xi \wedge \xi'}$
$=\alpha^{\ell',\, \xi'} $.
This  means 
(since \mbox{$\M_1\left(\cX\right) = \medcup_{(\ell, \xi)} \M_{\ell,\, \xi}  $})
that a unique distribution $\alpha$
is the attractor.
In particular, 
there cannot be a QSD 
different from $\alpha$. 

For any initial condition $\mu$:
$\hcm{0.5} \limInf{\tp} \PR_\mu(X_\tp \in dy\, \big| \, \tp < \ext ) = \alpha(dy),\quad$
where the convergence holds in the weak topology 
(ie $\alpha$ is a \textit{quasi-limiting distribution}).
One can then easily adapt the proof of Lemma 7.2 in \cite{Cat} to deduce that $\alpha$ is effectively a QSD and $\frl{\tp \ge 0}$ 
$\PR_\alpha (\tp < \ext) = e^{-\lambda\, \tp}$.
By letting $t_{h, 2}\rightarrow \infty$ in \Req{U:CNX}, with $\mu_2 = \mu_1 = \mu \in \M_{\ell,\, \xi}  $:
\begin{align*}
	\|\, \PR_\mu \left[\, X_{t} \in dx \; 
	| \; \tp < \ext \right]  - \alpha(dx) \, \|_{TV}
	\le C(\ell, \xi) \; e^{-\gamma \; \tp} 
	\end{align*}
This ends the proof 
of  Theorem \ref{U:AllPho} (up to Appendix E).
\hfill $\square$

\subsection{Proof of Theorem~\ref{U:EtaECV}:}
\label{U:sec:eta}

\subsubsection*{Step 1: proof of the uniform convergence to $\heig$}
\label{U:sec:Eta}

Considering the arguments in the proof of Theorem~\ref{U:PersW},
it is easily seen that 
for any probability measure $\mu$,
there exists $c'\iPs$, $t'\iPs$
such that:
	\begin{align*}
	&\frl{x\in \cX}
	\frl{\tp \ge t'\iPs}
	\quad
	\PR_{x} \left( \tp < \ext\right)  
	\le c'\iPs\; \PR_{\mu} \left( \tp < \ext\right).
\end{align*}
Here, we need this estimate for $\mu:= \alpha$.
%
To achieve this, we only need to apply \textup{$(A1)$}\  and adjust the value for $c\iPs$:
$\quad 
c'\iPs:= 
c\iPs\, e^{-\lambda\,t\iMix} \, /\, (\alpha(\cD_{\ell\iMix})\, c\iMix),
\quad$
where $t\iMix, c\iMix$ are given by \textup{$(A1)$}\, for initial condition in $\cD_{\ell\iMix}$.
This can be translated in term of a uniform bound on $\heig$ by: 
\begin{align*}
	&	\| \heig_\bullet \|_\infty:= \Tsup{\tp \ge 0}
\| \heig_\tp \|_\infty \le c'\iPs \vee e^{\lambda\, t\iPs}  < \infty.
	\EQn{U:etaInf}{\| \heig_\bullet \|_\infty}
\end{align*}
Like in the proof of Proposition~2.3 in \cite{ChQSD}, we deduce that, 
for any $s, t >0$, $\mu \in \M_{\ell,\, \xi}  $:
\begin{align*}
	|\langle\mu|\heig_{t}\rangle - \langle\mu|\heig_{t+s}\rangle| 
	= \langle\mu|\heig_{t}\rangle \,
	|\langle\alpha - \mu|\heig_{s}\rangle| 
	\le \| \heig_\bullet \|_\infty^2\; C(\ell, \xi) \, e^{-\gamma\, \tp}.
	\EQn{U:Hcauchy}{}
\end{align*}
The constant $C$ can actually 
be taken independently 
of $\ell, \xi$.
Indeed, because the previous expression 
is linear in $\mu$ 
and $\langle \alpha \, \big| \, \heig_t\rangle \equiv 1$:
$$|\langle\mu\, \big| \, \heig_{t} - \heig_{t+s}\rangle|
= 2|\langle\bar{\mu}\, \big| \, \heig_{t} - \heig_{t+s}\rangle|,
\quad \text{ where } \bar{\mu}:= (\mu + \alpha)/2.
$$
By choosing $\ell$ sufficiently large 
to ensure $\xi:= \alpha(\cD_{\ell\iMix})/2 >0$,
we deduce that for any $\mu \in \M_1\left(\cX\right)$, 
$\bar{\mu} \in \M_{\ell,\, \xi}  $.
The inequality \Req{U:Hcauchy} 
is thus uniform in $\mu \in \M_1\left(\cX\right)$,
so that $(\heig_\tp)$ defines 
a Cauchy sequence for the uniform norm. 
We deduce that $\heig_\tp$ converges 
to some unique function $\heig$,
whose norm is also bounded by $\| \heig_\bullet \|_\infty$.
\hfill $\square$

\subsubsection*{Step 2: Characterization of the survival capacity $\heig$}
\indent
The rest of the proof is directly taken from \cite{ChQSD}.
As the punctual limit of $(\heig_\tp)$, and since for any $\tp \ge 0$, $\heig_\tp$ vanishes on $\partial$, this also hold for $\heig$. With the uniform bound \Req{U:etaInf}, we deduce that $\heig$ is also bounded.
As stated in the beginning 
of this Subsection~\ref{U:sec:eta}, 
we can replace $\zeta$ 
by any probability measure $\mu$ in (\ref{U:eqPS})\, 
(with specific values for 
$ c\iPs(\mu) = c\iPs(\ell)/\xi, 
t\iPs(\mu) = t\iPs(\ell) >0$). 
In particular, 
for $\mu = \delta_x$, 
with $x\in \cD_\ell$:
\begin{align*}
	&\frl{\tp \ge t\iPs(\ell) } \quad
	\PR_\alpha (\tp <\ext) 
	\le c\iPs(\ell) \, \PR_x (\tp <\ext)
	\quad\text{thus } 
	\frl{\tp \ge t\iPs(\ell) }\quad  
	\heig_\tp(x)\ge c\iPs(\ell)>0.
\end{align*}
This proves that $\heig$ has a positive lower-bound on any $\cD_\ell$.
By the Markov property and \Req{U:LBz}:
\begin{align*}
	\frl{u>0}\quad 
	&P_u \heig(x) 
	= \limInf{t} \dfrac{\E_x\left[\PR_{X_u}( \tp <\ext)\right]}{
		\PR_\alpha(\tp <\ext)}
	= e^{- \lambda\, u} \limInf{t} 
	\dfrac{\PR_x(t+u<\ext)}{\PR_\alpha(t+u<\ext)}
	= e^{- \lambda\, u} \, \heig(x).
\end{align*}
From this and \Req{U:etaInf}, we immediately deduce that $\heig$ is in the domain of $\mathcal{L}$ and $\mathcal{L} \, \heig = -\lambda\, \heig$.

\hfill $\square$

\subsection{Proof of Theorem~\ref{U:QECV}:}
\label{U:sec:QECv}
Except for $(iii)$, for which we will prove  \Req{U:ECvBeta}, 
and for the uniqueness of the stationary distribution,
the proof is almost the same as in \cite{ChQSD}.

\subsubsection*{Step 1: Proof that the $Q$-process is well-defined and characterization}

Let $\Lambda_s$ be a $\mathcal{F}_s$-measurable set and $\mu\in \M_1\left(\cX\right)$. By the Markov property:
\begin{align*}
	&\PR_\mu(\Lambda_s \, \big| \, \tp < \ext) 
	= \E_\mu\left[ e^{\lambda s} \,\heig_{t-s}(X_s) / \langle\mu\, \big| \, \heig_\tp \rangle;\; s<\ext,\; \Lambda_s\right].
\end{align*}
By Theorem~\ref{U:EtaECV}, the random variable
$
M_s^t:= \idc{s<\ext} \,e^{\lambda s}\, 
\heig_{t-s}(X_s)\ /\ \langle\mu\, \big| \, \heig_\tp \rangle,$
(where $\tp \ge s$)
converges a.s. to:
\begin{align*}
	M_s:=  \idc{s<\ext} \,e^{\lambda s}\, 
	\heig(X_s)\ /\ \langle\mu\, \big| \, \heig \rangle,
\end{align*}
where $\langle\mu\, \big| \, \heig \rangle >0$ because $\heig$ is positive on $\cX$. 
For $t$ sufficiently large (a priori depending on $\mu$), 
we deduce from \Req{U:etaInf} 
and  the convergence of $\langle \mu\, \big| \, \eta_t\rangle$ to $\langle \mu\, \big| \, \eta\rangle$:
%
\begin{equation}\EQn{U:Mst}{}
	0 \le M_s^t 
	\le 2\, e^{\lambda s}\,\| \heig_\bullet \|_\infty\ /\ \langle\mu\, \big| \, \heig \rangle.
\end{equation}
Thus, by the dominated convergence Theorem, we obtain that
$\qquad E_\mu (M_s) = 1.$

By the penalisation’s theorem of Roynette, Vallois and Yor 
(cf Theorem~2.1 in \cite{Roy06}) these two conditions imply that $M$ is a martingale under $P_\mu$ and
that \mbox{$\PR_\mu(\Lambda_s \, \big| \, \tp < \ext) $}
converges to $E_\mu (M_s;\; \Lambda_s)$ for all $\Lambda_s \in \mathcal{F}_s$ when $t \rightarrow \infty$. 
In particular for $\mu = \delta_x$,
this means that $\mathbb{Q}_x$ is well defined and:
\begin{align*}
	\dfrac{d\mathbb{Q}_x}{d\PR_x} _{\, \big| \,\mathcal{F}_s}  
	=  \idc{s<\ext} \,e^{\lambda s}\, \dfrac{\heig(X_s)}{\heig(x)}.
	\EQn{dQdP}{d\PR/d\mathbb{Q}}
\end{align*}
\Req{dQdP}\, implies directly \Req{U:qt}. Concerning \Req{U:etaQ}:
\begin{align*}
	&\mu  B[\heig] Q_\tp (dy)
	= \int \dfrac{\heig(x)} {\langle \mu\, \big| \, \heig \rangle} \mu(dx)\, \dfrac{\heig(y)}{\heig(x)}\,
	e^{\lambda\, \tp}\, p(x; t; dy)
	\\&\hcm{1}
	= \dfrac{\heig(y)} {\langle \mu\, \big| \, P_\tp\, \heig \rangle}\, \mu P_\tp(dy)
	= \mu P_\tp B[\heig]
	\hcm{1} \text{ by } \Req{U:EFeta},
	\\&\hcm{1}
	= \dfrac{\heig(y)\, \PR_\mu(\tp <\ext)} {\langle \mu\, P_\tp\, \big| \, \heig \rangle}
	\times \dfrac{\mu P_\tp(dy)}  {\PR_\mu(\tp <\ext)}
	= \dfrac{\heig(y)} {\langle \mu\, A_\tp\, \big| \, \heig \rangle}\, \mu A_\tp(dy)
	= \mu   A_\tp B[\heig].
\end{align*}

For a more general convergence, with $\mu$ as initial condition and $\Lambda_s \in \mathcal{F}_s$, we deduce:
\begin{align*}
	&\limInf{t}
	\PR_\mu(\Lambda_s \, \big| \, \tp < \ext) 
	= E_\mu \left( e^{\lambda s}\, \heig(X_s)\, /\, \langle\mu\, \big| \, \heig \rangle 
	;\; s<\ext, \Lambda_s\right)
	\\&\hspace{1cm}
	=  \int_\cX \mu(dx) \dfrac{\heig(x)}{\langle\mu\, \big| \, \heig \rangle}
	E_x \left( e^{\lambda s}\, \dfrac{\heig(X_s)}{\heig(x)}
	;\; s<\ext, \Lambda_s\right)
	\\&\hspace{1cm}
	= \int_\cX \mu(dx) \dfrac{\heig(x)}{\langle\mu\, \big| \, \heig \rangle} \mathbb{Q}_x(\Lambda_s)
	= \mathbb{Q}_{\mu B[\heig]}(\Lambda_s).
\end{align*}
by $\Req{dQdP}\,$ and the definition of $\mu  B[\heig]$ in 
\Req{U:etaEt}.

Moreover, the convergence holds in fact in total variation over $\mathcal{F}_\spr$, as we prove it now. By the previous calculations, \Req{U:Mst}\, and \Req{U:etaInf}, for any $\epsilon>0$:
\begin{align*}
	&
	\left\|    \PR_\mu(dw \, \big| \, \tp < \ext)  
	- \mathbb{Q}_{\mu  B[\heig]}(dw)
	\right\|_{TV, \mathcal{F}_\spr}
	\le \E_\mu |M_\spr^\tp - M_\spr|
	\\&\hspace{1cm} 
	\le 4\, e^{\lambda \spr}\, \dfrac{\| \heig_\bullet \|_\infty}{\langle\mu\, \big| \, \heig \rangle}\,
	\PR_\mu(|M_\spr^\tp - M_\spr| \ge \epsilon) 
	+ \epsilon,
	\\&
	\text{so }\quad
\underset{\tp\rightarrow \infty}{\limsup}\,
	\left\|    \PR_\mu(dw \, \big| \, \tp < \ext)  
	- \mathbb{Q}_{\mu B[\heig]}(dw)
	\right\|_{TV, \mathcal{F}_\spr} 
	\le \epsilon.
\end{align*}
By letting $\epsilon\rightarrow 0$, we conclude:
\begin{align*}
	\frl{\spr\in \R_+}\quad
	\left\|    \PR_\mu(dw \, \big| \, \tp < \ext)  
	- \mathbb{Q}_{\mu B[\heig]}(dw)
	\right\|_{TV, \mathcal{F}_\spr}
\underset{\tp\rightarrow \infty} {\longrightarrow}0.
\end{align*}

For the proof that $X$ defines a strong Markov process under $(\mathbb{Q}_x)_{x\in \cX}$, we refer again to the proof in \cite{ChQSD}.

\subsubsection*{Step 2: The invariant distribution for $X$ under $\mathbb{Q}$}
For all $\tp \ge 0$ and $f\in \B_b(\cX)$, with \Req{dQdP}:
\begin{align*}
	\langle\beta \, \big| \, Q_\tp\, f \rangle
	&= \langle\alpha \, \big| \, \heig\!\times\! Q_\tp\, f \rangle
	= e^{\lambda\, \tp}\,\langle\alpha \, \big| \, P_\tp\, (\heig\!\times\! f) \rangle
	\\&
	= \langle\alpha \, \big| \, \heig\!\times\! f \rangle 
	= \langle\beta \, \big| \,  f \rangle,
	\EQn{U:BQ}{\beta Q_\tp}
\end{align*}
where we used $\Req{U:LBz}$. We prove the uniqueness with the next subsection.

\subsubsection*{Step 3: Proof of  \Req{U:ECvBeta} }

Exploiting \Req{U:etaQ}, we deduce from our definitions:
\begin{align*}
	\| (\mu B[\heig]) Q_t
	- \beta \|_{\frac{1}{\heig}}
	&=  \left \| \dfrac{\mu A_t}{\langle \mu A_t\, \big| \, \heig\rangle}
	- \alpha \right\|_{TV}
	\\&
	\le  \langle \mu A_t\, \big| \, \heig\rangle^{-1}\ltm
	[\left\|  \mu A_t- \alpha \right\|_{TV}
	+ |\langle \mu A_t\, \big| \, \heig\rangle - 1|].
	\EQn{U:EBeta}{}
\end{align*}
To ensure a lower-bound on $\langle \mu A_t\, \big| \, \heig\rangle$,
we exploit \Req{U:EFeta} and write:
\begin{equation*} 
	\langle \mu A_t\, \big| \, \heig\rangle
	= \dfrac{e^{\lambda t} \langle \mu P_t\, \big| \, \heig\rangle}
	{\langle \mu\, \big| \, \heig_t\rangle}
	= \dfrac{\langle \mu\, \big| \, \heig\rangle}
	{\langle \mu\, \big| \, \heig_t\rangle}.
\end{equation*}
We already know that $\heig_t$ is uniformly upper-bounded 
and $h$ has a lower-bound on any $\cD_\ell$.
Since $|\langle \mu A_t\, \big| \, \heig\rangle - 1|
= |\langle \mu A_t - \alpha\, \big| \, \heig\rangle|
\le \left\|  \mu A_t- \alpha \right\|_{TV}\, \| \heig_\bullet \|_\infty,$
and exploiting \Req{U:EBeta} and \Req{U:ECvAl}, we conclude that there exists $C' = C'(\ell, \xi)>0$ such that:
\begin{equation*} 
	\quad\frl{\tp  > 0}
	\frlq{\mu \in \M_{\ell,\, \xi}  }
	\| \mathbb{Q}_{\mu B[\heig]} ( X_{\tp } \in dx)
	- \beta(dx) \|_{\frac{1}{\heig}}
	\le C'\; e^{-\gamma \; \tp }.
	\tag*{$\square$}
\end{equation*}

\subsubsection*{Step 4: Convergence with initial condition for the Q-process}

When $\mu_Q$ is the initial condition of the Q-process, it is in general not possible to interpret it as $\mu  B[\heig]$. Indeed, we should expect in this case $\mu(dx)$ to be proportional to $\heig(x)^{-1}\, \mu_Q(dx)$, which may not be integrable. Thus, the convergence to $\beta$ might in general not be exponential. 

However, it is exponential for measures 
with support in any of the $\cD_\ell$, 
in particular Dirac masses.
Indeed, we have a lower-bound of $\heig$:
$\quad \heig^{(\ell)}:= \inf\left \lbrace \heig_x;\; x\in \cD_\ell \right \rbrace,\quad $
which is positive because of \textup{$(A1)$}\, and \Req{U:eqPS}.
Thus, if $\mu_Q\in \M_1\left(\cX\right)$ has support on $\cD_\ell$, 
$\langle \mu_Q\, \big| \, 1/\heig \rangle \le 1/\heig^{(\ell)}< \infty,$ so:
\begin{align*}
	\mu_Q = \mu  B[\heig],
	\quad \text{ with } \mu(dx):= \mu_Q  B(1/\heig)
	:= \mu_Q(dx) \, /\, (\heig(x)\ltm\langle \mu_Q \, \big| \, 1/\heig\rangle ).
\end{align*}
Now, $\mu$ has the same support as $\mu_Q$, thus $\mu(\cD_\ell) = 1$, i.e. $\mu \in \M_{\ell, 1}$. By \Req{U:ECvBeta}:
\begin{align*}
\left\|  \mu_Q Q_\tp - \beta \right\|_{TV}
	= \left\|  \mu  B[\heig] Q_\tp -  \beta \right\|_{TV}
	\le C(\ell, 1)\, e^{-\gamma\, \tp}.
\end{align*}
More generally, 
since the Q-process is linear 
with its initial condition,
and by \textup{$(A0)$},
the property of uniqueness of the stationary distribution $\beta$ holds. \\

Besides, to have exponential convergence, it suffices that: 
$\langle \mu_Q\, \big| \, 1/\heig \rangle < \infty$. It can be deduced from 
$\sum_{\ell \ge 1} \mu_Q(\cD_\ell \setminus \cD_{\ell-1})\, /\heig^{(\ell)} < \infty$ 
(note that one has lower-bounds of $\heig^{(\ell)}$).
In any case, the convergence still holds in total variation.
\hfill $\square$

\section*{Appendices: }
\setcounter{eq}{0}
\stepcounter{section}
\subsubsection*{\textbf{Appendix A: }
	Combine all the inequalities to prove Proposition \ref{A2XN}}

We shall first prove that 
an upper-bound of the global supremum
can be deduced from the upper-bounds in
Propositions  \ref{U:cEinf}-5.
So we start by assuming that the inequalities
 derived in these propositions 
hold for some parameters $\epsilon^X$, $\epsilon^0$, $C^N_\infty$, $C^X_\infty$  and $C_0$
($C^N_\infty$ coming from Proposition \ref{U:cEinf};
$\epsilon^X$ and $C^X_\infty$ from  Proposition \ref{U:cEXinf};
  $\epsilon^0$ and $C_0$ from  Proposition \ref{U:cE0})
and explain how these inequalities 
can imply the global supremum in \textup{$(A2)$}.
This implication shall hold
at least for $\epsilon^X$ and $\epsilon^0$ sufficiently small, 
which is obtained with $n\iET$ sufficiently large. 
The constraints on $\epsilon^X$ and $\epsilon^0$ are mentioned 
while we handle the inequalities.
We prove next that we can indeed find suitable choices
of $\epsilon^X$, $\epsilon^0$, $C^X_\infty$, $C^N_\infty$ and $C_0$
for the upper-bounds in Propositions  \ref{U:cEinf}-5
to hold with these constraints.

$t_{h}$ is introduced to make sure 
that $\cE^X_\infty\vee \cE^N_\infty\vee \cE_0 
< \infty$ ($\le \exp[\rho\, t_{h}]$).
It is needed to justify 
the following inequalities,
but this specific upper-bound plays no role.
By the upper-bounds in Propositions \ref{U:cEXinf} and \ref{U:cEinf}:
\begin{align*}
	&\cE^X_\infty
	\le C^X_\infty \, \left( 1+ \cE_0 \right) 
	+ \epsilon^X\; C_{\infty}^N  (1 +\cE^X_\infty)\\
	&(1- \epsilon^X\; C_{\infty}^N ) \, \cE^X_\infty
	\le C^X_\infty  + \epsilon^X\; C_{\infty}^N   +C^X_\infty  \,\cE_0.
\end{align*}
Provided that:
$\epsilon^X\le  (2\,C_{\infty}^N )^{-1}$, 
recalling that $C_{\infty}^N \wedge C^X_\infty \ge 1$, 
and combining it with the upper-bounds of Proposition \ref{U:cE0}, it yields:
\begin{align*}
	&
	\cE^X_\infty
	\le 3\,C^X_\infty + 2\,C^X_\infty\, \cE_0,
	\hcm{1}
	\cE^N_\infty
	\le 4\, C_{\infty}^N \, C^X_\infty 
	+ 2\, C_{\infty}^N \, C^X_\infty \, \cE_0
	\\ &\hspace{1cm}
	\cE_0 \le  C_0 
	+ 7\, \epsilon^0\, C_{\infty}^N \, C^X_\infty 
	+ 4\, \epsilon^0\, C_{\infty}^N \, C^X_\infty  \cE_0
	\\ &\hspace{.5cm}
	\text{thus }
	\left( 1 - 4\,\epsilon^0C_{\infty}^N  \,  C^X_\infty\right)  \; \cE_0
	\le C_0 + 7\,\epsilon^0C_{\infty}^N  \,  C^X_\infty.
\end{align*}
Provided:
$\,\epsilon^0 \le  \left( 8\,C_{\infty}^N \, C^X_\infty\right)^{-1},$ 
and recalling that $C_{\infty}^N \wedge C_0 \ge 1$, we deduce:
\begin{align*}
	&
	\cE_0 
	\le 4\,C_0,\;\qquad
	\cE^X_\infty
	\le 11\, C^X_\infty \, C_0,\;
	\qquad
	\cE^N_\infty
	\le 12\, C_{\infty}^N \,C^X_\infty \, C_0.
\end{align*}
Finally, provided:
$\,\epsilon^X\le  (2\,C_{\infty}^N )^{-1},\,
\epsilon^0 
\le  \left( 8\,C_{\infty}^N \, C^X_\infty\right)^{-1},$
conditions which we can satisfy 
and restrict the choices of $n_{\infty}$ and $n\iET> n_{\infty}$,
we deduce:
\begin{equation*} 
	sup_{(x, n)}\left \lbrace \E_{(x, n)}[\exp(\rho \; \widehat{\tau}_{E})]\right \rbrace 
	\le 12\, C_{\infty}^N \,C^X_\infty \, C_0
	< \infty.
	\EQn{U:eT}{}
\end{equation*}

More precisely, for any $\rho$, we obtain from Proposition \ref{U:cEinf}
the constants $n_\infty$ and $C_{\infty}^N $, so that we can set $\epsilon^X:=(2\,C_{\infty}^N )^{-1}$.
We then deduce, thanks to Proposition \ref{U:cEXinf}, some value for $n\iET$ and $C^X_\infty$.
Setting $\epsilon_0=(8 C^N_\infty C^X_\infty)^{-1}$,
we can choose, according to Proposition \ref{U:cE0}, 
some value $n_0> 0$ and $C_0$. 
Taking the limit in \Req{U:eT} as $t_{h} \rightarrow \infty$
(recall that $\widehat{\tau}_{E}:= \tau_{E} \wedge \ext\wedge t_{h}$)
and choosing $n:=n_\infty\vee n_E\vee  n_0$
conclude the proof of Proposition \ref{A2XN}.

\hfill $\square$

\subsubsection*{\textbf{Appendix B:} 
	Descent from infinity, proof of Proposition \ref{U:cEinf}} 
\label{U:sec:cEinf}

\paragraph{Lemma \ref{U:ydP} implies Proposition \ref{U:cEinf}}
\textcolor{white}{:}

We obtain by induction and the Markov property:
$\quad \frl{ n > 0} \quad
\PR_{n} (k\, t < \tau^D_{\downarrow} ) 
\le \epsilon^k. $\\
Thus, by choosing $\epsilon$ sufficiently small (for any given value of $t>0$), we ensure:
\begin{align*}
	&\hspace{1cm} 
	C_{\infty}^N 
	:=\Tsup{n > 0}\left \lbrace 
	\E_{n}[\exp(\rho \tau^D_{\downarrow} )]\right \rbrace 
	< +\infty.
\end{align*}
\noindent
A fortiori with $
T_{\downarrow} 
:= \inf\left \lbrace t, N_t \le  n_{\infty} \right \rbrace \wedge \tau_{E}
\le  \tau^D_{\downarrow} ,$
$$
\sup_{(x, n)}\left \lbrace \E_{(x, n)}[\exp(\rho \, T_{\downarrow} )]\right \rbrace 
\le  C_{\infty}^N  < \infty.$$
At time  $T_{\downarrow} $, the process is either in $E$ or in $\cT^X_\infty$. Thus:
\begin{align*}
	&\E_{(x, n)}[\exp(\rho \, \widehat{\tau}_{E} )]
	\le \E_{(x, n)}[\exp(\rho \, T_{\downarrow} ) \,
	; \; (x, n)_{T_{\downarrow}} \in E]\\
	&\hspace{2cm}
	+\E_{(x, n)}\Big[\exp(\rho \, T_{\downarrow} )
	\E_{(X, N)_{T_{\downarrow}} }
	[\exp(\rho \, \widehat{\tau}_{E} )] \,
	; \;(X, N)_{T_{\downarrow}} \in \cT^X_\infty \Big],
\end{align*}
with the Markov property
and the fact that $(\widehat{\tau}_{E} - \cT^X_\infty)_+ \le t_{h}$
on the event\\
$\{(X, N)_{T_{\downarrow}} \in \cT^X_\infty\}$.
Therefore:
\begin{align*}
	\cE^N_\infty \le   C_{\infty}^N\, (1+ \cE^X_\infty).
	\hcm{2}
	\tag*{$\square$}
\end{align*}

\paragraph{Proof of Lemma \ref{U:ydP}}
\textcolor{white}{:}

The proof of this Lemma 
relies mainly on the same arguments as in \cite{BP12}, part 6,
related to the descent from infinity.
Let $Z_t:= \sigma/2 \times \sqrt{N^D_t}$.
It is solution to the following EDS:
\begin{align*}  
	& Z_t
	:= z +  \int_{0}^{t} \psi \left( Z_s \right) \;ds  
	+ B_t, \text{ where }
	\psi(z):= - \frac{1}{2\, z} +\frac{r\, z}{2} 
	- c\, z^3.
	\EQn{U:Z2}{}
\end{align*}
%
As long as $Z$ is very large 
and $|B|$ not exceptionally large, 
the leading term $- c\, Z_t^3$
indeed makes the process come down in finite time.
Let $V:= Z - B$. It is the solution of the ODE:
\begin{align*}
	&\dfrac{dV_t}{dt} 
	= - \dfrac{1}{2 (V_t + B_t) } + \dfrac{r\; (V_t + B_t)}{2} 
	- \cY\; (V_t + B_t)^3,
	\EQn{U:Vt}{V_t}
	\\&
	\text{Let }
	z_2 \ge z_1:= \sup\left \lbrace z >0, \;
	\left|- \frac{1}{2 z} + \frac{r z}{2} \; \right| 
	\ge \frac{\cY  z^3}{2}\right \rbrace,
	\EQn{U:wyinf}{z_1}
	\\&
	T_B:= \inf\left \lbrace t >0, B_t \notin 
	[- z_2,\;  2\, z_2] \right \rbrace,
	\hspace{1.2cm}
	T_V:= \inf\left \lbrace t >0, V_t < 2 z_2 \right \rbrace,
	\EQn{U:TV}{T_V}
\end{align*}
where we consider w.l.o.g. 
an initial condition $z$ strictly bigger than $2\, z_2$, so that $T_V$ is positive a.s.
Then, as in \cite{BP12}, 
we get on the time interval $[0, T_B \wedge T_V]$: 
\begin{align*}
	&\hcm{0.5}
	B_t \ge - z_2 \ge - V_t/2, \text{ implying } \;
	V_t + B_t \ge  V_t/2 \text{ and } V_t + B_t  \ge  z_2,     
	\\& \hcm{1.5}
	\left| - \dfrac{1}{2 (V_t + B_t) } + \dfrac{r\; (V_t + B_t)}{2} \right|
	\le \dfrac{\cY}{2} (V_t + B_t)^3,
	\\&
	\dfrac{d}{dt} \left[ (V_t)^{-2} \right] 
	= \frac{-2}{(V_t)^3}\dfrac{dV}{dt}
	\ge 2 \times \left( \cY 
	- \dfrac{\cY}{2} \right)
	\times 	\left( \frac{V_t + B_t}{V_t} \right) ^3
	\ge \dfrac{\cY}{8} , 
	\\&
	\text{ thus }
	V_t^{(-2)} - z^{(-2)}
	\ge \cY\, t / 8  \quad
	\text{ and in particular } \qquad
	V_t \le \sqrt{8  / (\cY \ltm t)}.
\end{align*}
Thus, $\left \lbrace t \le T_B \right \rbrace
\subset \left \lbrace T_V \le t\right \rbrace
\cup \left \lbrace V_{t} \le \sqrt{8  / (\cY \ltm t)}\right \rbrace$.

By \Req{U:TV}, let $z_2$ be sufficiently big to ensure:
$\, \PR(T_B < t) \le \epsilon.$
Then, denote: \\$z_\infty:= 
\left( \sqrt{8  / (\cY \ltm t)} + 2\, z_2 \right) \vee (4\, z_2).$
We deduce that, on the event  $\left \lbrace t \le T_B \right \rbrace$,
either $Z_t \le z_\infty$ 
or $T_V \le t$ while
	$Z_{T_V} \le 4 z_2 \le z_\infty$.
In any case, $\tau^D_{\downarrow}  \le t$. Hence:
$\, \frl{z > 0} \PR_{z} (t < \tau^D_{\downarrow} )  
\le \epsilon.$
\hfill $\square$

\subsubsection*{\textbf{Appendix C:}
	Mal-adaptation too large, 
	proof of Proposition \ref{U:cEXinf}:}
\label{U:sec:cEXinf}

\paragraph*{Lemma \ref{U:N^DP} implies Proposition \ref{U:cEXinf}}
\textcolor{white}{:}

Let $\rho,\, \epsilon,\, n_{\infty} > 0$ ($\cY>0$ is the same as for the definition of $Z$). 
For simplicity, 
we choose $\tDo:= \log(2)/\rho >0$ (i.e. $\exp\left[ \rho\, \tDo \right] = 2$),
and assume w.l.o.g. $\tDo < t_h$.
We choose $r_\vee \in \R$ according to Lemma \ref{U:N^DP} 
such that:
\begin{align*}
	&
	\frl{n > 0}
	\frl{r \le r_\vee}\quad
	\PR_{n} \left(  \tDo < \ext^D \right) 
	\le e^{-\rho\, \tDo}/2 = 1/4,
	\\&\frlq{r \le r_\vee}
	\PR_{n_{\infty}} \left( T^D_{\infty} \le \tDo \right)  
	+ \PR_{n_{\infty}} \left( N^D_{\tDo} \ge n_{\infty} \right) 
	\le  \epsilon /4.
\end{align*}
Since $\limsup_{\|x\| \rightarrow \infty} r(x) = -\infty$,
with $n\iET$ chosen sufficiently large:
$$\frl{x\notin B(0, n\iET)} 
r(x) \le~r_\vee.$$

Let $(X, N)$ with initial condition $(x, n)\in \cT^X_\infty$.
In the following, we denote:
\begin{align*}
	&T^N_\infty: = \inf \left \lbrace t\ge 0,\, N_t \ge n_c \right \rbrace,
	\quad
	\tau_0:= \inf\left \lbrace t> 0,\, (X,\,N)_t \in \cT_0 \right \rbrace,
	\\&
	\Tex:= \tDo \wedge T^N_\infty \wedge \tau_0 \wedge \tau_{E} \wedge \ext.
	\EQn{U:Uout}{\Tex}
\end{align*}
Since, on the event $\left \lbrace \Tex = \tDo \right \rbrace$, 
either $N_{\tDo} \ge n_{\infty}$
or $(X,Y)_{\tDo} \in \cT^X_\infty$:
\begin{align*}
	&\E_{(x, n)}[\exp(\rho \, \widehat{\tau}_{E} )]
	= \E_{(x, n)}[\exp(\rho \, \Tex ) \,;\; \Tex = \widehat{\tau}_{E}]
	+\E_{(x, n)}[\exp(\rho \, \widehat{\tau}_{E} )\,;\; \Tex = \tau_0]
	\\& \hspace{3cm}
	+\E_{(x, n)}[\exp(\rho \, \widehat{\tau}_{E} ) \,;\; \Tex = \tDo]
	+\E_{(x, n)}[\exp(\rho \, \widehat{\tau}_{E} ) \,;\; \Tex = T^N_\infty]
	\\& \hspace{1cm}
	\le  \exp(\rho \, \tDo ) \, \left( 1+ \cE_0  \right)
	+ \exp(\rho \, \tDo ) \, \PR_{(x, n)}[\Tex = \tDo] \cE^X_\infty
	\\& \hspace{1.2cm}
	+ \exp(\rho \, \tDo )\; \left( \PR_{(x, n)}[\Tex = T^N_\infty]
	+ \PR_{(x, n)}[N_{\tDo} \ge n_{\infty},\, \Tex = \tDo]\right)
	\, \cE^N_\infty,
\end{align*}
by the Markov property. 
Now, by \Req{U:Uout}, 
$N^D$ is an upper-bound of $N$ before $\Tex$.
Thus,  by our definitions of $\tDo, n\iET, r_\vee$:
\begin{align*}
	&\E_{(x, n)}[\exp(\rho \, \widehat{\tau}_{E} )]
	\le  2 \,  \left( 1+ \cE_0  \right)
	+ \cE^X_\infty / 2
	+ \epsilon
	\, \cE^N_\infty / 2.
\end{align*}
Taking the supremum over $(x, n)\in \cT^X_\infty$ 
in the last inequality yields:

$\quad \cE^X_\infty
\le  4 \,  \left( 1+ \cE_0  \right)
+ \epsilon\, \cE^N_\infty.
$\hfill $\square$

\paragraph{Proof of Lemma \ref{U:N^DP}}
\textcolor{white}{:}

We recall our definition of $Z$
and $\psi$ in \Req{U:Z2}.
In the following, we condider $\psi$
as a function of $r$,
thus the notation 
$\psi_r(z):= - 1/(2\, z) +(r z)/2
- c\, z^3$.

\textbf{Step 1:} $\sup_{z>0} \, \psi_r(z) \underset{r\rightarrow -\infty}{\longrightarrow}-\infty$.

Let $A > 0$, $z_A:= \frac{2}{A}$ 
and $r_\vee:= -A^2$.
Then: 
\begin{align*}
	&\frl{z \le z_A} 
	\frlq{r \le 0}
	&&\psi_r(z) \le -1/(2\, z_A) = -A,&
	\\& \frl{z \ge  z_A} 
	\frlq{r \le r_\vee \le 0} 
	&&\psi_r(z) \le r_\vee \, z_A  = - 2 A.&
	\tag*{$\square$}
\end{align*}

\textbf{Step 2:} 
bound on $Z^{A}_t:= z - A\, t + B_t$ 
for $A$ large:

Let $\epsilon, t_D > 0$. 
We can choose $\Delta z > 0$ such that, 
with $N\sim \mathcal{N}\left(0, 1\right)$: 
\begin{align*}
	\PR\left(\Tsup{t\le t_D} B_t \ge \Delta z\right)
	= 2\, \PR\left( N \ge  \Delta z \, /\,\sqrt{t_D}\right) \le \epsilon. 
	\EQn{U:DeltaY}{\Delta z}
\end{align*}
Then, we can choose $A> 0$ (sufficiently big) such that:
$$\PR\left( B_{t_D} \ge  A\, t_D\right) 
= \PR\left( N \ge  A\,\sqrt{t_D}\right) \le~\epsilon. $$
We also choose $r_\vee$ thanks to step 1
such that: 
$$\frl{r \ge r_\vee}\quad \sup_{z>0} \, \psi_r(z)  \le -A\le 0.$$
We now assume that the initial condition 
of $Z$ satisfies $Z\le z_\infty$
($z_\infty = \sigma \sqrt{n_{\infty}} /2$).

For any 
$z\iET\ge z_{\infty} + \Delta z$ and
$r \le r_\vee,$
we deduce:
\begin{align*}
	&
	\Tsup{t\le t_D} Z_t 
	\le z_{\infty} + \Tsup{t\le t_D} B_t,
	\\&
	\PR\left( \Tsup{t\le t_D} Z_t  \ge z\iET \right)
	\le \PR\left(\Tsup{t\le t_D} B_t \ge \Delta z\right) 
	\le \epsilon
	\quad \text{ by } \Req{U:DeltaY},
	\\ &
	\PR\left( Z_{t_D}  \ge z_\infty \right)
	\le \PR\left( B_{t_D} \ge  A\, t_D\right) 
	\le \epsilon
	\quad\text{ by our choice of A and } r_\vee.
	\end{align*}
Thus 
	\begin{align*}
	\PR\left( T^D_{\infty}\le t_D  \right)
	\le \PR\left( Z_{t_D}  \ge z_{\infty} \right)
	+ \PR\left( \Tsup{t\le t_D} Z_t  \ge z\iET \right)
	\le 2\, \epsilon,
\end{align*}
with $n\iET = (2\, z\iET/\sigma)^2$.
It proves the second claim of the Lemma 
(up to a change of $\epsilon$ by $\epsilon/2$).
\\

\textbf{Step 3:} descent from infinity and extinction 

Now, we need to assume $\cY > 0$. Let again $\epsilon, t_D>0$.
Thanks to Lemma \ref{U:ydP} (for $r =0$ since $\PR(t_D < \tau_{\downarrow })$ is decreasing with $r$)
we choose $z_{\downarrow}>0$ such that, with\\
$\tau_{\downarrow } 
:=	\inf \left \lbrace t\ge 0; \; Z_t \le z_{\downarrow} \right \rbrace$:
\begin{align*}
	&\frl{r \le 0}
	\frl{ z> 0} \quad
	\PR_{z_{\infty}} (t_D < \tau_{\downarrow }) 
	\le \epsilon
	\EQn{U:Eyd}{z_{\downarrow}}
\end{align*}
Like in the previous step, we choose $A>0$ such that:
$$\PR\left( B_{t_D} \ge  A\, t_D - z_{\downarrow} \right) 
\le \epsilon. $$
Again, we choose $r_\vee$ thanks to step 1 such that:
$\quad 
\frl{r \le r_\vee}\quad
\sup_{z>0} \, \psi_r(z)  \le -A\le 0.$

Then, with $r \le r_\vee$, on the event $\left \lbrace \tau_{\downarrow } \le t_D \right \rbrace$, conditionally on 
$Z_{\tau_{\downarrow }}$:
\begin{align*}
	&\PR_{N^D_{\tau_{\downarrow }}} \left( 2\, t_D - \tau_{\downarrow } < \widetilde{\ext}^D \right)
	\le \PR_{Z_{\tau_{\downarrow }}} \left( \widetilde{Z}_{t_D} > 0\right)
	\le \PR\left( z_{\downarrow} - A\, t_D + B_{t_D} >0  \right)  \le \epsilon,
	\EQn{U:periclite}{Z_{\tau_{\downarrow }}}
\end{align*}
by our choices of $A$ and $r_\vee$. 
Finally, by the Markov property, for any $z> 0$:
\begin{align*}
	&\PR_{z_{\infty}} \left( 2\, t_D < \ext^D \right)
	\le  \PR_{z_{\infty}} \left( t_D < \tau_{\downarrow } \right) 
	+\E_{z_{\infty}} \left[ 
	\PR_{Z_{\tau_{\downarrow }}} \left(  2\, t_D - \tau_{\downarrow } < \widetilde{\ext^D} \right)
	;\;  \tau_{\downarrow } \le t_D \right]\\
	&\hspace{2cm} 
	\le 2\; \epsilon
	\quad \text{ with } \Req{U:Eyd},\, \Req{U:periclite}
\end{align*}
which proves the first claim of the Lemma
(replace $\epsilon$ by $\epsilon/2$ in the proof 
and take $t_D = t/2$).

\hfill $\square$

\subsubsection*{\textbf{Appendix D:} 
	Too few individuals,
	proof of Proposition \ref{U:cE0}:}
\label{U:sec:cE0}

For $(x, n)\in \cT_0$, with $n_0$ sufficiently small, we would like to say that mortality is so strong in this area that it overcomes an exponential growth at rate $\rho$.
In order to get an estimate of mortality in $\cT_0$, 
we will use some coupling with branching processes 
and consider the process after a time $t_D = 1$
(arbitrary).
In practice, we prove
that for any $\rho, \epsilon'>0$, 
there exists $C'\ge 1$
such that for any $n\iET$ sufficiently large:
\begin{equation*} 
	\cE_0 \le C'  + \epsilon' \, \left( \cE^N_\infty +\cE^X_\infty + \cE_0 \right).
\end{equation*}
By taking $\epsilon' = (\epsilon\wedge 1) /2$, $C = 2\,C'$, it clearly implies
Proposition \ref{U:cE0}.
\\

The equation 
$\quad 
\textstyle{N^{U}_t
	= n_0 +  \int_{0}^{t} r_+ N^{U}_s  ds  
	+ \sigma \int_{0}^{t} \sqrt{N^{U}_{s}} \; dB^N_{s}}
\quad $
defines an upper-bound of $N$ on $[0, t_D]$
provided $n\le n_0$,
while $N^{U}$ is a classical branching process.
The survival of $(X, N)$ beyond $t_D$ clearly implies 
the survival of $N^{U}$ beyond $t_D$. 
Let us 
define $\rho_0$ by the relation:
$\quad \PR_{n_0}\left( t_D < \ext^{U} \right)  
=: \exp(- \rho_0\; t_D).\quad $
For a branching process like $N^{U}$, it is classical that:
$\rho_0\rightarrow \infty$ 
as $n_0 \rightarrow 0$.
Indeed, 
with $u(t, \lambda)$ the Laplace exponent of $N^U$
(cf e.g. \cite{P16} Subsection 4.2, notably Lemma 5):
$\quad 
\PR_{n_0}\left( \ext^{U} \le t_D  \right)  
=  \exp[ - n_0\, \lim_{\lambda \rightarrow \infty} u(t_D, \lambda)]
\rightarrow 1,$ as $n_0\rightarrow 0$.

So we can impose that $\rho_0 > \rho$, 
and even that $\exp(- (\rho_0-\rho)\; t_D)$ is sufficiently small 
to make transitions from $\cT_0$ 
to $\cT_0$, $\cT^N_\infty$ or $\cT^X_\infty$ of little incidence.
\begin{align*}
	\E_{(x, n)}[\exp(\rho \widehat{\tau}_{E} )]
	&\le 
	\E_{(x, n)}\Big[
	\exp(\rho \widehat{\tau}_{E} );\;
	\widehat{\tau}_{E} < t_D
	\Big]
	\\&\hspace{2.5cm}
	+	\E_{(x, n)}\Big[
	\exp(\rho \widehat{\tau}_{E} );
	(x, n)_{t_D} \in \cT_0\cup \cT^N_\infty \cup \cT^X_\infty \Big]
	\\&
	\le \exp[\rho\; t_D]
	+ \exp(\rho\; t_D)
	\; (\cE_0 + \cE^N_\infty + \cE^X_\infty)
	\; \PR_{(x, n)} (t_D < \ext) 
	\\ &
	\le C'
	+ \epsilon'
	\; (\cE_0 + \cE^N_\infty + \cE^X_\infty),
	\qquad 
	\text{  where } C':= \exp[\rho\; t_D]
	\\
	\text{ and } \epsilon':= &\exp(- (\rho_0-\rho)\; t_D)\rightarrow 0
	\text{ as } n_0\rightarrow 0.
	\tag*{$\square$}
\end{align*}

\subsubsection*{\textbf{Appendix E:} 
	Proof of Fact \ref{U:ser}}
Like in the proof of Lemma \ref{U:dec} with Fact \ref{U:ser}:
\begin{align*}
	&\mu A_{t_h } =  \dfrac{\PR_\mu (j\, t\iDB< \ext)}{\PR_\mu (t_h < \ext)}
	\; \mu A_{j\, t\iDB} \cdot \PR_{t_h -j\, t\iDB}
	\\&\hspace{.7cm}
	= \sum_{k = 1}^j \amu(k, j\, t\iDB)\times 
	\dfrac{\PR_\mu (j\, t\iDB< \ext)}{\PR_\mu (t_h < \ext)} 
	\times \dfrac{\PR_{\zeta} (t_h -k\, t\iDB< \ext)}
	{\PR_{\zeta} ( [j -k]\, t\iDB< \ext)}
	\zeta A_{t_h -k\, t\iDB}
	\\ &\hspace{1.3cm}
	+ r_j
	\times \dfrac{\PR_\mu (j\, t\iDB< \ext)}{\PR_\mu (t_h < \ext)} 
	\times \PR_{\nu_j} (t_h -j\, t\iDB< \ext) \;
	\nu_j A_{t_h -j\, t\iDB}
	\EQn{muAth}{\nu_j}
\end{align*}
Yet, by \Req{U:AMU}:
$\quad \amu(k, j\, t\iDB)\times 
\dfrac{\PR_\mu (j\, t\iDB< \ext)}{\PR_\mu (t_h < \ext)} 
\times \dfrac{\PR_{\zeta} (t_h -k\, t\iDB< \ext)}
{\PR_{\zeta} ( [j -k]\, t\iDB< \ext)}
=  \dfrac{c\iDB}{c\iPs}\, \left( 1-\dfrac{c\iDB}{c\iPs} \right)^{k-1},$

\noindent
so that we obtain, by evaluating the measures in \Req{muAth}\, on $\cX$:
\begin{align*}
	&r_j
	\times \dfrac{\PR_\mu (j\, t\iDB< \ext)}{\PR_\mu (t_h < \ext)} 
	\times \PR_{\nu_j} (t_h -j\, t\iDB< \ext)
	= 1 -  \sum_{k = 1}^{j}   \dfrac{c\iDB}{c\iPs}\, \left( 1-\dfrac{c\iDB}{c\iPs} \right)^{k-1}
	= \left( 1-\text{\large{$\,^{c\iDB}\!/_{c\iPs}$} } \right)^{j}.
	\tag*{$\square$}
\end{align*}
\section*{Aknowledgment}
I am very grateful to Etienne Pardoux, my PhD supervisor, 
for his great support all along the redaction of this article.
The comments of my reviewer 
have also been 
particularly detailed and helpful. I express my sincere thanks to him.
I wish to thank also Nicolas Champagnat
for their advices and support.

I would like finally to thank the very inspiring meetings and discussions brought about by the Chair “Modélisation Mathématique et Biodiversité” of VEOLIA-Ecole Polytechnique-MnHn-FX.


\end{document}